\documentclass[11pt]{article}
\usepackage[greek,english]{babel}
\usepackage[utf8]{inputenc}
\usepackage[T1]{fontenc}
\usepackage{lmodern}
\usepackage{geometry}
\usepackage{ragged2e}
\usepackage{pgfplots}
\usepackage{amssymb,amsmath,latexsym,color,textcomp,marvosym}
\usepackage{ulem,url}
\usepackage{multirow}
\usepackage[counterclockwise]{rotating}
\usepackage{dsfont}
\usepackage{mathrsfs}
\usepackage{comment}
\usepackage{tabulary}
\usepackage{mathtools}

\newcommand{\bveps}{\boldsymbol{\varepsilon}}

\newcommand{\bd}{\boldsymbol{d}}
\newcommand{\bb}{\boldsymbol{b}}

\newcommand{\bc}{\boldsymbol{c}}

\newcommand{\bA}{\boldsymbol{A}}

\newcommand{\bD}{\boldsymbol{D}}
\newcommand{\bx}{\boldsymbol{x}}

\newcommand{\bu}{\boldsymbol{u}}
\newcommand{\bv}{\boldsymbol{v}}
\newcommand{\bX}{\boldsymbol{X}}
\newcommand{\bXt}{\boldsymbol{X_t}}

\newcommand{\bS}{\boldsymbol{S}}
\newcommand{\bs}{\boldsymbol{s}}

\newcommand{\be}{\boldsymbol{e}}

\newcommand{\bmu}{\boldsymbol{\mu}}
\newcommand{\bm}{\boldsymbol{m}}

\newcommand{\scal}[1]{\langle #1 \rangle}

\usepackage{etoolbox}
\patchcmd{\thebibliography}{\section*{\refname}}{}{}{}

\renewcommand{\thefootnote}{\raise.8ex\hbox{\scriptsize{\arabic{footnote}}}}
\renewcommand{\baselinestretch}{1.5}\small\normalsize

\newcommand{\sign}{\text{sign}}
\newcommand{\tg}{\text{tg}}

\def\cqfd{\null\hfill{$\Box$}\par\vspace*{0.2cm}}

\def\1g{1\hskip -3pt \mbox{l}}

\newtheorem{lem}{Lemma}

\newtheorem{prop}{Proposition}
\newtheorem{cor}{Corollary}
\newtheorem{rem}{Remark}

\newtheorem{exemple}{Example}
\newtheorem{definition}{Definition}

\numberwithin{equation}{section}
\numberwithin{theo}{section}
\numberwithin{lem}{section}
\numberwithin{cor}{section}
\numberwithin{exemple}{section}
\numberwithin{rem}{section}
\numberwithin{property}{section}
\numberwithin{prop}{section}
\numberwithin{definition}{section}

\begin{document}

	\renewcommand{\baselinestretch}{1.5}
	\selectlanguage{english}
	\title{{\LARGE Path prediction of aggregated $\alpha$-stable moving averages\\ using semi-norm representations}}
	\author{
		{\sc Sébastien Fries\footnote{ENSAE-CREST, 5 Avenue Henri Le Chatelier, 91120 Palaiseau, France. E-Mail:
				sebastien.fries@gmail.com \hspace*{1cm} I am grateful to the Agence Nationale de la Recherche (ANR), which supported this work via the Project
				MultiRisk  (ANR CE26 2016 - CR)}
	}}
	
	\date{}
	
	\maketitle

	\renewcommand{\baselinestretch}{1}
	
	\renewcommand{\baselinestretch}{1.5}

	\begin{center}
		
		\vspace*{-0.75cm}
		
		\Large{09/2018}
		
		\vspace*{0.5cm}
		
		\textbf{Abstract}
		
	\end{center}

	{\justify
		For $(X_t)$ a two-sided $\alpha$-stable moving average, this paper studies the conditional distribution of future paths  given a piece of observed trajectory when the process is far from its central values.
		Under this framework, vectors of the form $\bXt=(X_{t-m},\ldots,X_t,X_{t+1},\ldots,X_{t+h})$, $m\ge0$, $h\ge1$, are multivariate $\alpha$-stable and the dependence between the past and future components is encoded in their spectral measures.
		A new representation of stable random vectors on unit cylinders --sets $\{\bs\in\mathbb{R}^{m+h+1}: \hspace{0.3cm} \|\bs\|=1\}$ for $\|\cdot\|$ an adequate semi-norm-- is proposed in order to describe the tail behaviour of vectors $\bXt$ when only the first $m+1$ components are assumed to be observed and large in norm. 
		Not all stable vectors admit such a representation and $(X_t)$ will have to be <<anticipative enough>> for $\bXt$ to admit one. 
		The conditional distribution of future paths can then be explicitly derived using the regularly varying tails property of stable vectors and has a natural interpretation in terms of pattern identification.
		The approach extends to processes resulting from the linear combination of stable moving averages and applied to several examples.
		
	}
	
	\vspace*{1cm}
	
	\noindent \textit{Keywords:} Anticipative processes, Noncausal processes, Stable processes, Stable random vectors,\\ 
	\hspace*{1.8cm} Spectral representation, Pattern identification, Prediction \\
	
	\vspace*{0.1cm}
	
	\noindent \textit{MSC classes:}  	60G52, 60E07, 60G25
	
	\restoregeometry
	\renewcommand{\baselinestretch}{1.5}

	\section{Introduction}
	
	Stochastic processes depending on the future values of an i.i.d. sequence, often referred to as \textit{anticipative}, have witnessed a recent surge of attention from the statistical and econometric literatures.
	This gain of interest is driven in particular by their convenience for modelling exotic patterns in time series, such as explosive bubbles in financial prices \cite{cav18,cht17,fri17,f18,gjm16,gz17,hec16,hec17a,hec17b, hen15} (see also \cite{and09,beh11a,beh11b,che17,gj16,gj18,lan11,lan13,ss16}).
	The attractive flexibility of anticipative processes cannot yet be fully leveraged however, as their dynamics, and especially the conditional distribution of future paths given the observed past trajectory, remains largely mysterious. 
	A remarkable exception is that of the anticipative $\alpha$-stable AR(1) for which partial results were obtained in \cite{gz17} and further completed in \cite{f18}. 
	Even in this simplest case within the family of anticipative processes however, future realisations feature a complex dependence on the observed past, which is reflected in the functional forms of the conditional moments obtained in \cite{f18}.
	Interestingly, the dynamics simplifies when the anticipative stable AR(1) is far from its central values, where it appears to follow an explosive exponential path with a determined killing probability. 
	This naturally raises the question of whether and under which form such a behaviour could be found in more general anticipative linear processes.
	
	For $X_t=\sum_{k\in\mathbb{Z}}d_k\varepsilon_{t+k}$ a two-sided moving average with $(\varepsilon_t)$ an independent and identically distributed (i.i.d.) $\alpha$-stable sequence, this paper analyses the conditional distribution of future paths given the observed trajectory, say $(X_{t+1},\ldots,X_{t+h})$ given $(X_{t-m},\ldots,X_t)$, $m\ge0$, $h\ge1$, when the process is far from its central values.  
	Under this framework, any vector of the form $\bXt=(X_{t-m},\ldots,X_{t+h})$ is multivariate $\alpha$-stable and its distribution is characterised by a unique finite measure $\Gamma$ on the Euclidean unit sphere $S_{m+h+1}=\{\bs\in\mathbb{R}^{m+h+1}:\hspace{0.2cm}\|\bs\|_e=1\}$, where $\|\cdot\|_e$ denotes the Euclidean norm (Theorem 2.3.1 in \cite{st94}).  
	The measure $\Gamma$ in particular completely describes 
	the conditional distribution of the normalised paths $\bXt/\|\bXt\|_e$, the <<shape>> of the trajectory, when $\bXt$ is large according to the Euclidean norm and given some information about the observed first $m+1$ components. 
	A straightforward application of Theorem 4.4.8 by Samorodnitsky and Taqqu (1994) \cite{st94} indeed shows that
	\begin{align}
		\label{eqintro:regvar}
		\mathbb{P}\Big(\boldsymbol{X}_t/\|\boldsymbol{X}_t\|_e \in A \hspace{0.2cm} \Big| \hspace{0.2cm} \|\boldsymbol{X}_t\|_e>x \hspace{0.2cm} \text{and} \hspace{0.2cm} \boldsymbol{X}_t/\|\boldsymbol{X}_t\|_e\in B\Big) \underset{x\longrightarrow\infty}{\longrightarrow} \dfrac{\Gamma(A\cap B)}{\Gamma(B)},
	\end{align}
	for any appropriately chosen Borel sets $A,B\subset S_{m+h+1}$.
	As such however, \eqref{eqintro:regvar} is of little value for prediction purposes where only $X_{t-m},\ldots,X_t$ are assumed to be observed, given that the conditioning generally depends on the future realisations $X_{t+1},\ldots,X_{t+h}$. 
	The idea developed here is to obtain a version of \eqref{eqintro:regvar} where the Euclidean norm is replaced by a semi-norm $\|\cdot\|$ satisfying
	\begin{align}\label{asu:semin_intro}
		\|(x_{-m},\ldots,x_0,x_1,\ldots,x_h)\|=\|(x_{-m},\ldots,x_0,0,\ldots,0)\|,
	\end{align}
	for any $(x_{-m},\ldots,x_h)\in\mathbb{R}^{m+h+1}$. In this view, a new representation of stable random vectors on the <<unit cylinder>> $C_{m+h+1}^{\|\cdot\|}:=\{\bs\in\mathbb{R}^{m+h+1}: \hspace{0.2cm} \|\bs\|=1\}$ is thus explored, where $\|\cdot\|$ is such a semi-norm.
	Contrary to representations involving norms (see Theorem 2.3.8 in \cite{st94}), not all stable random vectors admit representations on unit cylinders and a characterisation is provided. 
	It is shown that only if $(X_t)$ is <<anticipative enough>> will $\bXt$ admit a representation by a measure $\Gamma^{\|\cdot\|}$ on $C_{m+h+1}^{\|\cdot\|}$. 
	Property \eqref{eqintro:regvar} is then shown to hold with an adequate semi-norm and with $\Gamma$ (resp. $S_{m+h+1}$) replaced by $\Gamma^{\|\cdot\|}$ (resp. $C_{m+h+1}^{\|\cdot\|}$). The problem finally boils down to choosing appropriate 
	Borels $B$ in \eqref{eqintro:regvar} reflecting that only the past <<shape>> $(X_{t-m},\ldots,X_{t})/\|\bXt\|$ is observed.
	
	The use of \eqref{eqintro:regvar} to infer about the future paths of $(X_t)$ has connections with the so-called \textit{spectral process} introduced by Basrak and Segers (2009) \cite{bs09} which has opened a fruitful line of research (see for instance \cite{bps16,dhs17,j17,js14,ms10,ps17}). 
	This spectral process is defined as the limit in distribution of a vector of observations of a multivariate regularly varying time series conditionnally on the first observation being large.
	The approach followed here differs in that it operates at the representation level of $\alpha$-stable vectors, establishing a link between the spectral representation and the tail conditional distribution of stable linear processes and shedding light on the (un)predictability of their extremes.
	A natural interpretation of path prediction in terms of pattern identification emerges from Property \eqref{eqintro:regvar} applied to stable linear processes, similar to what Janssen (2017) \cite{j17} pointed out in a framework close to that of Basrak and Segers \cite{bs09}.
	The results are extended to encompass processes resulting from the linear combination of $\alpha$-stable moving averages, coined \textit{stable aggregates}, and illustrated on several examples. 
	Contrary to non-aggregated moving averages, which trajectories recurrently feature the same pattern from one extreme episode to another, stable aggregates appear flexible enough to accomodate trajectories exhibiting various patterns through time.
	
	Section \ref{sec:stabgen} characterises the representation of general $\alpha$-stable vectors on semi-norm unit cylinders and shows that Property \eqref{eqintro:regvar} can be restated under this new representation. 
	Focusing first on $\alpha$-stable moving averages and then on linear combination thereof, Section \ref{sec:MA_representability} studies under which condition on the process $(X_t)$ the vector $(X_{t-m},\ldots,X_{t+h})$ admits a representation on the unit cylinder $C_{m+h+1}^{\|\cdot\|}$. 
	The anticipativeness of $(X_t)$ surprisingly arises as a necessary condition for such a representation to exist.
	Section \ref{sec:application} then exploits Property \eqref{eqintro:regvar} to analyse the tail conditional distribution of general stable aggregates and of some particular processes: the aggregation of anticipative AR(1), the anticipative AR(2) and the anticipative fractionally integrated process. 
	Section \ref{sec:bivariate_AR1} finally considers a simple bivariate process to illustrate an extension to vector moving averages. 
	New properties emerge in higher dimensions where, in particular, the presence of a non-anticipative component does not rule out the existence of adequate semi-norm representations.
	Section \ref{sec:conclude} concludes and provides perspectives for future work.
	Proofs are collected in Section \ref{sec:proofs}.

	\section{Stable random vectors representation on unit cylinders}
	\label{sec:stabgen}
	
	This section starts by recalling the characterisation of stable random vectors on the Euclidean unit sphere before exploring the case of unit cylinders relative to semi-norms and reformulating the regularly varying tails property.
	
	\begin{definition}
		A random vector $\bX=(X_1,\ldots,X_d)$ is said to be a stable random vector in $\mathbb{R}^d$ if and only if for any positive numbers $A$ and $B$ there is a positive number $C$ and a non-random vector $\bD\in\mathbb{R}^d$ such that
		$$
		A\bX^{(1)}+B\bX^{(2)}\stackrel{d}{=}C\bX+\bD,
		$$
		where $\bX^{(1)}$ and $\bX^{(2)}$ are independent copies of $\bX$. Moreover, if $\bX$ is stable, then there exists a constant $\alpha\in(0,2]$ such that the above holds with $C=(A^\alpha+B^\alpha)^{1/\alpha}$, and $\bX$ is then called $\alpha$-stable.
	\end{definition}
	The Gaussian case ($\alpha=2$) is henceforth excluded. 
	For $0<\alpha<2$, the vector $\bX=(X_1,\ldots,X_d)$ is an $\alpha$-stable random vector if and only if there exists a unique pair $(\Gamma,\bmu^0)$, $\Gamma$ a finite measure on $S_d$ and $\bmu^0$ a non-random vector in $\mathbb{R}^d$, such that, 
	\begin{align}\label{def:char_fun_vec}
		\mathbb{E}\Big[e^{i\scal{\bu,\bX}}\Big] = \exp\Bigg\{-\int_{S_d}|\scal{\bu,\bs}|^\alpha \bigg(1-i\,\sign(\scal{\bu,\bs})w(\alpha,\scal{\bu,\bs})\bigg)\Gamma(d\bs)+i\,\scal{\bu,\bmu^0}\Bigg\}, \hspace{0.3cm} \forall\bu\in\mathbb{R}^d,
	\end{align}
	where $\scal{\cdot,\cdot}$ denotes the canonical scalar product, $w(\alpha,s) = \tg\left(\frac{\pi \alpha}{2}\right)$, if $\alpha \ne 1$, and $w(1,s) = - \frac{2}{\pi}\ln|s|$ otherwise, for $s\in\mathbb{R}$.
	The pair $(\Gamma,\bmu^0)$ is called the \textit{spectral representation} of the stable vector $\bX$, $\Gamma$ is its \textit{spectral measure} and $\bmu^0$ its \textit{shift vector}. 
	In particular, $\bX$ is symmetric if and only if  $\bmu^0=0$ and $\Gamma(A)=\Gamma(-A)$ for any Borel set $A$ in $S_d$ (Theorem 2.4.3 in \cite{st94}), and in that case
	\begin{align}\label{def:char_fun_vec_sas}
	\mathbb{E}\Big[e^{i\scal{\bu,\bX}}\Big] = \exp\Bigg\{-\int_{S_d}|\scal{\bu,\bs}|^\alpha\Gamma(d\bs)\Bigg\}, \hspace{0.3cm} \forall\bu\in\mathbb{R}^d.
	\end{align}
	In the univariate case, \eqref{def:char_fun_vec} boils down to
	\begin{align*}
		\mathbb{E}\Big[e^{i u X}\Big] = \exp\Bigg\{-\sigma^{\alpha}|u|^{\alpha}\bigg(1-i\beta \,\sign(u)w(\alpha,u)\bigg)+iu\mu\Bigg\}, \hspace{0.3cm} \forall u\in\mathbb{R},
	\end{align*}
	for some $\sigma>0$, $\beta\in[-1,1]$ and $\mu\in\mathbb{R}$. The representations \eqref{def:char_fun_vec} and \eqref{def:char_fun_vec_sas} of a stable random vector involves integration over all directions of $\mathbb{R}^d$,\footnote{\label{foot:direc} By \textit{direction} of $\mathbb{R}^d$, it is meant the equivalence classes of the relation <<$\equiv$>> defined by: $\bu\equiv\bv$ if and only if there exists $\lambda>0$ such that $\bu=\lambda\bv$, for $\bu,\bv\in\mathbb{R}^d$.} here parameterised by the unit sphere relative to the Euclidean norm. 
	Proposition 2.3.8 in \cite{st94} shows that the unit sphere relative to any norm can be used instead, provided a change of spectral measure and shift vector. 
	We study alternative representations where integration is performed over a unit cylinder relative to a semi-norm. 
	For a given semi-norm, not all stable vectors admit such a representation, which motivates the following definition. 
	
	\begin{definition}\rm
		Let $\|\cdot\|$ be a seminorm on $\mathbb{R}^d$, $C_d^{\|\cdot\|}:=\{\bs\in\mathbb{R}^d:\hspace{0.2cm} \|\bs\|=1\}$ be the corresponding {unit cylinder}, and let $\bX=(X_1,\ldots,X_d)$ be an $\alpha$-stable random vector. \\
		\textit{(Asymmetric case)} In the case where $\bX$ is not symmetric, we say that $\bX$ \textit{is representable on} $C_d^{\|\cdot\|}$ if there exists a non-random vector $\bmu_{\|\cdot\|}^0\in\mathbb{R}^d$ and a Borel measure $\Gamma^{\|\cdot\|}$ on $C_d^{\|\cdot\|}$ satisfying for all  $\bu\in\mathbb{R}^d$
		\begin{align}
    	\int_{C_d^{\|\cdot\|}}|\scal{\bu,\bs}|^\alpha \Gamma^{\|\cdot\|}(d\bs) & < +\infty,
    	\label{asu:integ_simple}
    	\end{align}
    	if $\alpha\ne1$, and if $\alpha=1$,
    	\begin{align}
		\int_{C_d^{\|\cdot\|}}|\scal{\bu,\bs}|\Big|\ln|\scal{\bu,\bs}|\hspace{0.03cm}\Big| \Gamma^{\|\cdot\|}(d\bs) & < +\infty,
		\label{asu:integ_complique}
		\end{align}
		such that the joint characteristic function of $\bX$ can be written as in \eqref{def:char_fun_vec} with $(S_d,\Gamma,\bmu^0)$ replaced by $(C_d^{\|\cdot\|},\Gamma^{\|\cdot\|},\bmu^0_{\|\cdot\|})$.\\
		\textit{(Symmetric case)} In the case where $\bX$ is symmetric $\alpha$-stable (S$\alpha$S), $0<\alpha<2$, we say that $\bX$ {is representable on} $C_d^{\|\cdot\|}$ if there exists a symmetric Borel measure $\Gamma^{\|\cdot\|}$ on $C_d^{\|\cdot\|}$ satisfying \eqref{asu:integ_simple}
		such that the joint characteristic function of $\bX$ can be written as in \eqref{def:char_fun_vec_sas} with $(S_d,\Gamma)$ replaced by $(C_d^{\|\cdot\|},\Gamma^{\|\cdot\|})$.
	\end{definition}
	
	\begin{rem}
	\rm As unit cylinders are unbounded sets, the integrability conditions \eqref{asu:integ_simple}-\eqref{asu:integ_complique} 
	ensure the sanity of the above definition.
	\end{rem}
	
	\noindent We start by characterising stable random vectors that are representable on a given semi-norm unit cylinder.
	
	\begin{prop}\label{prop238}
		Let $\|\cdot\|$ be a seminorm on $\mathbb{R}^d$ and $C_d^{\|\cdot\|}$ be the corresponding unit cylinder. 
		Denote $K^{\|\cdot\|}=\{x\in S_d:\|x\|=0\}$. 
		Let also $\bX$ be an $\alpha$-stable random vector on $\mathbb{R}^d$ with spectral representation $(\Gamma,\bmu^0)$ on the Euclidean unit sphere (with $\bmu^0=0$ if $\bX$ is S$\alpha$S). 
		If $\alpha\ne1$ or if $\bX$ is S1S, then
		$$
		\bX \text{ is representable on } C_d^{\|\cdot\|} \iff \Gamma( K^{\|\cdot\|})=0.
		$$
		If $\alpha=1$ and $\bX$ is not symmetric, then
		\begin{align*}
			& \bX \text{ is representable on } C_d^{\|\cdot\|} \iff \int_{S_d}\Big|\ln\|\bs\|\hspace{0.03cm}\Big|\Gamma(d\bs) < +\infty.
		\end{align*}
		Moreover, if $\bX$ is representable on $C_d^{\|\cdot\|}$, its spectral representation is then given by $(\Gamma^{\|\cdot\|},\bmu^0_{\|\cdot\|})$ where
		$$
		\Gamma^{\|\cdot\|} (d\bs) = \|\bs\|_e^{-\alpha}\, \Gamma\,\circ\,T_{\|\cdot\|}^{-1}(d\bs)
		$$
		with $T_{\|\cdot\|}:S_d\setminus K^{\|\cdot\|}\longrightarrow C_d^{\|\cdot\|}$ defined by $T_{\|\cdot\|}(\bs)=\bs/\|\bs\|$, and
		$$
		\bmu^0_{\|\cdot\|} = \left\{
		\begin{array}{ll}
		\bmu^0, & \text{if} \quad \alpha \ne 1 \hspace{0.45cm} \text{or} \hspace{0.45cm} \text{if} \hspace{0.3cm} \bX \hspace{0.3cm} \text{is} \hspace{0.3cm} S1S,\\
		\bmu^0 + \tilde{\bmu}, & \text{if} \quad \alpha = 1 \hspace{0.3cm} \text{and} \hspace{0.3cm} \bX \hspace{0.3cm} \text{is not symmetric},
		\end{array}
		\right.
		$$
		$$
		\tilde{\bmu} =(\tilde{\mu}_j), \hspace{1cm} \text{and} \hspace{1cm} \tilde{\mu}_j = -\dfrac{2}{\pi}\int_{S_d\setminus K^{\|\cdot\|}}s_j\ln\|\bs\|\Gamma(d\bs), \quad j=1,\ldots,d.
		$$
	\end{prop}
	
	\begin{rem}\rm
	The representability condition in the case [$\alpha=1$ and $\bX$ not symmetric] is slightly stronger than that in the other cases.
	Indeed, $\int_{K^{\|\cdot\|}}\Big|\ln\|\bs\|\hspace{0.03cm}\Big|\Gamma(d\bs)\le \int_{S_d}\Big|\ln\|\bs\|\hspace{0.03cm}\Big|\Gamma(d\bs) <+\infty$ necessarily implies that $\Gamma(K^{\|\cdot\|})=0$ since $\Big|\ln\|\bs\|\hspace{0.03cm}\Big|=+\infty$ for $\bs\in K^{\|\cdot\|}$.
	\end{rem}
	
	\begin{rem}\label{rem:interp_repres}
		\rm
		The case $d=2$ is insightful. 
		In view of \eqref{eqintro:regvar}, the spectral measure of the $\alpha$-stable vector $(X_1,X_2)$ describes its likelihood of being in any particular direction
		of $\mathbb{R}^2$ when it is large in norm.
		As unit spheres relative to norms span all the directions of $\mathbb{R}^2$, spectral measures on such spheres can describe any potential tail dependence of $(X_1,X_2)$. 
		Unit cylinders however do not span all directions of $\mathbb{R}^2$ and spectral measures thereon necessarily encode less information. 
		Consider for instance the unit cylinder $C_2^{\|\cdot\|}=\{(s_1,s_2)\in\mathbb{R}^2: \hspace{0.2cm} |s_1|=1\}$ associated to the semi-norm such that $\|(x_1,x_2)\|=|x_1|$ for all $(x_1,x_2)\in\mathbb{R}^2$. It is easy to see that $C_2^{\|\cdot\|}$ spans all directions of $\mathbb{R}^2$ but the ones of $(0,-1)$ and $(0,+1)$. 
		A stable vector $(X_1,X_2)$ will admit a representation on $C_2^{\|\cdot\|}$ provided these directions are irrelevant to characterise its distribution, that is, if $\Gamma\Big(\{(0,-1),(0,+1)\}\Big)=0$. 
		In terms of tail dependence, the latter condition intuitively means that realisations $(X_1,X_2)$ where $X_2$ is extreme and $X_1$ is not almost never occur (i.e., occur with probability zero).\footnote{The conditions $\Gamma\big(\{(0,-1),(0,+1)\}\big)=0$ and $\int_{S_2}\big|\ln\|\bs\|\hspace{0.03cm}\big|\Gamma(d\bs) < +\infty$ can also be related to the stronger condition ensuring the existence of conditional moments of $X_2$ given $X_1$ obtained in \cite{ct94,ct98} (see also Theorem 5.1.3 in \cite{st94}) and which requires $\Gamma$ not to be too concentrated around the points $(0,\pm1)$. 
		Namely, assuming $\int_{S_2} |s_1|^{-\nu}\Gamma(ds)<+\infty$ for some $\nu\ge0$, then $\mathbb{E}[|X_2|^\gamma|X_1]<+\infty$ for $\gamma<\min(\alpha+\nu,2\alpha+1)$, despite the fact that $\mathbb{E}[|X_2|^\alpha]=+\infty$.
		If the previous holds for some $\nu>0$, then necessarily both of the aforementioned conditions are satisfied.} 
	\end{rem}
	Provided the adequate representation exists, Property \eqref{eqintro:regvar} then holds with semi-norms instead of norms, providing the cornerstone for studying the tail conditional distribution of stable processes.
	\begin{prop}\label{prop:cond_tail} 
		Let $\bX=(X_1,\ldots,X_d)$ be an $\alpha$-stable random vector and let $\|\cdot\|$ be a seminorm on $\mathbb{R}^d$. If $\bX$ is representable on $C_d^{\|\cdot\|}$, then for every Borel sets $A,B\subset C_d^{\|\cdot\|}$ with $\Gamma^{\|\cdot\|}\Big(\partial (A\cap B)\Big)=\Gamma^{\|\cdot\|}\big(\partial B\big)=0$, and $\Gamma^{\|\cdot\|}(B)>0$,
		\begin{align}\label{eq:condtail}
			\mathbb{P}_x^{\|\cdot\|}(\bX,A|B) & \underset{x\rightarrow+\infty}{\longrightarrow} \dfrac{\Gamma^{\|\cdot\|}(A\cap B)}{\Gamma^{\|\cdot\|}(B)},
		\end{align}
		where $\partial B$ (resp. $\partial (A\cap B)$) denotes the boundary of $B$ (resp. $A\cap B$), and
		$$
		\mathbb{P}_x^{\|\cdot\|}(\bX,A|B) := \mathbb{P}\bigg(\dfrac{\bX}{\|\bX\|}\in A\bigg| \|\bX\|>x, \dfrac{\bX}{\|\bX\|}\in B\bigg).
		$$
	\end{prop}
	
	\section{Unit cylinder representation for paths of stable linear processes}
	\label{sec:MA_representability}
	
	Given a semi-norm, Proposition \ref{prop:cond_tail} is only applicable to stable vectors that are representable on the corresponding unit cylinder. 
	This section investigates under which condition on an stable moving average $(X_t)$ vectors of the form $(X_{t-m},\ldots,X_t,X_{t+1},\ldots,X_{t+h})$ admit such representations. 
	A characterisation is proposed and is then extended to linear combination of stable moving averages. 
	Any semi-norm satisfying \eqref{asu:semin_intro} could be relevant for the prediction framework mentioned in introduction. However to fix ideas and avoid numerous cases with respect to all the possible kernels, 
	we restrict to semi-norms such that
	\begin{align}
		& \|(x_{-m},\ldots,x_0,x_1,\ldots,x_h)\|=0  \iff  x_{-m}=\ldots=x_0=0, \label{asu:semin}
	\end{align}
	for any $(x_{-m},\ldots,x_h)\in\mathbb{R}^{m+h+1}$,	which in particular satisfy \eqref{asu:semin_intro}.
	
	\begin{exemple}\rm
	Semi-norms on $\mathbb{R}^{m+h+1}$ satisfying \eqref{asu:semin} can be naturally obtained from norms on the $m+1$ first components of vectors.
	For any $p\in[1,+\infty]$, one can consider for instance semi-norms $\|\cdot\|$ defined by
	$$
	\|(x_{-m},\ldots,x_0,x_1,\ldots,x_h)\| = \Big(\sum_{i=-m}^0 |x_i|^p\Big)^{1/p},
	$$
	for any $(x_{-m},\ldots,x_0,x_1,\ldots,x_h)\in\mathbb{R}^{m+h+1}$ with by convention $\big(\sum_{i=-m}^0 |x_i|^p\big)^{1/p} = \sup\limits_{-m\le i \le 0} |x_i|$ for $p=+\infty$.
	
	\end{exemple}
	
	\subsection{The case of moving averages}
	
	Consider $(X_t)$ the $\alpha$-stable moving average defined by
	\begin{align}
		\label{def:ma}
		X_t   = \sum_{k\in\mathbb{Z}}d_k\varepsilon_{t+k}, \hspace{1cm} \varepsilon_t\stackrel{i.i.d.}{\sim}\mathcal{S}(\alpha,\beta,\sigma,0)
	\end{align}
	with $(d_{k})$ a real deterministic sequence such that 
	\begin{equation}
		\label{assump:coef}
	    \text{if} \hspace{0.3cm} \alpha\ne1 \hspace{0.3cm} \text{or} \hspace{0.3cm} (\alpha,\beta) = (1,0), \hspace{0.75cm}	0<\sum_{k\in\mathbb{Z}}|d_{k}|^s < +\infty, \hspace{0.75cm} \text{for some}\hspace{0.3cm} s\in(0,\alpha)\cap[0,1],
	\end{equation}
	and 
	\begin{align}\label{assump:coefa1}
	\text{if} \hspace{0.3cm} \alpha=1 \hspace{0.3cm} \text{and} \hspace{0.3cm} \beta\ne0,\hspace{1cm}
	0<\sum_{k\in\mathbb{Z}}|d_{k}| \Big|\ln|d_k|\hspace{0.05cm}\Big| < +\infty.
	\end{align}
	Letting for $m\ge0$, $h\ge1$,
	\begin{align}
		\label{def:multi}
		\bXt  = (X_{t-m},\ldots,X_t,X_{t+1},\ldots,X_{t+h}),
	\end{align}
	it follows from Proposition 13.3.1 in Brockwell and Davis (1991)
	that the infinite series converge almost surely and both $(X_t)$ and $\bX_t$ are well defined. 
	The random vector $\bXt$ is multivariate $\alpha$-stable: denoting $\bd_{k}:=(d_{k+m},\ldots,d_{k},d_{k-1},\ldots,d_{k-h})$ for $k\in\mathbb{Z}$, 
	the spectral representation of $\bXt$ on the Euclidean sphere reads $(\Gamma,\bmu^0)$ with
	\begin{align}\label{def:spectral}
		\Gamma & = \sigma^\alpha\sum_{\vartheta\in S_1}\sum_{k\in\mathbb{Z}} w_\vartheta \|\bd_{k}\|_e^\alpha \delta_{\left\{\dfrac{\vartheta\bd_{k}}{\|\bd_{k}\|_e}\right\}}, \\
		\bmu^0 & = - \mathds{1}_{\{\alpha=1\}}\frac{2}{\pi}\beta\sigma\sum_{k\in\mathbb{Z}}\bd_{k}\ln\|\bd_{k}\|_e,\nonumber 
	\end{align}
	where $w_\vartheta = (1+\vartheta\beta)/2$, $S_1=\{-1,+1\}$, $\delta$ is the dirac mass and by convention, if for some $k\in\mathbb{Z}$, $\bd_{k}=\boldsymbol{0}$, i.e. $\|\bd_k\|_e=0$, then the $k$th term vanishes from the sums.
	Notice in particular that for $\beta=0$, it holds that $w_{-1} = w_{+1} =1/2$, $\bmu^0=\boldsymbol{0}$, and both the measure $\Gamma$ and the random vector $\bX_t$ are symmetric.
	The next result characterises the representability of $\bXt$ on a unit cylinder for fixed $m$ and $h$. 
	
	\begin{lem}\label{le:ma_representable} 
		Let $\bXt$ satisfy \eqref{def:ma}-\eqref{def:multi}  
		and let $\|\cdot\|$ be a semi-norm on $\mathbb{R}^{m+h+1}$ satisfying \eqref{asu:semin}. 
		For $\alpha\ne1$ and $(\alpha,\beta)=(1,0)$, the vector $\bXt$ is representable on $C_{m+h+1}^{\|\cdot\|}$ if and only if
		\begin{align}\label{eq:ma_representable}
			\forall k\in\mathbb{Z}, \quad \Big[(d_{k+m},\ldots,d_k)=\boldsymbol{0} \quad \Longrightarrow \quad\forall \ell \le k-1, \quad d_{\ell}=0\Big].
		\end{align}
		For $\alpha=1$ and $\beta\ne0$, the vector $\bXt$ is representable on $C_{m+h+1}^{\|\cdot\|}$ if and only if in addition to \eqref{eq:ma_representable}, it holds that
        \begin{align}\label{eq:additionallourd}
        \sum_{k\in\mathbb{Z}}\|\bd_k\|_e \bigg|\ln\Big(\|\bd_k\|/\|\bd_k\|_e\Big)\bigg| < +\infty.
        \end{align}
	\end{lem}
	In the cases $\alpha\ne1$ and $(\alpha,\beta)=(1,0)$, the representability of $\bXt$ on a semi-norm unit cylinder depends on the number of observation $m+1$ but not on the prediction horizon $h$. 
	Moreover, it is easy to see that if \eqref{eq:ma_representable} is true for some $m\ge0$, it then holds for any $m'\ge m$. 
	The case $\alpha=1$, $\beta\ne0$ is more intricate, the roles of $m$ and $h$ in the validity of the additional requirement \eqref{eq:additionallourd} not being as clear-cut.\\
	\indent A key distinction appears between moving averages according to whether finite length paths admit semi-norm representations. 
	This distinction especially matters for the applicability of Proposition \ref{eq:condtail} when studying the conditional dynamics of a given process.
	The following definition thus introduces the notion of \textit{past-representability} of a stable moving average.
	\begin{definition}\label{def:ma_representable}\rm
		Let $(X_t)$ be an $\alpha$-stable moving average satisfying \eqref{def:ma}-\eqref{assump:coefa1}. 
		We say that the stable process $(X_t)$ is \textit{past-representable} if there exists at least one pair $(m,h)$, $m\ge0$, $h\ge1$, such that $\bXt=(X_{t-m},\ldots,X_t,X_{t+1},\ldots,X_{t+h})$ is representable on $C_{m+h+1}^{\|\cdot\|}$ for some semi-norm satisfying \eqref{asu:semin}. 
		For any such pair $(m,h)$, we will say that $(X_t)$ is $(m,h)$-past-representable.
	\end{definition}
	\begin{rem}
		\rm It can be noticed that if $\bXt=(X_{t-m},\ldots,X_t,X_{t+1},\ldots,X_{t+h})$ is representable on $C_{m+h+1}^{\|\cdot\|}$ for some semi-norm satisfying \eqref{asu:semin}, then it is representable on unit cylinders relative to any other semi-norms satisfying \eqref{asu:semin}. This holds because \eqref{asu:semin} ensures that all these semi-norms have the same kernel. The notion of past-representability can thus be defined independently of the particular choice of a semi-norm.\footnote{This will not be true in general under the weaker assumption \eqref{asu:semin_intro} and different notions of representability of a process could emerge depending on the kernels of the semi-norms.}
	\end{rem}
	The following proposition provides a characterisation of past-representability.
	\begin{prop}\label{prop:ma_pastrepres}
		Let $(X_t)$ be an $\alpha$-stable moving average satisfying \eqref{def:ma}-\eqref{assump:coefa1}.\\ 
		$(\iota)$ With the set $\mathcal{M} = \{m\ge1: \hspace{0.2cm} \exists k\in\mathbb{Z}, \hspace{0.2cm} d_{k+m}=\ldots=d_{k+1}=0, \hspace{0.2cm} d_{k}\ne0\}$, define
		\begin{equation}\label{def:m0}
			m_0 = \left\{
			\begin{array}{cc}\sup\,\mathcal{M}, &\text{if} \hspace{0.5cm} \mathcal{M} \ne \emptyset,\\
				0 , & \text{if} \hspace{0.5cm} \mathcal{M} = \emptyset.
			\end{array}
			\right.
		\end{equation}
		\indent $(a)$ For $\alpha\ne1$ and $(\alpha,\beta)=(1,0)$, the process $(X_t)$ is past-representable if and only if
		\begin{align}\label{eq:ma_pastrepres}
			m_0 < +\infty.
		\end{align}
		\indent Moreover, letting $m\ge 0$, $h\ge1$, the process $(X_t)$ is $(m,h)$-past-representable if and only if \indent \eqref{eq:ma_pastrepres} holds and $m\ge m_0$.\\ 
		
		\indent $(b)$ For $\alpha=1$ and $\beta\ne0$, the process $(X_t)$ is past-representable if and only if in addition to \indent \eqref{eq:ma_pastrepres}, there exist an $m\ge m_0$ and an $h\ge1$ such that \eqref{eq:additionallourd} holds.
		If such a pair $(m,h)$ exists, \indent $(X_t)$ is then $(m,h)$-past-representable.\\
		
		\noindent $(\iota\iota)$ Let $\|\cdot\|$ a seminorm satisfying \eqref{asu:semin} and assume that $(X_t)$ is $(m,h)$-past-representable for some $m\ge0$, $h\ge1$. 
		The spectral representation $(\Gamma^{\|\cdot\|},\bmu^{\|\cdot\|})$ of the vector $\bXt=(X_{t-m},\ldots,X_t,X_{t+1},\ldots,X_{t+h})$ on $C_{m+h+1}^{\|\cdot\|}$ is then given by \eqref{def:spectral} with the Euclidean norm $\|\cdot\|_e$ replaced by the semi-norm $\|\cdot\|$.
	\end{prop}
	\begin{rem}\label{foot:0past}\rm
	Note in particular that $m_0=0$ if and only if for some $k_0\in\mathbb{Z}\cup \{-\infty\}$, $d_k\ne0$ for all $k\ge k_0$ and $d_k=0$ for all $k<k_0$.
	\end{rem}
	\begin{rem}\rm
	Proposition \ref{prop:ma_pastrepres} shows that for an $\alpha$-stable moving average to be past-representable, sequences of consecutive zero values in the coefficients $(d_k)$ have to be either of finite lengths, or infinite to the left. 
	This surprisingly places the anticipativeness of a stable moving average as a necessary --and sufficient for $\alpha\ne1$ and $(\alpha,\beta)=(1,0)$-- condition for its past-representability.
	The less anticipative a moving average is, in the sense of the larger the gaps of zeros in its forward-looking side, then the higher $m$ has to be chosen so as to have the representability of $(X_{t-m},\ldots,X_t,X_{t+1},\ldots,X_{t+h})$ on the appropriate unit cylinder. 
	Purely non-anticipative moving averages are in particular immediately ruled out.
	\end{rem}
	\begin{cor}\label{cor:nonanticiAR}\rm 
		Let $(X_t)$ an $\alpha$-stable moving average satisfying \eqref{def:ma}-\eqref{assump:coefa1}. 
		If $(X_t)$ is purely non-anticipative, i.e., $d_k=0$ for all $k\ge1$, then $(X_t)$ is not past-representable.
	\end{cor}
	
	\begin{rem}\label{rem:dar1}\rm
		This fault line between anticipativeness and non-anticipativeness sheds light on the predictability of extreme events in linear processes. 
		Consider for illustration the two following $\alpha$-stable AR(1) processes defined as the stationary solutions of
		\begin{align}
			X_t & = \rho X_{t+1} + \varepsilon_t, & \forall t\in\mathbb{Z}, & \label{def:ar1_a}\\
			Y_t & = \rho Y_{t-1} + \eta_t, & \forall t\in\mathbb{Z}, & \label{def:ar1_n}
		\end{align}
		where $0<|\rho|<1$, and $(\varepsilon_t)$, $(\eta_t)$ are independent i.i.d. stable sequences. 
		While $(X_t)$ generates bubble-like trajectories --explosive exponential paths eventually followed by sharp returns to central values--, the trajectories of $(Y_t)$ feature sudden jumps followed by exponential decays. 
		In both processes, an extreme event stems from a large realisation of an underlying error $\varepsilon_\tau$ or $\eta_\tau$, at some time $\tau$. 
		On the one hand for the non-anticipative AR(1) \eqref{def:ar1_n}, a jump does not manifest any early visible sign before its date of occurrence as it is independent of the past trajectory. 
		Jumps in the trajectory of $(Y_t)$ are unpredictable and one only has information about their unconditional likelihood of occurrence.
		On the other hand for the anticipative AR(1) \eqref{def:ar1_a}, extremes do manifest early visible signs and are gradually reached as their occurrence dates approach. The past trajectory is informative about future extreme events, and in particular more informative than their plain unconditional likelihood of occurrence. 
		Building on the <<information encoding>> interpretation of spectral measures given in Remark \ref{rem:interp_repres}, the fact that $(X_t)$ (resp. $(Y_t)$) is past-representable (resp. not past-representable) can be seen as a consequence of the dependence (resp. independence) of future extreme events on past ones.
	\end{rem}
	The condition for past-representability simplifies for ARMA processes and is equivalent to the autoregressive polynomial having at least one root located inside the unit circle.
	\begin{cor}\label{cor:mar}
		\rm
		Let $(X_t)$ be the strictly stationary solution of
		\begin{equation*}
			\psi(F)\phi(B)X_t=\Theta(F)H(B)\varepsilon_t, \hspace{1cm} \varepsilon_t\stackrel{i.i.d.}{\sim}\mathcal{S}(\alpha,\beta,\sigma,0),
		\end{equation*}
		where $\psi$, $\phi$, $\Theta$, $H$ are polynomials of arbitrary finite degrees with roots located outside the unit disk and $F$ (resp. $B$) is the forward (resp. backward) operator: $FX_t:=X_{t+1}$ (resp. $BX_t:=X_{t-1}$). 
		We suppose furthermore that $\psi$ and $\Theta$ (resp. $\phi$ and $H$) have no common roots. 
		Then, for any $\alpha\in(0,2)$ and $\beta\in[-1,1]$, the following statements are equivalent:\\
		$(\iota)  \hspace{0.4cm}(X_t) \hspace{0.1cm} \text{is past-representable},$ \\
		$(\iota\iota) \hspace{0.3cm} \deg(\psi)\ge1,$\\
		$(\iota\iota\iota) \hspace{0.3cm} m_0 < +\infty,$\\
		with $m_0$ as in \eqref{def:m0}.
		Moreover, letting $m\ge0$, $h\ge1$, the process $(X_t)$ is $(m,h)$-past-representable if and only if $m\ge m_0$ with $m_0<+\infty$. 
	\end{cor}
	\begin{rem}\rm
	For ARMA processes, we can notice in particular that the discrepancy between the cases [$\alpha\ne1$ or $(\alpha,\beta)=(1,0)$] and [$\alpha=1$, $\beta\ne0$] vanishes.
	Also, only the roots of the AR polynomial matter for past-representability, the MA part having no role.
	\end{rem}
	
	\subsection{Aggregation of moving averages}
	
	As will be seen in the next section, stable moving averages of the form \eqref{def:ma} generate trajectories bound to feature the same pattern $t\mapsto c d_{\tau-t}$ (up to a scaling $c$ and a time shift $\tau$) recurrently through time. 
	This can be seen as a strong limitation when it comes to time series modelling as argued by Gouriéroux and Zakoian (2017) \cite{gz17} in the context of explosive bubbles. 
	They suggest to alleviate this restriction by considering processes resulting from the linear combination of different models. These aggregations feature richer dynamics but little results are available to describe them (see for instance \cite{f18} for the aggregation of stable anticipative AR(1)).
	Linear combinations of stable moving averages will fit naturally into our framework and the results will extend.
	\begin{definition}\rm
		\label{def:aggma}
		Let $(X_{1,t}),\ldots,(X_{J,t})$ be $J\ge1$ stable moving averages, each satisfying \eqref{def:ma}-\eqref{assump:coefa1}, 
		for some coefficients sequences $(d_{j,k})_k$ and mutually independent error sequences $\varepsilon_{j,t}\stackrel{i.i.d.}{\sim}\mathcal{S}(\alpha,\beta_j,1,0)$, $j=1,\ldots,J$. 
		Let also $(\pi_j)_{j=1,\ldots,J}$ be positive numbers and define $(X_t)$ as
		\begin{align*}
			X_t = \sum_{j=1}^J \pi_j X_{j,t}, \hspace{1cm} \text{for} \hspace{0.3cm} t\in\mathbb{Z}.
		\end{align*}
		We will call such process $(X_t)$ a stable \textit{aggregated moving average}, an \textit{aggregated process}, or simply, a \textit{stable aggregate}, and call $(X_{j,t})$, $j=1,\ldots,J$ the \textit{latent} moving averages of $(X_t)$.
	\end{definition}
	We provide the spectral representation of paths of the aggregated process $(X_t)$ on the Euclidean unit sphere in the next lemma.
	\begin{lem}\label{le:spec_repres_agg}
	Let $(X_t)$ an $\alpha$-stable aggregate with latent moving averages $(X_{1,t}),\ldots,(X_{J,t})$ as in Definition \ref{def:aggma}, and let $\bX_t$ as in \eqref{def:multi} for any $m\ge0$, $h\ge1$.
	Then, $\bX_t$ is $\alpha$-stable and its spectral representation $(\Gamma,\bmu^0)$ on the Euclidean unit sphere $S_{m+h+1}$ writes
	\begin{align}
     \Gamma & = \sum_{j=1}^J \sum_{\vartheta\in S_1}\sum_{k\in\mathbb{Z}} w_{j,\vartheta}\pi_j^\alpha \|\bd_{j,k}\|_e^\alpha \delta_{\left\{\dfrac{\vartheta\bd_{j,k}}{\|\bd_{j,k}\|_e}\right\}}, \label{eq:spec_repres_agg} \\
     \bmu^0 & = - \mathds{1}_{\{\alpha=1\}}\frac{2}{\pi}\sum_{j=1}^J\sum_{k\in\mathbb{Z}}\pi_j\beta_j\bd_{j,k}\ln\|\pi_j\bd_{j,k}\|_e,\nonumber
	\end{align}
	where $\bd_{j,k}=(d_{j,k+m},\ldots,d_{j,k},d_{j,k-1},\ldots,d_{j,k-h})$, $w_{j,\vartheta}=(1+\vartheta\beta_j)/2$, for any $k\in\mathbb{Z}$, $j=1,\ldots,J$, $\vartheta\in S_1$, and if $\bd_{j,k}=\boldsymbol{0}$, the term vanishes by convention from the sums.
	\end{lem}
	\begin{rem}
	\rm 	Notice that $\Gamma = \sum_{j=1}^J\pi_j^\alpha \Gamma_j$, where $\Gamma_j$ denotes the spectral measure of the path $\bX_{j,t}$ from the moving average $(X_{j,t})$, $j=1,\ldots,J$, which is of the form \eqref{def:spectral}.
	\end{rem}
	If all the $\bX_{j,t}$'s are symmetric ($\beta_j=0$ for all $j$), then $\bX_t$ and $\Gamma$ are symmetric as well, but the reciprocal however does not hold true.
	The measure $\Gamma$ will be symmetric if and only if $\sum_{j=1}^J\pi_j^\alpha \Big(\Gamma_j(A)-\Gamma_j(-A)\Big)=0$ for any Borel set $A\subset S_{m+h+1}$.
	The latter condition is necessary and sufficient for $\bX_t$ to be symmetric in the case where $\alpha\ne1$, whereas for $\alpha=1$, it guarantees that $\bX_t$ will be symmetric up to an additive shifting, as $\bmu^0$ may be non-zero.
	The symmetry of paths intervenes in the representability conditions provided in the following lemma.

	\begin{lem}\label{le:spec_repres_agg_vec}
	Let $(X_t)$ an $\alpha$-stable aggregate with latent moving averages $(X_{1,t}),\ldots,(X_{J,t})$ as in Definition \ref{def:aggma}. Let $m\ge0$, $h\ge1$, $\bX_t$ as in \eqref{def:multi}, and $\|\cdot\|$ be a semi-norm on $\mathbb{R}^{m+h+1}$ satisfying \eqref{asu:semin}.\\
	When either $\alpha\ne1$ or $\bX_t$ S1S, the vector $\bX_t$ is representable on $C_{m+h+1}^{\|\cdot\|}$ if and only if \eqref{eq:ma_representable} holds with $m$ for all sequences $(d_{j,k})_k$, $j=1,\ldots,J$.\\
	For $\alpha=1$ and $\bX_t$ asymmetric, the vector $\bX_t$ is representable on $C_{m+h+1}^{\|\cdot\|}$ if and only if \eqref{eq:ma_representable} and \eqref{eq:additionallourd} hold with $m$ and $h$ for all sequences $(d_{j,k})_k$, $j=1,\ldots,J$.
	\end{lem}
	The notion of past-representability in Definition \ref{def:ma_representable} straightforwardly encompasses the case of stable aggregated processes and the next proposition provides a characterisation. 
	In view of Lemma \eqref{prop238}, the condition for the representability of a path $\bX_t$ on a unit cylinder changes according to whether it is symmetric or not in the case $\alpha=1$. 
	\begin{prop}
		\label{prop:aggma}
		Let $(X_t)$ an $\alpha$-stable aggregate with latent moving averages $(X_{1,t}),\ldots,(X_{J,t})$ as in Definition \ref{def:aggma}.\\
		$(\iota)$ Define for $j=1,\ldots,J$ the sets $\mathcal{M}_j = \{m\ge1: \hspace{0.2cm} \exists k\in\mathbb{Z}, \hspace{0.2cm} d_{j,k+m}=\ldots=d_{j,k+1}=0, \hspace{0.2cm} d_{j,k}\ne0\}$, and
		\begin{equation}\label{def:m0j}
			m_{0,j} = \left\{
			\begin{array}{cc}\sup\,\mathcal{M}_j, &\text{if} \hspace{0.5cm} \mathcal{M}_j \ne \emptyset,\\
				0 , & \text{if} \hspace{0.5cm} \mathcal{M}_j = \emptyset.
			\end{array}
			\right.
		\end{equation}
		
		$(a)$ For $\alpha\ne1$, the aggregated process $(X_t)$ is past-representable if and only if $(X_{j,t})$ is past-\indent representable for all $j=1,\ldots,J$, i.e.,
		\begin{align}\label{cond:foralljmoj}
		\sup\limits_{j=1,\ldots,J} m_{0,j} < +\infty.
		\end{align}
		\indent Moreover, letting $m\ge0$, $h\ge1$, $(X_t)$ is $(m,h)$-past-representable if and only if \eqref{cond:foralljmoj} holds \indent and $m\ge \max\limits_{j=1,\ldots,J}m_{0,j}$.\\
		
		$(b)$ For $\alpha=1$, the process $(X_t)$ is past-representable if and only if \eqref{cond:foralljmoj} holds and there exists \indent a pair $(m,h)$, $m\ge \max\limits_{j=1,\ldots,J}m_{0,j}$, $h\ge1$ such that either
		$$
		\bX_t \hspace{0.2cm} \text{is} \hspace{0.2cm} S1S, \hspace{1cm} \text{or}, \hspace{1cm} \bX_t \hspace{0.1cm} \text{asymmetric and \eqref{eq:additionallourd} holds for all sequences} \hspace{0.1cm} (d_{j,k})_k,
		$$
		\indent where $\bX_t$ generically denotes a vector as in \eqref{def:multi}. If such a pair exists, then the process $(X_t)$ \indent is $(m,h)$-past-representable.\\
		
		\noindent $(\iota\iota)$ Let $\|\cdot\|$ a semi-norm satisfying \eqref{asu:semin} and assume that $(X_t)$ is $(m,h)$-past-representable for some $m\ge0$, $h\ge1$. 
		The spectral representation $(\Gamma^{\|\cdot\|},\bmu^{\|\cdot\|})$ of the vector $\bXt=(X_{t-m},\ldots,X_t,X_{t+1},\ldots,X_{t+h})$ on $C_{m+h+1}^{\|\cdot\|}$ is then given by \eqref{eq:spec_repres_agg} with the Euclidean norm $\|\cdot\|_e$ replaced by the semi-norm $\|\cdot\|$.
	\end{prop}
	\begin{rem}
	\rm The necessary condition \eqref{cond:foralljmoj} extends what was noticed in the case of non-aggregated moving averages, namely, that anticipativeness is a minimal requirement for past-representability.
	A single non-anticipative latent moving average is enough to render the aggregated process not past-representable, regardless of the other latent components.
	\end{rem}
	\begin{rem}
		\rm For $\alpha\ne1$, the past-representability of an aggregated process is equivalent to that of its latent moving averages, but this does not seem to hold in general for $\alpha=1$.
		In the latter case however, if all the latent moving averages are symmetric, that is, $\beta_1=\ldots=\beta_J=0$, then the paths $\bX_t$ are S1S for any $m\ge0$, $h\ge1$ and $(\iota)(b)$ collapses to $(\iota)(a)$.
	\end{rem}
	The representability condition also simplifies in the case of aggregated ARMA processes and requires each latent ARMA process to be anticipative.
	\begin{cor}\label{cor:agg_arma}\rm
	For any $j=1,\ldots,J$, let $(X_{j,t})$ be the ARMA strictly stationary solution of $\psi_j(F)\phi_j(B)X_{j,t}=\Theta_j(F)H_j(B)\varepsilon_{j,t}$, as in Corollary \ref{cor:mar}, with mutually independent sequences $\varepsilon_{j,t}\stackrel{i.i.d.}{\sim}\mathcal{S}(\alpha,\beta_j,1,0)$.
	Define $X_t=\sum_{j=1}^J\pi_jX_{j,t}$ for any positive scalings $(\pi_j)_j$. Then, for any $\alpha\in(0,2)$, $(\beta_1,\ldots,\beta_J)\in[-1,1]^J$, the following statements are equivalent:\\
	$(\iota)\hspace{0.45cm}(X_t) \hspace{0.1cm} \text{is past-representable},$\\
	$(\iota\iota)\hspace{0.3cm}\inf\limits_{j}\deg(\psi_j)\ge1,$ \\
	$(\iota\iota\iota)\hspace{0.2cm}\sup\limits_{j}m_{0,j}<+\infty,$\\
	with the $m_{0,j}$'s as in \eqref{def:m0j}.
	Moreover, letting $m\ge0$, $h\ge1$, the aggregated process $(X_t)$ is $(m,h)$-past-representable if and only if for any $j=1,\ldots,J$, $m_{0,j}<+\infty$, and $m\ge \max\limits_{j} m_{0,j}$.
	\end{cor}
	
	\section{Conditional tail distribution of stable aggregates}
	\label{sec:application}
	
	In this section, we will derive the tail conditional distribution of linear stable processes for which Proposition \ref{prop:cond_tail} will be applicable. 
	The case of a general past-representable stable aggregate is considered as well as particular examples.\\
	To be relevant for the prediction framework, the Borel set $B$ appearing in Proposition \ref{prop:cond_tail} has to be chosen such that the conditioning event $\{\|\bXt\|>x\}\cap\{\bXt/\|\bXt\|\in B\}$ is independent of the future realisations $X_{t+1},\ldots,X_{t+h}$. For $\|\cdot\|$ a semi-norm on $\mathbb{R}^{m+h+1}$ satisfying \eqref{asu:semin}, denote $S_{m+1}^{\|\cdot\|}=\{(s_{-m},\ldots,s_0)\in\mathbb{R}^{m+1}: \hspace{0.3cm} \|(s_{-m},\ldots,s_0,0,\ldots,0)\|=1\}$.\footnote{The set $S_{m+1}^{\|\cdot\|}$ corresponds to the unit sphere of $\mathbb{R}^{m+1}$ relative to the restriction of $\|\cdot\|$ to the first $m+1$ dimensions.}
	Then, for any Borel set $V\subset S_{m+1}^{\|\cdot\|}$, define the Borel set $B(V)\subset C_{m+h+1}^{\|\cdot\|}$ as
	\begin{align*}
		B(V) = V \times \mathbb{R}^{h}. 
	\end{align*}
	Notice in particular that for $V=S_{m+1}^{\|\cdot\|}$, we have $B(V) = C_{m+1}^{\|\cdot\|}$.
	In the following, we will use Borel sets of the above form to condition the distribution of the complete vector $\bXt/\|\bXt\|$ on the observed <<shape>> of the past trajectory. The latter information is contained in the Borel set $V$, which we will typically assume to be some small neighbourhood on $S_{m+1}^{\|\cdot\|}$.
	It will be useful in the following to notice that
	\begin{align*}
	V \times \mathbb{R}^{h} = \left\{\bs\in C_{m+h+1}^{\|\cdot\|}: \hspace{0.3cm} f(\bs) \in V \right\},
	\end{align*}
	where $f$ the function defined by 
	\begin{align}\label{def:f}
		f:\begin{array}{ccc}
			\mathbb{R}^{m+h+1} & \longrightarrow & \mathbb{R}^{m+1}  \\
			(x_{-m},\ldots,x_0,x_1,\ldots,x_h) & \longmapsto & (x_{-m},\ldots,x_0)
		\end{array}
		. 
	\end{align}
	
	\subsection{Stable aggregates: general case}
	
	Let $(X_t)$ an $\alpha$-stable aggregate as in Definition \ref{def:aggma} (possibly a moving average if $J=1$).
	Assume $(X_t)$ is $(m,h)$-past-representable, for some $m\ge 0$, $h\ge1$ and let $\bXt$ as in \eqref{def:multi}. 
	Denoting $\Gamma^{\|\cdot\|}$ the spectral measure of $\bX_t$ on the unit cylinder $C_{m+h+1}^{\|\cdot\|}$ for some semi-norm satisfying \eqref{asu:semin}, we know by Proposition \ref{prop:aggma} $(\iota\iota)$, that $\Gamma^{\|\cdot\|}$ is of the form
	\begin{equation}\label{eq:specmeasu}
		\Gamma^{\|\cdot\|} = \sum_{j=1}^J\sum_{\vartheta\in S_1}\sum_{k\in\mathbb{Z}}w_{j,\vartheta}\pi_j^\alpha \|\bd_{j,k}\|^\alpha \delta_{\left\{\dfrac{\vartheta \bd_{j,k}}{\|\bd_{j,k}\|}\right\}}. 
	\end{equation}
	\begin{prop}\label{prop:relativprob}
		Under the above assumptions, we have
		\begin{align}\label{eq:relativprob}
			\mathbb{P}_x^{\|\cdot\|}\Big(\bXt,A\Big|B(V)\Big) \underset{x\rightarrow+\infty}{\longrightarrow} \dfrac{\Gamma^{\|\cdot\|}\Bigg(\bigg\{\dfrac{\vartheta\bd_{j,k}}{\|\bd_{j,k}\|}\in A: \hspace{0.2cm} \dfrac{\vartheta f(\bd_{j,k})}{\|\bd_{j,k}\|}\in V\bigg\}\Bigg)}{\Gamma^{\|\cdot\|}\Bigg(\bigg\{\dfrac{\vartheta\bd_{j,k}}{\|\bd_{j,k}\|}\in C_{m+h+1}^{\|\cdot\|}: \hspace{0.2cm} \dfrac{\vartheta f(\bd_{j,k})}{\|\bd_{j,k}\|}\in V\bigg\}\Bigg)},
		\end{align}
		for any Borel sets $A\subset C_{m+h+1}^{\|\cdot\|}$, $V\subset S_{m+1}^{\|\cdot\|}$ such that $\bigg\{\dfrac{\vartheta\bd_{j,k}}{\|\bd_{j,k}\|}\in C_{m+h+1}^{\|\cdot\|}: \hspace{0.2cm} \dfrac{\vartheta f(\bd_{j,k})}{\|\bd_{j,k}\|}\in V\bigg\}\ne\emptyset$, $\Gamma^{\|\cdot\|}\Big(\partial(A\cap B(V))\Big)=\Gamma^{\|\cdot\|}(\partial B(V))=0$, 
		where $B(V)=V \times \mathbb{R}^{h}$ and $f$ is as in \eqref{def:f}. 
	\end{prop}
	\begin{rem}
		\rm
		
		$(\iota)$ Setting $V=S_{m+1}^{\|\cdot\|}$, 
		and $A$ an arbitrarily small closed neighbourhood of all the points $(\vartheta \bd_{j,k}/\|\bd_{j,k}\|)_{\vartheta,j,k}$, we can see that $\lim\limits_{x\rightarrow+\infty} \mathbb{P}\Big(\bXt/\|\bXt\|\in A\Big|\|\bXt\|>x\Big) = 1$. 
		In other terms, when far from central values, the trajectory of process $(X_t)$ necessarily features patterns of the same shape as some $\vartheta \bd_{j,k}/\|\bd_{j,k}\|$, which is a finite piece of a moving average's coefficient sequence. 
		The index $j$ indicates from which of the $J$ underlying moving averages the pattern stems from, the index $k$ points to which piece $(d_{j,k+m},\ldots,d_{j,k},d_{j,k-1},\ldots,d_{j,k-h})$ of this moving average it corresponds, and $\vartheta\in\{-1,+1\}$ indicates whether the pattern is flipped upside down (in case the extreme event is driven by a negative value of an error $(\varepsilon_{j,\tau})$).
		The likelihood of a pattern $\vartheta \bd_{j,k}/\|\bd_{j,k}\|$ can be evaluated by setting $A$ to be a small neighbourhood of that point. 
		In particular, only one pattern $\bd_k/\|\bd_k\|$ can appear through time for $J=1$ (up to a time shift and sign flipping). This is no longer the case in general for $J\ge2$, where the shape of each extreme event appears as if being drawn from a collection of patterns.\\
		$(\iota\iota)$ In view of point $(\iota)$, the observed path $(X_{t-m},\ldots,X_{t-1},X_t)/\|\bXt\|$ will \textit{a fortiori} be of the same shape as some $\vartheta (d_{j,k+m},\ldots,d_{j,k+1},d_{j,k})/\|\bd_{j,k}\|$ when an extreme event will approach in time.
		Observing the initial part of the pattern can give information about the remaining unobserved piece: the conditional likelihood of the latter can be assessed by setting $V$ to be a small neighbourhood of the observed pattern. 
	\end{rem}
	\begin{rem}\label{rem:uncertainty}
		\rm
		The tail conditional distribution given in \eqref{eq:relativprob} highlights three types of uncertainty/approximation for prediction:\footnote{The considerations developed in this remark focus solely on the probabilistic uncertainty of the prediction assuming that the process $(X_t)$ is entirely known, that is, no parameter nor any sequence $(d_{j,k})$ has to be inferred from data.}\\
		$(\iota)$ In practice, events of the type $\{(X_{t-m},\ldots,X_{t-1},X_t)/\|\bXt\|=\vartheta (d_{j,k+m},\ldots,d_{j,k+1},d_{j,k})/\|\bd_{j,k}\|\}$ have probability zero of occurring, and only noisy observations such as $(X_{t-m},\ldots,X_{t-1},X_t)/\|\bXt\|\approx\vartheta (d_{j,k+m},\ldots,d_{j,k+1},d_{j,k})/\|\bd_{j,k}\|$ are available on a realised trajectory. 
		The choice of an adequate conditioning neighbourhood $V$ in \eqref{eq:relativprob} given a piece of trajectory will thus have to rely on a statistical approach.
		One could envision tests of hypotheses to determine whether a piece of realised (noisy) trajectory <<is more similar>> to a certain pattern 1 or to an other pattern 2, or whether it <<is more similar>> to a certain pattern 1 rather than any patterns in a certain collection.\\
		$(\iota\iota)$ Even for an arbitrarily small neighbourhood $V$ --that is, even if the observed path can be confidently identified with a particular pattern-- uncertainty regarding the future trajectory may remain.
		It could indeed be that several patterns $\vartheta \bd_{j,k}/\|\bd_{j,k}\|$ coincide on their first $m+1$ components, but differ by the last $h$.
		The stable anticipative AR(1) and its aggregated version are typical examples of this phenomenon that will be studied in the next section.\\
		$(\iota\iota\iota)$ The tail conditional distribution \eqref{eq:relativprob} is an asymptotic behaviour as the (semi-)norm of $\bXt$ grows infinitely large. 
		It is thus only an approximation of the true dynamics during extreme events.
		It would be interesting to obtain a finer asymptotic development in $x$ of the above convergence to gauge the approximation error of the true conditional distribution.
		It would be especially useful to quantify how far from/how variable around the predicted patterns the future path can be.
	\end{rem}
	
	\subsection{Aggregation of AR(1)}
	
	We now consider $(X_t)$ the aggregation of stable anticipative AR(1) processes introduced in \cite{gz17} defined by
	\begin{align}\label{def:aggar1}
		X_{t} = \sum_{j=1}^J\pi_jX_{j,t},\quad  X_{j,t} & = \rho_j X_{j,t+1}+\varepsilon_{j,t}, \quad 0<|\rho_j|<1, \quad j=1,\ldots,J
	\end{align}
	where $\pi_j>0$ for any $j$, and $(\varepsilon_{j,t})_{t\in\mathbb{Z}}\stackrel{i.i.d.}{\sim}\mathcal{S}(\alpha,\beta_j,1,0)$ are mutually independent i.i.d. sequences. We assume without loss of generality that the $\rho_j$'s are distinct. 
	For each anticipative AR(1) with parameter $\rho_j$, the moving average coefficients are of the form $(\rho_j^k\mathds{1}_{\{k\ge0\}})_k$, and thus, $m_{0,j}=0$ for all $j$, where the $m_{0,j}$'s are given in \eqref{def:m0j}.
	By Corollary \eqref{cor:agg_arma}, we know for any $m\ge0$, $h\ge1$, the aggregated process $(X_t)$ is $(m,h)$-past-representable.
	The spectral measures of paths $\bX_t$ simplify and charge finitely many points. 
	Their forms are given in the next lemma.
	\begin{lem}
		\label{le:djk_aggar1}
		Let $(X_t)$ be an aggregation of $\alpha$-stable anticipative AR(1) processes as in \eqref{def:aggar1}. 
		Letting $\bXt$ as in \eqref{def:multi} for $m\ge0$, $h\ge1$, its spectral measure on $C_{m+h+1}^{\|\cdot\|}$ for a seminorm satisfying \eqref{asu:semin} is given by
		\begin{align}\label{eq:specaggar1}
			\Gamma^{\|\cdot\|} & =  \sum_{\vartheta\in S_1} \Bigg[w_\vartheta\delta_{\{(\vartheta,0,\ldots,0)\}} + \sum_{j=1}^J  \pi_j^\alpha  \bigg( w_{j,\vartheta}\sum_{k=-m+1}^{h-1}\|\bd_{j,k}\|^\alpha \delta_{\left\{\dfrac{\vartheta \bd_{j,k}}{\|\bd_{j,k}\|}\right\}}  + \dfrac{\Bar{w}_{j,\vartheta}}{1-|\rho_j|^\alpha} \|\bd_{j,h}\|^\alpha  \delta_{\left\{\dfrac{\vartheta \, \bd_{j,h}}{\|\bd_{j,h}\|}\right\}}\bigg)\Bigg],
		\end{align}
		where for all $\vartheta\in S_1$, $j\in\{1,\ldots,J\}$ and $-m+1\le k \le h$,
		\begin{align*}
		\bd_{j,k} & = (\rho_j^{k+m}\mathds{1}_{\{k\ge-m\}},\ldots,\rho_j^{k}\mathds{1}_{\{k\ge0\}},\rho_j^{k-1}\mathds{1}_{\{k\ge1\}},\ldots,\rho_j^{k-h}\mathds{1}_{\{k\ge h\}}),\\
		w_{j,\vartheta} & = (1+\vartheta\beta_j)/2,\\
		w_\vartheta & = \sum_{j=1}^J\pi_j^\alpha w_{j,\vartheta},\\
		\Bar{w}_{j,\vartheta} & = (1+\vartheta\Bar{\beta}_j)/2,\\
		\Bar{\beta}_j & = \beta_j\dfrac{1-\rho_j^{<\alpha>}}{1-|\rho_j|^\alpha},
		\end{align*}
		and if $h=1$ and $m=0$, the sum $\sum_{k=-m+1}^{h-1}$ vanishes by convention.
	\end{lem}
	The next proposition provides the tail conditional distribution of future paths in the case where the $\rho_j$'s are positive. 
	Let us first introduce useful neighbourhoods of the distinct charged points of $\Gamma^{\|\cdot\|}$.
	Denote $\bd_{0,-m}=(\overbrace{1,0,\ldots,0}^{m+h+1})$ so that the charged points of $\Gamma^{\|\cdot\|}$ are all of the form $\vartheta \bd_{j,k}/\|\bd_{j,k}\|$ with indexes $(\vartheta,j,k)$ in the set 
	$\mathcal{I}:=S_1\times\Big(\{1,\ldots,J\}\times\{-m,h\}\cup \{(0,-m)\}\Big)$.
	With $f$ as in \eqref{def:f}, define for any $(\vartheta_0,j_0,k_0)\in \mathcal{I}$, the set $V_0$ as any closed neighbourhood of $\vartheta_0f(\bd_{j_0,k_0})/\|\bd_{j_0,k_0}\|$ such that
	\begin{equation}\label{eq:v0}
		\forall (\vartheta',j',k')\in\mathcal{I}, \hspace{1.4cm} \dfrac{\vartheta' f(\bd_{j',k'})}{\|\bd_{j',k'}\|}\in V_0 \hspace{0.3cm} \Longrightarrow \hspace{0.3cm} \dfrac{\vartheta' f(\bd_{j',k'})}{\|\bd_{j',k'}\|} = \dfrac{\vartheta_0 f(\bd_{j_0,k_0})}{\|\bd_{j_0,k_0}\|},
	\end{equation}
	In other terms, $V_0\times\mathbb{R}^d$ is a subset of $C_{m+h+1}^{\|\cdot\|}$ in which the only points charged by $\Gamma^{\|\cdot\|}$ all have the first $(m+1)$\textsuperscript{th} coinciding with $\vartheta_0f(\bd_{j_0,k_0})/\|\bd_{j_0,k_0}\|$. 
	Define also $A_{\vartheta,j,k}$ for any $(\vartheta,j,k)$ as any closed neighbourhood of $\vartheta\bd_{j,k}/\|\bd_{j,k}\|$ which does not contain any other charged point of $\Gamma^{\|\cdot\|}$, that is,
	\begin{equation}\label{eq:avjk}
		\forall (\vartheta',j',k')\in\mathcal{I}, \hspace{2cm} \dfrac{\vartheta'\bd_{j',k'}}{\|\bd_{j',k'}\|}\in A_{\vartheta,j,k} \hspace{0.3cm} \Longrightarrow \hspace{0.3cm} (\vartheta',j',k') = (\vartheta,j,k).
	\end{equation}
	\begin{prop}
		\label{prop:aggar1_pred}
		Let $(X_t)$ be an aggregation of $\alpha$-stable anticipative AR(1) processes as in \eqref{def:aggar1} with $\rho_j\in(0,1)$ for all $j$'s.
		Let $\bXt$, the $\bd_{j,k}$'s and the spectral measure of $\bXt$ be as given in Lemma \ref{le:djk_aggar1}, for any $m\ge0$, $h\ge1$.
		Let $V_0$ be any small closed neighbourhood of $\vartheta_0 f(\bd_{j_0,k_0})/\|\bd_{j_0,k_0}\|$ in the sense of \eqref{eq:v0} for some $(\vartheta_0,j_0,k_0)\in\mathcal{I}$ and let $B(V_0)=V_0\times \mathbb{R}^h$. 
		Then, with $A_{\vartheta,j,k}$ an arbitrarily small neighbourhood of some $\vartheta\bd_{j,k}/\|\bd_{j,k}\|$ as in
		\eqref{eq:avjk}, 
		the following hold.\\
		$(\iota)$ Case $m\ge1$.\\ 
		\indent $(a)$ If $0\le k_0 \le h$:
		\begin{align*}
			\mathbb{P}_x^{\|\cdot\|}\Big(\bXt,A_{\vartheta,j,k}\Big|B(V_0)\Big) & \underset{x\rightarrow\infty}{\longrightarrow}  
			\left\{
			\begin{array}{lc}
				|\rho_{j_0}|^{\alpha k} (1-|\rho_{j_0}|^\alpha)\delta_{\vartheta_0}(\vartheta)\delta_{j_0}(j),&  0\le k \le h-1,\\
				& \\
				|\rho_{j_0}|^{\alpha h}\delta_{\vartheta_0}(\vartheta)\delta_{j_0}(j), &  k=h.
			\end{array}
			\right.
		\end{align*}
		\indent $(b)$ If $-m\le k_0 \le -1$:
		\begin{align*}
			\mathbb{P}_x^{\|\cdot\|}\Big(\bXt,A_{\vartheta,j,k}\Big|B(V_0)\Big) \underset{x\rightarrow\infty}{\longrightarrow}  \delta_{\vartheta_0}(\vartheta)\delta_{j_0}(j)\delta_{k_0}(k).
		\end{align*}
		$(\iota\iota)$ Case $m=0$.\\
		\begin{align*}
			\mathbb{P}_x^{\|\cdot\|}\Big(\bXt,A_{\vartheta,j,k}\Big|B(V_0)\Big) \underset{x\rightarrow\infty}{\longrightarrow} \left\{
			\begin{array}{lc}
			    \dfrac{\sum_{i=1}^J \pi_{i}^\alpha w_{i,\vartheta_0}}{\sum_{i=1}^J p_{i,\vartheta_0}}\delta_{\{\vartheta_0\}}(\vartheta), &  k=0\\
			    & \\
				\dfrac{p_{j,\vartheta_0}}{\sum_{i=1}^J p_{i,\vartheta_0}}|\rho_{j}|^{\alpha k}(1-|\rho_{j}|^\alpha)\delta_{\{\vartheta_0\}}(\vartheta), &  1\le k \le h-1,\\
				& \\
				\dfrac{p_{j,\vartheta_0}}{\sum_{i=1}^J p_{i,\vartheta_0}}|\rho_j|^{\alpha h}\delta_{\{\vartheta_0\}}(\vartheta), &  k=h,
			\end{array}
			\right.
		\end{align*}
		with $p_{j,\vartheta_0} = \pi_{j}^\alpha w_{j,\vartheta_0}/(1-|\rho_{j}|^\alpha)$. 
	\end{prop}
	\begin{rem}\rm
		For $m\ge1$, that is, if the observed path is assumed to be of length at least 2, there is a significant difference between whether  $k_0\in\{0,\ldots,h\}$ or $k_0\in\{-m,\ldots,-1\}$. 
		For the latter, the asymptotic probability of the whole path $\bXt/\|\bXt\|$ being in an arbitrarily small neighbourhood of $\vartheta \bd_{j,k}/\|\bd_{j,k}\|$ is 1 if and only if $\vartheta=\vartheta_0$, $j=j_0$, $k=k_0$: given the observed path, the shape of the future trajectory is fully determined.
		For the former, this probability is strictly positive if and only if $\vartheta=\vartheta_0$ and $j=j_0$, but the observed pattern is compatible with several distinct future paths.
		One can see why this is the case from the form of the sequences $\bd_{j,k}/\|\bd_{j,k}\|$ and of their restrictions to the first $m+1$ components $f(\bd_{j,k})/\|\bd_{j,k}\|$. On the one hand (omitting $\vartheta$),
		\begin{align*}
			\dfrac{\bd_{j,k}}{\|\bd_{j,k}\|} & = 
			\left\{
			\begin{array}{lll}
				\dfrac{(\overbrace{\rho_j^{k+m},\ldots,\rho_j^{k}}^{m+1},\overbrace{\rho_j^{k-1},\ldots,\rho_j,1,0,\ldots,0}^{h})}{\|(\rho_j^{k+m},\ldots,\rho_j^{k},\rho_j^{k-1},\ldots,\rho_j,1,0,\ldots,0)\|}, & \text{for} & \hspace{0.5cm} k\in\{0,\ldots,h\},\\
				& & \\
				\dfrac{(\rho_j^{k+m},\ldots,\rho_j,1,0,\ldots,0,0,\ldots,0)}{\|(\underbrace{\rho_j^{k+m},\ldots,\rho_j,1,0,\ldots,0}_{m+1},\underbrace{0,\ldots,0}_{h})\|}, & \text{for} & \hspace{0.5cm} k\in\{-m,\ldots,-1\}.
			\end{array}
			\right.
		\end{align*}
		We can notice that all the above sequences are pieces of  explosive exponentials, terminated at some coordinate.
		For $k\in\{0,\ldots,h\}$, the first zero component --the <<crash of the bubble>>--, is situated at or after the $(m+2)$\textsuperscript{th} component, whereas for $k\in\{-m,\ldots,-1\}$, it is situated at or before the $(m+1)$\textsuperscript{th}.
		Using the homogeneity of the semi-norm and \eqref{asu:semin_intro}, we have on the other hand that
		\begin{align*}
			\dfrac{f(\bd_{j,k})}{\|\bd_{j,k}\|} & = 
			\left\{
			\begin{array}{lll}
			\dfrac{(\overbrace{\rho_j^{m},\ldots,\rho_j,1}^{m+1})}{\|(\underbrace{\rho_j^{m},\ldots,\rho_j,1}_{m+1},\underbrace{0,\ldots,0,0,\ldots,0}_{h})\|}, & \text{for} & \hspace{0.5cm} k\in\{0,\ldots,h\},\\
				& & \\
				\dfrac{(\overbrace{\rho_j^{k+m},\ldots,\rho_j,1,0,\ldots,0}^{m+1})}{\|(\underbrace{\rho_j^{k+m},\ldots,\rho_j,1,0,\ldots,0}_{m+1},\underbrace{0,\ldots,0}_{h})\|}, & \text{for} & \hspace{0.5cm} k\in\{-m,\ldots,-1\}.
			\end{array}
			\right.
		\end{align*}
		Thus, conditioning the trajectory on the event $\{f(\bXt)/\|\bXt\| \approx f(\bd_{j_0,k_0})/\|\bd_{j_0,k_0}\|\}$ for some $k_0\in\{-m,\ldots,-1\}$ amounts to condition on the burst of a bubble being observed in the past trajectory with no new bubble forming yet, which allows to identify exactly the position of the pattern on the $j$\textsuperscript{th} moving average's coefficient sequence.\\
		When conditioning with $k_0\in\{0,\ldots,h\}$ however, the crash date is not observed and can happen either in the next $h-1$ periods, or after the $h$\textsuperscript{th}. 
		However, the shape of the observed path is that of a piece of exponential with growth rate $\rho_j^{-1}$ regardless of the remaining time before the burst, which leaves several future paths possible. 
		One can quantify the likelihood of each potential scenario: the quantity $|\rho_j|^{\alpha k}(1-|\rho_j|^{\alpha})$ corresponds to the probability that the bubble will peak in exactly $k$ periods ($0\le k<h$), and $|\rho_j|^{\alpha h}$ corresponds to the probability that the bubble will last at least $h$ more periods. 
	\end{rem}
	\begin{rem}
		\rm
		$(\iota)$ The previous remark confirms the interpretation of the conditional moments proposed in \cite{f18} for the stable anticipative AR(1) case ($J=1$). 
		It also extends it in two ways: $(\iota)$ by accounting for paths rather than point prediction, $(\iota\iota)$ by showing that the aggregation of AR(1) processes also features killed exponential explosive episodes but with various growth rates and crash probabilities. 
		Proposition \ref{prop:aggar1_pred} furthermore shows that asymptotically, as few as two observations are sufficient to identify the growth rate $\rho_j^{-1}$ of an ongoing extreme episode,\footnote{ This holds asymptotically in the (semi-)norm of the observed path, but in practice it can be expected that the noise surrounding the trajectory will make this identification difficult with only two observations. 
			Longer path lengths (higher $m$) may provide robustness to the identification, but could also incorporate some bias by taking into account past extreme events, such as now-collapsed bubbles.
			One can suspect a bias-variance trade-off when searching for an optimal choice of $m$.} and the conditional dynamics within this given event will be similar to that of a simple AR(1) with corresponding parameter.
		An identification of the growth rate in the early developments of the bubble appears possible, allowing to infer in advance the odds of crashes.\\
		$(\iota\iota)$ Notice that for $m=0$ (only the present value is assumed to be observed), no pattern can be observed but only the sign of the shock. 
		Hence, the growth rate $\rho_{j_0}^{-1}$ of the ongoing event is unidentifiable, which is reflected in the fact that the asymptotic probabilities of paths with growth rates $\rho_{j}^{-1}$, $j\ne j_0$, are positive (case $(\iota\iota)$ of Proposition \ref{prop:aggar1_pred}).
	\end{rem}
	
	\subsection{Two examples: the anticipative AR(2) and fractionally integrated AR}
	
	We focus here on two processes which both share the peculiar property of having a $0$-$1$ tail conditional distribution whenever the observed path is of length at least 2 (i.e., $m\ge1$): the anticipative AR(2) and the anticipative fractionally integrated AR.
	For an adequate choice of the parameters, the former can generate bubble-like trajectories with accelerating or decelerating growth rate and the latter can accommodate hyperbolic bubbles.
	In contrast with the anticipative AR(1), these bubbles do not display an exponential profile but still feature an inflation-peak-collapse behaviour.
	The two processes are defined as follows.
	
	\noindent\textbf{Anticipative AR(2)}
	
	\noindent The anticipative AR(2) is the strictly stationary solution of
	\begin{align}\label{def:ar2}
		(1-\lambda_1F)(1-\lambda_2F)X_t = \varepsilon_t, \hspace{1cm} \varepsilon_t \stackrel{i.i.d.}{\sim} \mathcal{S}(\alpha,\beta_,\sigma,0),
	\end{align}
	where $\lambda_i\in\mathbb{C}$ and $0<|\lambda_i|<1$ for i=1,2. 
	In case $\lambda_i\in\mathbb{C}\setminus\mathbb{R}$, $i=1,2$, we impose that $\lambda_1=\Bar{\lambda}_2$ to ensure $(X_t)$ is real-valued. 
	We further assume that $\lambda_1+\lambda_2\ne0$, to exclude the cases where $(X_{2t})$ and $(X_{2t+1})$ are independent anticipative AR(1) processes. 
	The solution of \eqref{def:ar2} admits the moving average representation $X_t=\sum_{k\in\mathbb{Z}} d_k \varepsilon_{t+k}$ with
	\begin{align}
		d_k = \left\{
		\begin{array}{cl}
			\dfrac{\lambda_1^{k+1}-\lambda_2^{k+1}}{\lambda_1-\lambda_2} \hspace{0.1cm}\mathds{1}_{\{k\ge0\}}, & \hspace{0.5cm} \text{if} \hspace{0.5cm} \lambda_1\ne\lambda_2, \\
			(k+1)\lambda^k\hspace{0.1cm}\mathds{1}_{\{k\ge0\}}, & \hspace{0.5cm} \text{if} \hspace{0.5cm} \lambda_1=\lambda_2=\lambda.
		\end{array}
		\right.\label{eq:ar2dk}
	\end{align}
	
	\noindent\textbf{Anticipative fractionally integrated AR}
	
	\noindent The anticipative fractionally integrated AR process can be defined as the stationary solution of
	\begin{align}\label{def:frac}
		(1-F)^d X_t = \varepsilon_t, \hspace{1cm} \varepsilon_t \stackrel{i.i.d.}{\sim} \mathcal{S}(\alpha,\beta_,\sigma,0),
	\end{align}
	with $\alpha(d-1)<-1$. The solution of \eqref{def:frac} admits the moving average representation $X_t=\sum_{k=0}^{+\infty} d_k \varepsilon_{t+k}$ with
	\begin{align}\label{eq:fracdk}
		d_0 = 1, && \text{and} && d_k = 
		\dfrac{\Gamma(k+d)}{\Gamma(d)\Gamma(k+1)} \hspace{0.1cm}\mathds{1}_{\{k\ge0\}}, \hspace{0.5cm} \text{for}  \hspace{0.3cm} k\ne0,
	\end{align}
	where $\Gamma(\,\cdot\,)$ denotes --here only-- the Gamma function.\\
	
	\noindent It can be shown that both process are necessarily $(m,h)$-past-representable for $m\ge1$ and $h\ge1$.
	The $0$-$1$ tail conditional distribution property when the observed path is of length at least 2 is exhibited in the next proposition.
	\begin{prop}
		\label{prop:ar2_pred}
		Let $(X_t)$ be the $\alpha$-stable anticipative AR(2) (resp. fractionally integrated AR) as in \eqref{def:ar2}-\eqref{eq:ar2dk} (resp. \eqref{def:frac}-\eqref{eq:fracdk}). 
		For any $m\ge1$ and $h\ge1$, let $\bXt$ as in \eqref{def:multi} and $\bd_{k}=(d_{k+m},\ldots,d_{k},d_{k-1},\ldots,d_{k-h})$ where $(d_k)$ is as in \eqref{eq:ar2dk} (resp. \eqref{eq:fracdk}). Let $V_0$ a small neighbourhood of $\vartheta_0\bd_{k_0}/\|\bd_{k_0}\|$ as in \eqref{eq:v0} --where we drop the indexes $j$-- for some $\vartheta_0\in S_1$, $k_0\ge-m$, and let $B(V_0)=V_0\times\mathbb{R}^h$.
		Then,
		\begin{align*}
			\mathbb{P}_x^{\|\cdot\|}\Big(\bXt,A\Big|B(V_0)\Big) & \underset{x\rightarrow\infty}{\longrightarrow} 
			\left\{ 
			\begin{array}{cc}
			    1, & \hspace{1cm} \text{if} \hspace{0.3cm} \dfrac{\vartheta_0 \bd_{k_0}}{\|\bd_{k_0}\|}\in A, \\
			    0, & \hspace{1cm} \text{otherwise},
			\end{array}
			\right.
		\end{align*}
		for any closed neighbourhood $A\subset C_{m+h+1}^{\|\cdot\|}$ such that $\partial A \cap \{\vartheta\bd_{k}/\|\bd_{k}\|:\hspace{0.1cm}\vartheta\in S_1,\hspace{0.1cm}k\ge-m\}=\emptyset$.
	\end{prop}
	\begin{rem}
		\rm
		Contrary to the anticipative AR(1), the trajectories of the anticipative AR(2) and fractionally integrated processes do not leave room for undeterminancy of the future path. 
		Asymptotically, given any observed path of length at least 2, the shape of the future trajectory can be deduced deterministically.
		This holds even if the peak/collapse of a bubble is not yet present in the observed piece of trajectory. 
		Therefore, provided the current pattern is properly identified,\footnote{See point $(\iota)$ of Remark \ref{rem:uncertainty}.} it appears possible in the framework of these models to infer in advance the peak and crash dates of bubbles with very high confidence --in principle, with certainty.
		
	\end{rem}
	
	\section{A step towards multivariate processes}
	\label{sec:bivariate_AR1}
	
	A simple bi-dimensional process is considered in this section to highlight that the approach developed in this paper can be brought to the multivariate framework and that new properties can also emerge. 
	In essence, the process considered is a vector where each univariate component consists respectively of a stable anticipative AR(1) and a stable non-anticipative AR(1), and dependence between both is allowed. Surprisingly, the presence of a non-anticipative component will not be pathological here contrary to the univariate case studied above, and Proposition \ref{prop:cond_tail} will be applicable. 
	Formally, define $(\bX_t)$ for all $t\in\mathbb{Z}$ as
	\begin{align}
		\left\{
		\begin{array}{rl}
			\bX_t & = (X_{1,t},X_{2,t})', \\
			X_{1,t} & = \rho_1 X_{1,t+1} + \varepsilon_{1,t},\\
			X_{2,t} & = \rho_2 X_{2,t-1} + \varepsilon_{2,t},\\
			\bveps_t & =(\varepsilon_{1,t},\varepsilon_{2,t})' \hspace{0.2cm} \text{i.i.d.} \hspace{0.2cm} S\alpha S\hspace{0.2cm}\text{with spectral measure} \hspace{0.2cm} \Gamma_2 \hspace{0.2cm} \text{on} \hspace{0.2cm} S_2 \hspace{0.2cm} \text{and zero shift vector},
		\end{array}
		\right.
		\label{def:bivarAR1}
	\end{align}
	where $0<|\rho_i|<1$, $i=1,2$.\footnote{The $S\alpha S$ assumption on the i.i.d. sequence $(\bveps_t)$ is made for the sake of simplicity and implies that $\Gamma_2$ is itself symmetric.} 
	We again have in mind applying Proposition \ref{prop:cond_tail} to a vector composed of past and future realisations of $(\bX_t)$. Limiting ourselves to the simplest $m=0$ and $h=1$ case, we will consider a vector of the form $\boldsymbol{\underline{X}}_t:=(\bX_t',\bX_{t+1}')'$, where $\bX_t$ is the present observation and $\bX_{t+1}$ the one-step ahead realisation to predict.
	The next result shows that $\boldsymbol{\underline{X}}_t$ is $\alpha$-stable, in fact $S \alpha S$, and it provides a necessary and sufficient condition on $\Gamma_2$ for its representability on an appropriate unit cylinder.\footnote{For expository purposes, the form of the spectral representations on the Euclidean unit sphere and on the unit cylinder are relegated in the proofs in Appendix.}
	\begin{prop}\label{prop:reprbility_multivar}
		Let $(\bX_t)$ as in \eqref{def:bivarAR1}, the semi-norm $\|\cdot\|$ on $\mathbb{R}^4$ such that $\|(x_1,x_2,x_3,x_4)\|=\sqrt{x_1^2+x_2^2}$ for any $(x_1,x_2,x_3,x_4)\in\mathbb{R}^4$, 
		and denote $C_4^{\|\cdot\|}$ its corresponding unit cylinder. The vector $\boldsymbol{\underline{X}}_t$ is then $S\alpha S$, and it is representable on $C_4^{\|\cdot\|}$ if and only if 
		\begin{equation}\label{asu:zeropoles}
			\Gamma_2 \Big(\{(0,-1),(0,+1)\}\Big) = 0.
		\end{equation}
	\end{prop}
	The representability condition \eqref{asu:zeropoles} appears in sharp contrast with Remark \ref{rem:dar1} and is also reminiscent of Remark \ref{rem:interp_repres}. 
	It intuitively means that the joint vector $\boldsymbol{\underline{X}}_t$ will admit a representation on the unit cylinder provided realisations $(\varepsilon_{1,t},\varepsilon_{2,t})$ where $\varepsilon_{2,t}$ is extreme and $\varepsilon_{1,t}$ is not occur with probability zero. 
	If this holds, then, intuitively, every jump in the trajectory of $(X_{2,t})$ necessarily coincides with a bubble peak in the trajectory of $(X_{1,t})$, and each incoming jump in the former is thus betrayed by the early build-up of a bubble in the latter.\footnote{Note that extreme realisations of $\varepsilon_{1,t}$ may nevertheless occur alongside non-extreme realisations of $\varepsilon_{2,t}$, as $\Gamma_2 \big(\{(-1,0),(+1,0)\}\big)$ can a priori be positive. Thus, intuitively, a bubble peak may be reached in $(X_{1,t})$ with no jump occurring in $(X_{2,t})$.}
	When considered univariately, $(X_{2,t})$ features sudden, unpredictable bursts --and is thus not-past-representable--, but this unpredictability appears to fade away when $(X_{2,t})$ is considered jointly with the <<informative>> process $(X_{1,t})$.
	The next proposition provides the tail conditional distribution of $\bX_{t+1}$ given (a large in norm) observation $\bX_t$, and shows that these heuristics are essentially correct. The anticipative component does inform about incoming jumps in the other component, and, quite surprisingly, the non-anticipative component also brings information about the anticipative one.
	For expository purposes, we distinguish several cases according to the conditioning event. 
	The proposition is followed by a detailed interpretation of each case.
	\begin{prop}
		\label{prop:biar1_asymp}
		Let $(\bX_t)$ as in \eqref{def:bivarAR1} and assume that \eqref{asu:zeropoles} holds.
		For $\eta_0>0$ and $\theta_0\in]-\pi,\pi]$, define
		$V_{0} = \{(\cos u, \sin u)\in S_2: \hspace{0.2cm} u\in[\theta_0-\eta_0,\theta_0+\eta_0]\}$ and let $B(V_0)=V_0\times\mathbb{R}^2$. Define also, for $\theta\in]-\pi,\pi]$, $\eta>0$, and $P$ any closed set of $\mathbb{R}^2$,\footnote{
		This ensures that $A_{\theta,\eta,P}=\{(0,0)\}\times P\hspace{0.3cm}+\hspace{0.3cm}\left\{ (\cos u, \sin u, 0, \rho_2 \sin u): \hspace{0.2cm} u \in [\theta-\eta,\theta+\eta]\right\}$ defines a proper Borel set, which could fail if $P$ was a general Borel set \cite{erdos}.
		One could assume more generally that $P$ is an $F_\sigma$ set, but a closed set will be enough for our purpose here.
		}
		$$
		A_{\theta,\eta,P} = \left\{ (\cos u, \sin u, 0, \rho_2 \sin u) + (0,0,x,y)\in C_4^{\|\cdot\|}: \hspace{0.2cm} u \in [\theta-\eta,\theta+\eta] \hspace{0.1cm} \text{and} \hspace{0.1cm} (x,y)\in P\right\},
		$$
		and $V_{\theta,\eta} = \{(\cos u, \sin u)\in S_2: \hspace{0.2cm} u\in[\theta-\eta,\theta+\eta]\}$.
		
		\noindent $(\iota)$ Assume that $V_0 \cap \{(\pm1,0),(0,\pm1)\}=\emptyset$. Then,
		\begin{align*}
			\mathbb{P}_x^{\|\cdot\|}\Big(\boldsymbol{\underline{X}}_t,A_{\theta,\eta,P}\Big|B(V_0)\Big) & \underset{x\rightarrow+\infty}{\longrightarrow} \dfrac{\Gamma_2\Big(V_{\theta,\eta}\cap V_0\Big)}{\Gamma_2\Big(V_0\Big)}\delta_{\{(0,0)\}}(P).
		\end{align*}

		\noindent $(\iota\iota)$ Assume $(0,\vartheta)\in V_0$, for some $\vartheta\in \{-1,+1\}$, and $V_0\cap \{\pm(1,0),(0,-\vartheta)\}=\emptyset$. Then, 
		\begin{align*}
			\mathbb{P}_x^{\|\cdot\|}\Big(\boldsymbol{\underline{X}}_t,A_{\theta,\eta,P}\Big|B(V_0)\Big) & \underset{x\rightarrow+\infty}{\longrightarrow} \dfrac{\dfrac{\sigma_2^\alpha}{2}\dfrac{|\rho_2|^\alpha}{1-|\rho_2|^\alpha}\delta_{\{(0,\vartheta)\}}(V_{\theta,\eta}) + \Gamma_2\Big(V_{\theta,\eta}\cap V_0\Big)}{\dfrac{\sigma_2^\alpha}{2}\dfrac{|\rho_2|^\alpha}{1-|\rho_2|^\alpha}+\Gamma_2\Big(V_0\Big)}\delta_{\{(0,0)\}}(P).
		\end{align*}
		
		\noindent $(\iota\iota\iota)$ Assume $(\vartheta,0)\in V_0$, for some $\vartheta\in \{-1,+1\}$, and $V_0\cap \{\pm(0,1),(-\vartheta,0)\}=\emptyset$. Then,
		\begin{align*}
			& \mathbb{P}_x^{\|\cdot\|}\Big(\boldsymbol{\underline{X}}_t,A_{\theta,\eta,P}\Big|B(V_0)\Big) \\
			& \\
			& \hspace{0.5cm} \underset{x\rightarrow+\infty}{\longrightarrow} \dfrac{\bigg(\dfrac{\sigma_1^\alpha}{2}\dfrac{ |\rho_1|^{2\alpha}}{1-|\rho_1|^\alpha}\delta_{\{0\}}(P_2) + \dfrac{|\rho_1|^{\alpha}}{2}{\sigma_1}_{|_{P_2}}^{\alpha} \bigg)\delta_{\{(\vartheta,0)\}}(V_{\theta,\eta})\delta_{\{\vartheta\rho_1^{-1}\}}(P_1) + \Gamma_2\Big(V_{\theta,\eta}\cap V_0\Big)\delta_{\{(0,0)\}}(P)}{\dfrac{\sigma_1^\alpha}{2}\dfrac{|\rho_1|^\alpha}{1-|\rho_1|^\alpha} + \Gamma_2(V_0)},
		\end{align*}
		where $P_1 = \{x:\hspace{0.1cm}(x,y)\in P\}$, $P_2 = \{y:\hspace{0.1cm}(x,y)\in P\}$ and ${\sigma_1}_{|_{P_2}}:=\bigg({\displaystyle\int_{S(P_2)}|s_1|^\alpha\Gamma_2(d\bs)}\bigg)^{1/\alpha}$ with $S(P_2):=\left\{\frac{(\rho_1^{-1},y)}{\sqrt{\rho_1^{-2}+y^2}}\in S_2:\hspace{0.2cm} y\in P_2\right\}$.
	\end{prop}
	
	\paragraph{Interpretation of Proposition \ref{prop:biar1_asymp}}\hfill
	
	\noindent In the spirit of this proposition, $V_0$ is typically a small neighbourhood on the unit sphere $S_2$ accounting for the observed realisation of $(X_{1,t},X_{2,t})/\sqrt{X_{1,t}^2+X_{2,t}^2}$, that is, the relative magnitudes of $X_{1,t}$ and $X_{2,t}$.\footnote{Recall that the results are always conditional on $\sqrt{X_{1,t}^2+X_{2,t}^2}$ being large: either $X_{1,t}$ is extreme, either $X_{2,t}$ is extreme, or both are.}
	The smaller the neighoubourhood $V_0$, the more <<accurately>> we assume to observe these relative magnitudes.
	The proposition considers three main scenarii: case $(\iota)$ $X_{1,t}$ and $X_{2,t}$ are of comparable magnitudes, case $(\iota\iota)$ $X_{2,t}$ is much larger -possibly infinitely larger- than $X_{1,t}$, and case $(\iota\iota\iota)$ $X_{1,t}$ is much larger -possibly infinitely larger- than $X_{2,t}$.
	Each of these three conditioning leads to different odds regarding the potential outcomes at $t+1$.\\
	
	\noindent \textbf{Case $\boldsymbol{(\iota)}$} To fix ideas, let us assume that $X_{1,t}$ and $X_{2,t}$ are observed to be of same signs and approximately of equal magnitudes, that is, $V_0$ is a small neighbourhood of $c(1,1)\in S_2$, with $c=2/\sqrt{2}$ (i.e., $\theta_0=\pi/4$ and $\eta_0>0$ small).
	Now, evaluating the tail conditional probability at $A_{\theta_0,\eta_0,P}$ for $P$ an arbitrarily small closed neighbourhood of $(0,0)$, for instance $P=[-\epsilon_1,\epsilon_1]\times [-\epsilon_2,\epsilon_2]$ for $\epsilon_1,\epsilon_2>0$ small, we obtain that
	$$
	\mathbb{P}_x^{\|\cdot\|}\Big(\boldsymbol{\underline{X}}_t,A_{\theta_0,\eta_0,P}\Big|B(V_0)\Big)  \underset{x\rightarrow+\infty}{\longrightarrow} 1.
	$$
	Intuitively during an extreme event, if $X_{1,t}$ and $X_{2,t}$ are observed to be of approximately equal magnitudes, then the vector $(X_{1,t},X_{2,t},X_{1,t+1},X_{2,t+1})/\sqrt{X_{1,t}^2+X_{2,t}^2}$ will belong with certainty to a small neighbourhood of $c(1,1,0,\rho_2)$.\footnote{The size of this neighbourhood will be commensurate to the accuracy of the observed relative magnitudes: for smaller $V_0$ (higher observation accuracy), smaller neighbourhoods around $c(1,1,0,\rho_2)$ will provide the same level of certainty. 
		For a fixed $V_0$, one can also evaluate the conditional probability over smaller neighbourhoods within $V_0$ by considering sets $A_{\theta,\eta,P}$ with $[\theta-\eta,\theta+\eta]\subset[\theta_0-\eta_0,\theta_0+\eta_0]$ for instance, which leads to the ratio in terms of $\Gamma_2$ as in the proposition. 
		Note that $P$ can be taken arbitrarily small regardless of $V_0$ without affecting the conditional probability, provided it contains $(0,0)$.}
	This straightforwardly extends to the case when $X_{1,t}$ and $X_{2,t}$ are of comparable magnitudes but not necessarily equal ones, i.e., when $V_0$ is instead a small neighbourhood of $(\cos \theta_0,\sin \theta_0)$. Then, $(X_{1,t},X_{2,t},X_{1,t+1},X_{2,t+1})/\sqrt{X_{1,t}^2+X_{2,t}^2}$ will belong with certainty to a small neighbourhood of $(\cos \theta_0,\sin\theta_0,0,\rho_2\sin\theta_0)$.\\
	This reveals that if at any date in time both series are simultaneously extreme, then, with certainty, the anticipative component will collapse at the immediately following date, and the non-anticipative component will decay by $\rho_2$.\\
	
	\noindent \textbf{Case $\boldsymbol{(\iota\iota)}$} Let us assume here that $V_0$ is an arbitrarily small neighbourhood of $(0,1)$, i.e., $X_{2,t}$ is observed positive and much larger in magnitude than $X_{1,t}$.
	Evaluating the conditional probability at $A_{\theta_0,\eta_0,P}$ with $P=[-\epsilon_1,\epsilon_1]\times [-\epsilon_2,\epsilon_2]$ an arbitrarily small neighbourhoods of $(0,0)$, we have that 
	$$
	\mathbb{P}_x^{\|\cdot\|}\Big(\boldsymbol{\underline{X}}_t,A_{\theta_0,\eta_0,P}\Big|B(V_0)\Big)  \underset{x\rightarrow+\infty}{\longrightarrow} 1.
	$$
	Thus, if during an extreme event, $X_{2,t}$ is observed to be much larger than $X_{1,t}$, then the vector $\boldsymbol{\underline{X}}_t/\|\boldsymbol{\underline{X}}_t\|$ will belong with certainty to a small neighbourhood of $(0,1,0,\rho_2)$.
	Observing at date $t$ an extreme in the non-anticipative component alongside a much smaller, possibly non-extreme value on the anticipative series indicates that at $t+1$, with certainty, the non-anticipative component will decay by $\rho_2$ whereas the anticipative component will remain small.\\

	\noindent \textbf{Case $\boldsymbol{(\iota\iota\iota)}$} Again to fix ideas, assume that $V_0$ is an arbitrarily small neighbourhood of $(1,0)$, i.e., $X_{1,t}$ is observed positive and much larger in magnitude than $X_{2,t}$. Contrary to $(\iota)$ and $(\iota\iota)$ where, practically, a single outcome captures all the probability mass, several clearly distinct potential outcomes share the likelihood in this case. \\
	Letting $A_{\theta_0,\eta_0,P}$ with $P=\Big([\rho_1^{-1}-\epsilon_1,\rho_1^{-1}+\epsilon_1]\cup[-\epsilon_1,\epsilon_1]\Big) \times \mathbb{R}$ for $\epsilon_1>0$ arbitrary small, we obtain after elementary computations that
	\begin{align*}
		\mathbb{P}_x^{\|\cdot\|}\Big(\boldsymbol{\underline{X}}_t,A_{\theta,\eta,P}\Big|B(V_0)\Big)\underset{x\rightarrow\infty}{\longrightarrow}
		1.
	\end{align*}
	Thus, the probability mass appears to be localised in a main region, which is a neighbourhood of the points $(1,0,0,z)$, $(1,0,\rho_1^{-1},z)$, $z\in\mathbb{R}$.
	Within this main region, the probability mass can be further localised into two distinct areas:
	\begin{align*}
		\mathbb{P}_x^{\|\cdot\|}\Big(\boldsymbol{\underline{X}}_t,A_{\theta,\eta,P}\Big|B(V_0)\Big)\underset{x\rightarrow\infty}{\longrightarrow}
		\left\{
		\begin{array}{cl}
			\hspace{-0.1cm} \dfrac{\Gamma_2(V_0)}{\dfrac{\sigma_1^\alpha}{2}\dfrac{|\rho_1|^\alpha}{1-|\rho_1|^\alpha} + \Gamma_2(V_0)},  & \hspace{0cm} \text{for} \hspace{0.3cm} P=[-\epsilon_1,\epsilon_1]\times[-\epsilon_2,\epsilon_2], \\
			& \\
			\dfrac{\dfrac{\sigma_1^\alpha}{2}\dfrac{|\rho_1|^\alpha}{1-|\rho_1|^\alpha}}{\dfrac{\sigma_1^\alpha}{2}\dfrac{|\rho_1|^\alpha}{1-|\rho_1|^\alpha} + \Gamma_2(V_0)}
			,  & \hspace{0cm} \text{for} \hspace{0.3cm} P=[\rho_1^{-1}-\epsilon_1,\rho_1^{-1}+\epsilon_1]\times \mathbb{R},
		\end{array}
		\right.
	\end{align*}
	for $\epsilon_1,\epsilon_2>0$ small.
	Given that the two areas have complementary probability masses, it appears that with certainty: $\boldsymbol{\underline{X}}_t/\|\boldsymbol{\underline{X}}_t\|$ will either $(1)$ belong to a small neighbourhood of $(1,0,0,0)$, or $(2)$ belong to a small neighbourhood of the points $(1,0,\rho_1^{-1},z)$, $z\in\mathbb{R}$. 
	The area corresponding to $(1)$ yields a straightforward interpretation:\\
	\indent $(1)$ The outcome $\Big\{\boldsymbol{\underline{X}}_t/\|\boldsymbol{\underline{X}}_t\|$ belongs to a small neighbourhood of $(1,0,0,0)\Big\}$ corresponds to an event in which the anticipative component is extreme at date $t$ and collapses at $t+1$ while the non-anticipative series is small both at $t$ and $t+1$. 
	The conditional likelihood of this outcome can be arbitrarily large or small according to how much weight $\Gamma_2$ charges on the neighbourhood $V_0$ of $(1,0)$.\\
	\indent (2) Contrary to the previous case, the probability mass on the area $(1,0,\rho_1^{-1},z)$, $z\in\mathbb{R}$ does not appear to be localised in an arbitrarily small neighbourhood but can in general be dispersed over all $z$ on the real line. 
	This family of events describes outcomes for which, from date $t$ to date $t+1$, the anticipative component increases by a factor $\rho_1^{-1}$ while the non-anticipative series remains either non-extreme or jumps to some extreme value.\\ 
	One can evaluate the probability mass of events corresponding to specific jumps sizes of the non-anticipative components. For any closed set $P_2\subset \mathbb{R}$, one has
	\begin{align*}
		\mathbb{P}_x^{\|\cdot\|}\Big(\boldsymbol{\underline{X}}_t,A_{\theta,\eta,P}\Big|B(V_0)\Big)\underset{x\rightarrow\infty}{\longrightarrow}
		\dfrac{\dfrac{\sigma_1^\alpha}{2}\dfrac{ |\rho_1|^{2\alpha}}{1-|\rho_1|^\alpha}\delta_{\{0\}}(P_2) + \dfrac{|\rho_1|^{\alpha}}{2}{\sigma_1}_{|_{P_2}}^{\alpha}}{\dfrac{\sigma_1^\alpha}{2}\dfrac{|\rho_1|^\alpha}{1-|\rho_1|^\alpha} + \Gamma_2(V_0)}, \hspace{0.3cm} \text{for} \hspace{0.2cm} P=[\rho_1^{-1}-\epsilon_1,\rho_1^{-1}+\epsilon_1]\times P_2.
	\end{align*}
	\begin{itemize}
		\item Taking for instance $P_2=[M,+\infty[$ (resp. $P_2=]-\infty,-M]\cup[M,+\infty[$), for some $M>0$, one gets the conditional likelihood of events $(1,0,\rho^{-1},z)$, for $|z|\ge M$, i.e., outcomes for which the anticipative component increases by a factor $\rho_1^{-1}$ and the non-anticipative jumps above the positive threshold $M$ (resp. outside the interval $]-M,M[$):
		\begin{align*}
			\mathbb{P}_x^{\|\cdot\|}\Big(\boldsymbol{\underline{X}}_t,A_{\theta,\eta,P}\Big|B(V_0)\Big)\underset{x\rightarrow\infty}{\longrightarrow}
			\dfrac{\dfrac{|\rho_1|^\alpha}{2}\displaystyle\int_{S(P_2)}|s_1|^\alpha\Gamma_2(d\bs)}{\dfrac{\sigma_1^\alpha}{2}\dfrac{|\rho_1|^\alpha}{1-|\rho_1|^\alpha} + \Gamma_2(V_0)},  & \hspace{0.5cm} \text{for} \hspace{0.2cm} P=[\rho^{-1}-\epsilon_1,\rho^{-1}+\epsilon_1]\times P_2,
		\end{align*}
		with $S(P_2):=\left\{\frac{(\rho_1^{-1},y)}{\sqrt{\rho_1^{-2}+y^2}}\in S_2:\hspace{0.2cm} y\in P_2\right\}$. 
		\item For $P_2=[-\epsilon_2,\epsilon_2]$ with $\epsilon_2>0$ small, one can gauge the conditional likelihood of the anticipative component increasing by a factor $\rho_1^{-1}$ from $t$ to $t+1$ while the non-anticipative component remains close to non-extreme.
		\begin{align*}
			\mathbb{P}_x^{\|\cdot\|}\Big(\boldsymbol{\underline{X}}_t,A_{\theta,\eta,P}\Big|B(V_0)\Big) & \underset{x\rightarrow\infty}{\longrightarrow}
			\dfrac{\dfrac{\sigma_1^\alpha}{2}\dfrac{ |\rho_1|^{2\alpha}}{1-|\rho_1|^\alpha} + \dfrac{|\rho_1|^\alpha}{2}\displaystyle\int_{S(P_2)}|s_1|^\alpha\Gamma_2(d\bs)}{\dfrac{\sigma_1^\alpha}{2}\dfrac{|\rho_1|^\alpha}{1-|\rho_1|^\alpha} + \Gamma_2((1,0))}
		\end{align*}
		for $P=[\rho^{-1}-\epsilon_1,\rho^{-1}+\epsilon_1]\times [-\epsilon_2,\epsilon_2]$.
	\end{itemize}
	\begin{rem}
	\rm Further insights can be drawn from case $(\iota\iota\iota)$ if $\Gamma_2((\pm1,0))=0$, i.e., if realisations $(\varepsilon_{1,t},\varepsilon_{2,t})$ where $\varepsilon_{1,t}$ is extreme and $\varepsilon_{2,t}$ is not almost never occur. Then, $\Gamma_2(V_0)$ becomes arbitrarily close to $0$ for $V_0$ an arbitrarily small neighbourhood of $(\vartheta,0)$.\footnote{This holds because $\Gamma_2$ is a finite measure.} 
	In that case, neglecting the difference and assuming  $\Gamma_2(V_0)=\Gamma_2((\vartheta,0))=0$, we have
		\begin{align*}
		\mathbb{P}_x^{\|\cdot\|}\Big(\boldsymbol{\underline{X}}_t,A_{\theta,\eta,P}\Big|B(V_0)\Big) \underset{x\rightarrow\infty}{\longrightarrow}
		\left\{
		\begin{array}{cl}
			\hspace{-0.1cm} 0,  & \hspace{0cm} \text{for} \hspace{0.3cm} P=[-\epsilon_1,\epsilon_1]\times[-\epsilon_2,\epsilon_2], \\
			& \\
			1,  & \hspace{0cm} \text{for} \hspace{0.3cm} P=[\rho_1^{-1}-\epsilon_1,\rho_1^{-1}+\epsilon_1]\times \mathbb{R},
		\end{array}
		\right.
	\end{align*}
	which indicates that a bubble in $(X_{1,t})$ necessarily reaches its peak at a jump date in $(X_{2,t})$. The bubble peak is always <<signaled>>.
	Observing $X_{1,t}$ extreme and $X_{2,t}$ non-extreme thus implies that the bubble will last at least one more period.\\
	Taking now $P_2=]-\infty,M]\cup[M,+\infty[$ for $M>0$ (resp. $P_2=[-\epsilon_2,\epsilon_2]$ for $\epsilon_2>0$) arbitrarily small, notice that the integral $\int_{S(P_2)}|\sigma_1|^\alpha\Gamma_2(d\bs)$ can be made arbitrarily close to $\sigma_1^\alpha-\Gamma_2((\pm1,0))$ (resp. $\Gamma_2((\pm1,0))$ ). Again, assuming that $\Gamma_2((\pm1,0))=0$ and neglecting the difference, this yields that
	\begin{align*}
		\mathbb{P}_x^{\|\cdot\|}\Big(\boldsymbol{\underline{X}}_t,A_{\theta,\eta,P}\Big|B(V_0)\Big) \underset{x\rightarrow\infty}{\longrightarrow}
		\left\{
		\begin{array}{cl}
			\hspace{-0.1cm} 1-|\rho_1|^\alpha,  & \hspace{0cm} \text{for} \hspace{0.3cm} P=[\rho_1^{-1}-\epsilon_1,\rho_1^{-1}+\epsilon_1]\times(]-\infty,M]\cup[M,+\infty[), \\
			& \\
			|\rho_1|^\alpha,  & \hspace{0cm} \text{for} \hspace{0.3cm} P=[\rho_1^{-1}-\epsilon_1,\rho_1^{-1}+\epsilon_1]\times [-\epsilon_2,\epsilon_2],
		\end{array}
		\right.
	\end{align*}
	for $M>0$ and $\epsilon_2>0$ arbitrarily small. One recognises the probability of a bubble surviving one or more period, and its complementary, in the univariate anticipative AR(1) model.\footnote{Here, it would be more accurate to speak about the probability of a bubble in $(X_{1,t})$ surviving at least two more periods, given that it will survive at least one more. }
	\end{rem}

	\noindent Table \ref{tab:summary1} summarises the potential outcomes of each specific conditioning event and illustrates the typical profile of the trajectory of $(\bX_t)$ in each case.
	
	\begin{table}[h!]
	\vspace*{-1cm}
	    \hspace*{-0.2cm}\begin{tabular}{| >{\centering\arraybackslash} m{5cm} | >{\centering\arraybackslash} m{4cm} | >{\centering\arraybackslash} m{2cm}|  >{\centering\arraybackslash} m{5.5cm}| } 
	     \hline\hline
	      Observation $\frac{(X_{1,t},X_{2,t})}{\sqrt{X_{1,t}^2+X_{2,t}^2}}\in V_0$ & Potential outcomes (neighbourhood of) & Conditional probability & Trajectorial interpretation \\
	      \hline
	       \begin{tabular}{c}\includegraphics[scale=0.55]{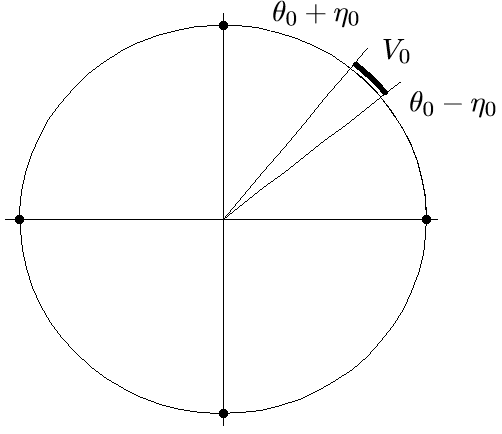}\\
	       $X_{1,t}$ and $X_{2,t}$ both extreme
	       \end{tabular}  &  \begin{tabular}{c}
	             $(\cos u,\sin u, 0, \rho_2 \sin u)$  \\
	             $u\in[\theta_0-\eta_0,\theta_0+\eta_0]$ 
	         \end{tabular} 
	         & 1 & \begin{tabular}{c}\includegraphics[scale=0.5]{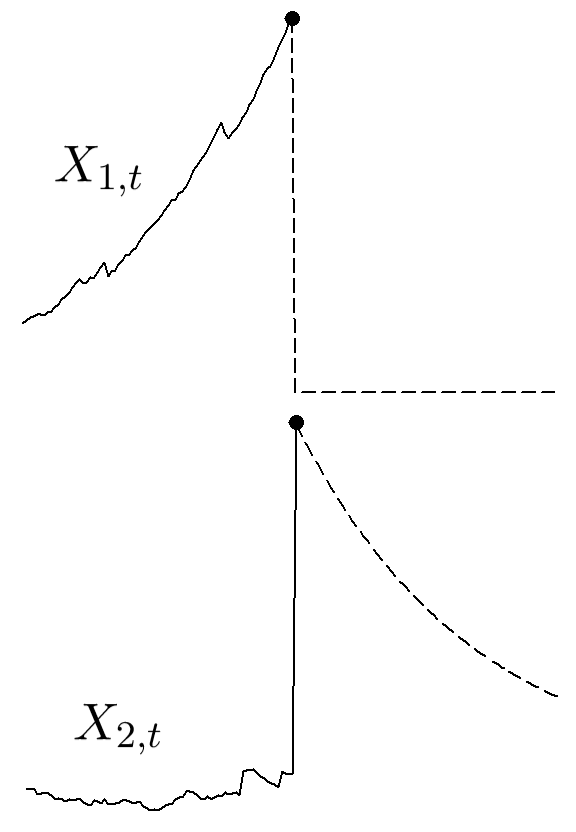}\\Bubble peak signaled by jump \end{tabular}\\
	         \hline
	         \begin{tabular}{c}
	        \includegraphics[scale=0.55]{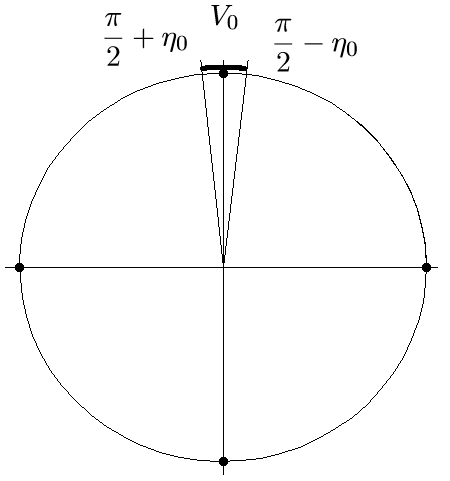}\\
	       $|X_{1,t}|<<|X_{2,t}|$
	       \end{tabular} & $(0,1, 0, \rho_2)$  & 1 & \begin{tabular}{c}\includegraphics[scale=0.5]{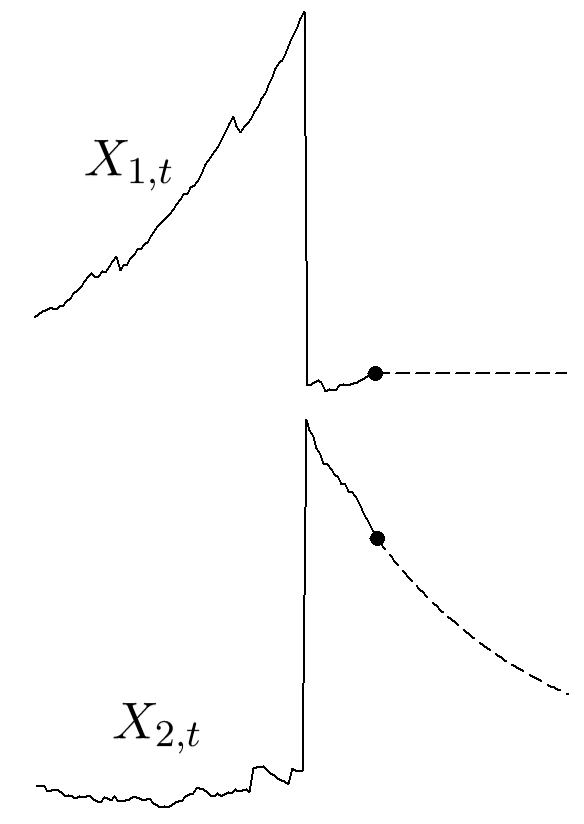}\\Post crash, jump decaying \end{tabular}\\
	       \hline\hline
	       \end{tabular}
	       \caption[Future paths odds and trajectorial interpretation of the bivariate stable process \eqref{def:bivarAR1}]{For each case considered in Proposition \ref{prop:biar1_asymp} (first column), the potential outcomes for $\boldsymbol{\underline{X}}_t/\|\boldsymbol{\underline{X}}_t\|$ is provided (second column), alongside the asymptotic conditional probability mass over the corresponding outcomes (third column). 
	       Each case can be related to  specific events in the trajectory of $(\bX_t)$, which are illustrated and labeled in the last column.
	       The solid lines represent past trajectories and the present dates are symbolised by points.
	       In the outcome $(1,0,\rho_1^{-1},z)$, $z\in\mathbb{R}$, the bubble survives at least one more period, but could survive more. Also, when the peak will be reached, a jump of \textit{a priori} any size (including zero) may occur and then decay.
	       Multiple potential paths are thus represented in dashed lines oriented by arrows, and the grey shaded area symbolises the jump size distribution. 
	       Here: $\omega_1=\Gamma_2(V_0)$ and  $\omega_2=\sigma_1^\alpha|\rho_1|^\alpha/\Big(2(1-|\rho_1|^\alpha)\Big)$.
	       }
	       \label{tab:summary1}
	       \end{table}
	\begin{table}[h!]
    \hspace*{-0.2cm}\begin{tabular}{| >{\centering\arraybackslash} m{5cm} | >{\centering\arraybackslash} m{4cm} | >{\centering\arraybackslash} m{2cm}|  >{\centering\arraybackslash} m{5.5cm}| } 
     \hline\hline
       Observation $\frac{(X_{1,t},X_{2,t})}{\sqrt{X_{1,t}^2+X_{2,t}^2}}\in V_0$ & Potential outcomes (neighbourhood of) & Conditional probability & Trajectorial interpretation \\
      \hline
          \multirow{2}{*}{
   \begin{tabular}{c} \\
   \\
   \\
   \\
   \\
   \includegraphics[scale=0.55]{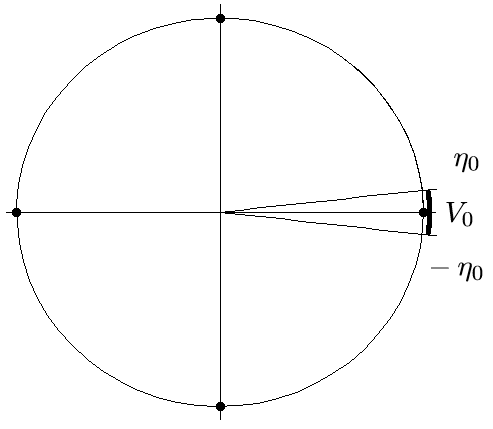}\\ $|X_{1,t}|>>|X_{2,t}|$
	       \end{tabular}} & $(1,0, 0, 0)$  & $\dfrac{\omega_1}{\omega_1+\omega_2}$ &  \begin{tabular}{c}\vspace{0.1cm}\includegraphics[scale=0.6]{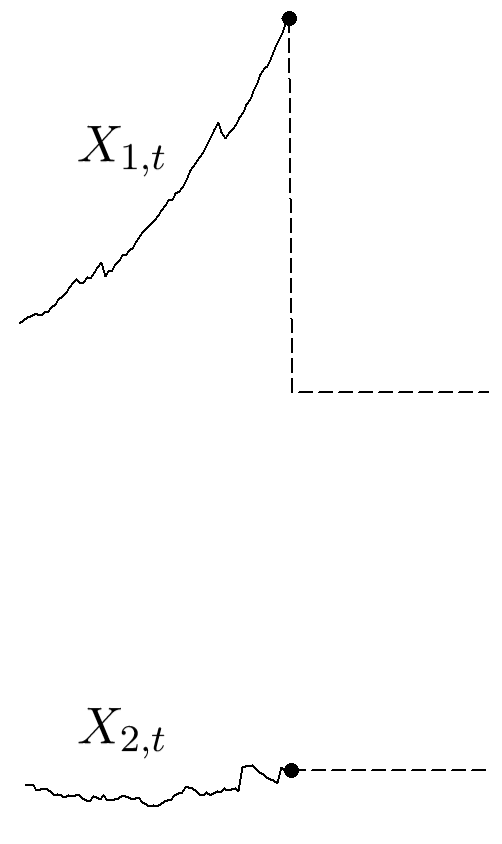}\\ Bubble peak, no jump signal\end{tabular}\vspace{0.1cm}\\ \cline{2-4}
	       & \begin{tabular}{c} $(1,0, \rho_1^{-1}, z)$\\ for $z\in\mathbb{R}$
	       \end{tabular}  & $\dfrac{\omega_2}{\omega_1+\omega_2}$ & \vspace{0.1cm} \begin{tabular}{c}\includegraphics[scale=0.6]{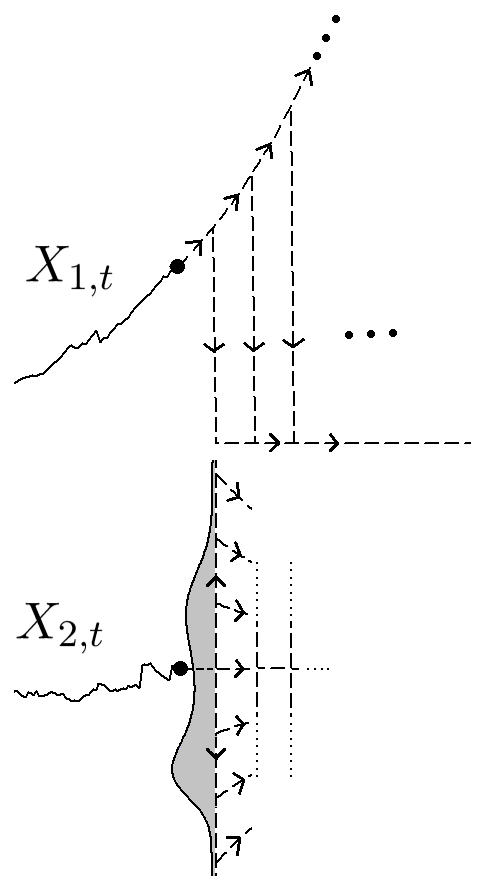}\\ Pre-peak bubble inflation,\\ potential peak/jump at $t+1$ \end{tabular}\\
	      \hline\hline
   	    \end{tabular}
	    \label{tab:summary2}
	\end{table}
	
	\clearpage

	\section{Concluding remarks and perspectives}
	\label{sec:conclude}

	\subsection{Conclusion}
	
\begin{figure}[h!]
    \centering
    \includegraphics[scale=0.5]{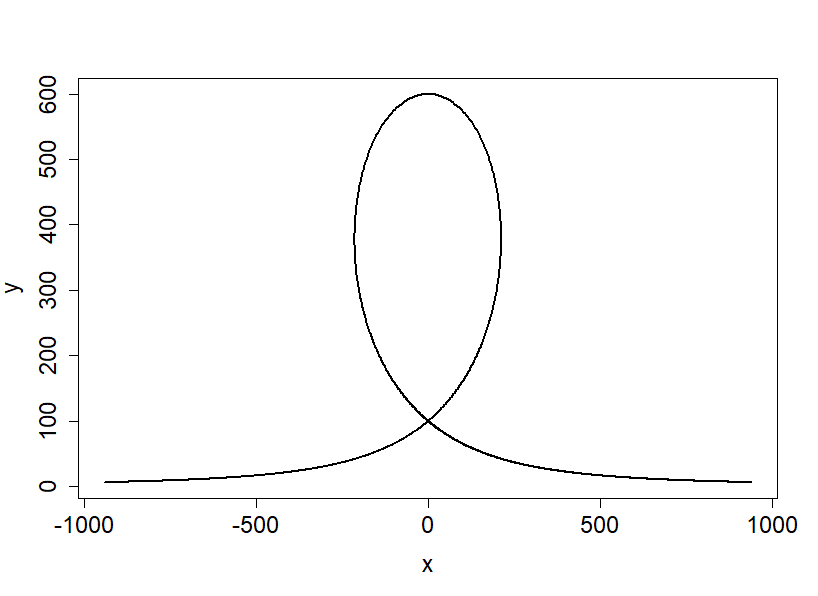}
    \caption{Horizontal strophoid for $a=100$, $b=5$.}
    \label{fig:my_strophoid}
	\end{figure}

\begin{figure}[h!]
		\centering
		\includegraphics[scale=1]{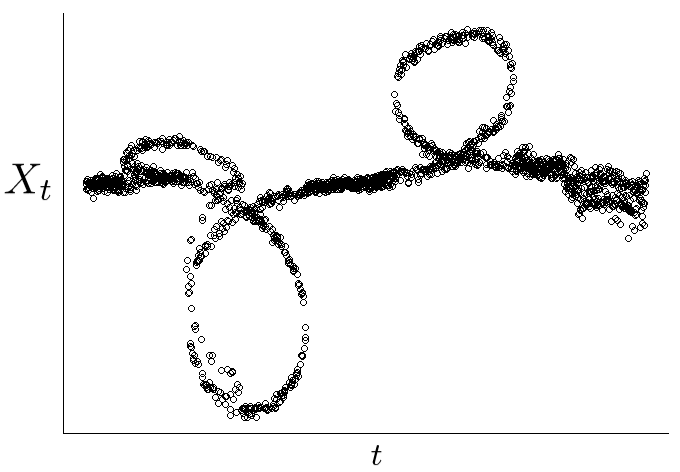}
		\caption{Sample path of a linear $1$-stable process featuring strophoidal patterns ($a=100$, $b=5$).}
		\label{fig:strophoid_traj}
	\end{figure}

    For $\alpha$-stable infinite moving averages and processes resulting from the linear combination thereof, the conditional distribution of future paths given the observed past trajectory during extreme events is obtained on the basis of a new spectral representation of stable random vectors on unit cylinders relative to semi-norms.
	Contrary to the case of norms, such representations yield a multivariate regularly varying tails property that is appropriate for prediction purposes, however not all stable random vectors can be represented on semi-norm unit cylinders. A characterisation is provided and finite length paths of $\alpha$-stable moving averages and aggregates, which are themselves multivariate $\alpha$-stable, are embedded into this framework.
	It is shown that paths of such a process admit semi-norm representations that are appropriate for prediction purposes if and only if the process is <<anticipative enough>>.\\
	\indent In this framework, our approach reveals that instead of their attractive <<causal>> interpretation, non-anticipative processes appear to rather presume, by construction, the unpredictability of extreme events. Anticipative processes however, instead of <<depending on the future>>, rather assume that future events feature early visible signs betraying their incoming occurrences. 
	These early signs take the form of emerging trends and patterns that an observer can identify and use to infer about future potential outcomes.
	Whether extreme events in some time series data feature early visible signs or not is arguably an intrinsic property of the natural phenomenon being measured rather than one of the modelling. 
	One can nevertheless see that enforcing a non-anticipative process on any given time series data mechanically leads to a model which assumes that extremes are not inferrable beyond their unconditional likelihood of occurrence.	It appears in addition that modelling a time series by a single, say, AR recursive equation, entails assuming that the considered series is determined by a single pattern appearing recurrently through time.\footnote{At least in the heavy-tailed framework. In lighter-tailed frameworks, patterns are more weakly observed, if at all, and the dynamics is dominated by the persistence of the past trajectory \cite{rootzen}.}\\
	\indent In the univariate framework, processes resulting from linearly aggregating anticipative processes thus circumvent two implicitly <<built-in>> limitations of classical time series modelling, at least from a probabilistic standpoint. 
	It is furthermore often argued that linear processes suffer intrinsic limitations with regards to their dynamics and the type of patterns they can capture or reproduce.
	Proposition \ref{prop:relativprob} and its subsequent remarks however show that not only are linear processes actually able to generate trajectories featuring any number of any kind of patterns through time, by the tuning of $J$ and of the sequences $(d_{j,k})$ upon which only very mild assumptions are imposed, but that their conditional dynamics is moreover tractable.
	Future developments could even extend the notion of stable aggregates from the linear combination of a finite number of moving averages to a countable or a continuum of moving averages, and to moving averages which coefficients are themselves stochastic.
	If the linearity assumption surely entails certain dynamical restrictions, the pattern-complexity of trajectories cannot be counted among these weaknesses.
	Before outlining the perspectives for future work, we provide as an illustration of this flexibility a linear process exhibiting strophoidal, looping-like patterns.\\
	\\~
	Consider for $a,b$ positive real numbers the horizontal strophoid $\mathcal{S}=\{(x(t),y(t))\in\mathbb{R}^2: \hspace{0.2cm}t\in\mathbb{R}\}$, where for any $t\in\mathbb{R}$,
	\begin{align*}
	x(t) & = -at\dfrac{b-t^2}{1+t^2},  & y(t) & = \dfrac{a(b+1)}{1+t^2}.
	\end{align*}
	Figure \ref{fig:my_strophoid} provides an illustration of the horizontal strophoid for $a=100$, $b=5$.
	Letting for any $(x,y)\in\mathbb{R}^2$
	\begin{align*}
	\Pi_x(y):=y^3 -a(b+3)y^2+(x^2+a^2(2b+3))y-a^3(b+1),
	\end{align*}
	a Cartesian equation of the locus of the strophoid is given by $\Pi_x(y)=0$. 
	Construct now a \textit{non-random} sequence $(d_k)$ in the following way: for a given $k\in\mathbb{Z}$, draw an element uniformly at random in the set $\{y\in\mathbb{R}:\hspace{0.2cm}\Pi_k(y)=0\}$ --which may contain either one, two or three elements-- and assign it to $d_k$.
	Define then the process $X_t = \sum_{k\in\mathbb{Z}}d_k\varepsilon_{t+k}$ for $(\varepsilon_t)$ an i.i.d. $\alpha$-stable sequence, $1/2<\alpha<2$.\footnote{An elementary analysis shows that $\lim\limits_{t\rightarrow t_1}|x(t)|=+\infty$, $\lim\limits_{t\rightarrow t_2}y(t)=0$ if and only if $t_1=\pm\infty$ and $t_2=\pm\infty$, and that $y(t)/x^2(t)\rightarrow (b+1)/a$ for $t\rightarrow\pm\infty$.
   	Thus, $d_k\underset{|k|\rightarrow\infty}{\sim} \text{const} \hspace{0.1cm} k^{-2}$ and $(X_t)$ is well defined for $1/2<\alpha<2$.}
	It can be checked that the process $(X_t)$ is $(m,h)$-past-representable for any $m\ge0$, $h\ge1$. 
	Proposition \ref{prop:relativprob} applied to $\bXt=(X_{t-m},\ldots,X_{t+h})$ with $V=S_{m+1}$ shows that $\bXt/\|\bXt\|$ is asymptotically of the form $\pm\bd_{k_0}/\|\bd_{k_0}\|$, for some $k_0\in\mathbb{Z}$. 
	Given the construction of $(d_k)$, we deduce that the linear process $(X_t)$ features looping-like patterns in its trajectories, as depicted on Figure \ref{fig:strophoid_traj} for the choice of parameters $a=100$ and $b=5$.
	Its thorough analysis is left for further research.
	
	\subsection{Perspectives}
	
	Future lines of research could focus on several opened questions, both on probabilistic and statistical aspects. 
	\begin{enumerate}
	    \item For practical use in applications, estimation/learning methods need to be developed to infer the structure of <<a best approximating>> aggregated moving average to some time series data. 
	    This requires identifying the coefficients sequences of the moving averages involved in the aggregation (which characterise the shape of patterns appearing during extreme events), the number of moving averages involved (the number of different patterns) and the distribution parameters of the i.i.d. errors driving the process.
	    Recovering the patterns amounts in the general case to estimating an infinite number of parameters and one is likely to seek instead a parsimonious, low-dimensional structure in the coefficients sequence.
	    This structure need not assume that the coefficients of the moving averages satisfy some linear recursive relation, as in the case of ARMA processes. 
	    Complex patterns can be achieved even with few parameters, as illustrated by the strophoid-generating process defined in the previous section, which coefficients sequence is characterised by only two parameters.
	    \item The conditional distribution of aggregated moving averages holds asymptotically in the (semi)norm of the observed trajectory being large. 
	    Could we evaluate the approximation error made when using the asymptotic distribution in lieu of the finite distance one for prediction? What could be said about the conditional distribution when the process is close to its central values? A perhaps dual question would be: could we evaluate how variable the future trajectory may be around the predicted deterministic paths ? One could expect in particular that this variability may increase with the prediction horizon.
	    \item The conditional distribution obtained obviously requires to provide a conditioning Borel set, which represents the information about the shape of the observed trajectory.
	    This piece of observed trajectory can be viewed as a single noisy realisation of a piece of pattern generated by the process, which leaves room for uncertainty in the identification of that pattern.
	    The choice of an adequate conditioning Borel should thus rely on a statistical method. One could envision for instance tests of hypotheses to determine whether the observed trajectory is <<more akin>> a certain pattern 1 or another pattern 2, or to any other pattern of a certain collection.
	    Moreover, the length of the observed trajectory does not have to remain fixed: shorter observation windows closer to the present date may contain more <<up-to-date>> information, less influenced by now vanished past extreme events on the one hand, but on the other hand could also be more subject to noise and make the pattern identification more difficult.
	    Conversely, wider observation windows may provide more robust pattern identification but may also incorporate biased information, being influenced by now irrelevant past events.
	    One could envision looking for event-driven optimal window length based on a bias-variance trade-off.
	    \item A multivariate extension is illustrated on a simple bivariate process with one purely anticipative and one purely non-anticipative component.
	    New properties already emerge, such as the fact that while univariate non-anticipative processes never induce paths representable on unit cylinders, their paths may nonetheless be representable in higher dimensions when considered alongside an <<informative>> anticipative process.
	    Both components are more predictable when considered jointly rather than univariately.
	    The general multivariate framework can be readily embedded into this framework but numerous potential interactions between univariate components render the task challenging.
	\end{enumerate}

	\clearpage
	
	\section{Postponed proofs}
	\label{sec:proofs}
	
	\subsection{Proof of Proposition \ref{prop238}}
	
	Consider first the case where either $\alpha\ne1$ or $\bX$ is S1S. We only provide the proof for $\alpha\ne1$ as it is similar under both assumptions.\\
	Assume that $\Gamma(K^{\|\cdot\|})=0$ and let us show that $\bX$ admits a representation of the unit cylinder $C_d^{\|\cdot\|}$ relative to the semi-norm $\|\cdot\|$.
	The characteristic function of $\bX$ writes for any $\bu\in\mathbb{R}^d$, with $a=\tg(\pi\alpha/2)$,
	\begin{align*}
		\varphi_{\bX}(\bu) & = \exp\Bigg\{-\int_{S_d}\Big(|\scal{\bu,\bs}|^\alpha-ia(\scal{\bu,\bs})^{<\alpha>}\Big)\Gamma(d\bs)+i\,\scal{\bu,\bmu^0}\Bigg\}\\
		& = \exp\Bigg\{-\int_{S_d\setminus K^{\|\cdot\|}}\Big(|\scal{\bu,\bs}|^\alpha-ia(\scal{\bu,\bs})^{<\alpha>}\Big)\Gamma(d\bs)+i\,\scal{\bu,\bmu^0}\Bigg\}\\
		& = \exp\Bigg\{-\int_{S_d\setminus K^{\|\cdot\|}}\bigg(|\scal{\bu,\dfrac{\bs}{\|\bs\|}}|^\alpha-ia(\scal{\bu,\dfrac{\bs}{\|\bs\|}})^{<\alpha>}\bigg)\|\bs\|^{\alpha}\Gamma(d\bs)+i\,\scal{\bu,\bmu^0}\Bigg\}\\
		& = \exp\Bigg\{-\int_{T_{\|\cdot\|}(S_d\setminus K^{\|\cdot\|})}\bigg(|\scal{\bu,\bs'}|^\alpha-ia(\scal{\bu,\bs'})^{<\alpha>}\bigg)\left\|\dfrac{\bs'}{\|\bs'\|_e}\right\|^{\alpha}\Gamma\circ T_{\|\cdot\|}^{-1}(d\bs')+i\,\scal{\bu,\bmu^0}\Bigg\}\\
		& = \exp\Bigg\{-\int_{C_d^{\|\cdot\|}}\bigg(|\scal{\bu,\bs}|^\alpha-ia(\scal{\bu,\bs})^{<\alpha>}\bigg)\underbrace{\|\bs\|_e^{-\alpha}\Gamma\circ T_{\|\cdot\|}^{-1}(d\bs)}_{\Gamma^{\|\cdot\|}(d\bs)}+i\,\scal{\bu,\bmu^0}\Bigg\}
	\end{align*}
	where we used the change of variable $\bs' = T_{\|\cdot\|}(\bs)=\bs/\|\bs\|$ between the third and fourth lines, which yields the representation on $X_d^{\|\cdot\|}$.\\
	Reciprocally, assume that $\bX$ is representable on $C_d^{\|\cdot\|}$.
	By definition of the representability of $\bX$ on $C_d^{\|\cdot\|}$, there exists a measure $\gamma^{\|\cdot\|}$ on $C_d^{\|\cdot\|}$ and a non-random vector $\bm_{{\|\cdot\|}}^0\in\mathbb{R}^d$ such that
	\begin{align*}
	\varphi_{\bX}(\bu) & = \exp\Bigg\{-\int_{C_d^{\|\cdot\|}}\Big(|\scal{\bu,\bs}|^\alpha-ia(\scal{\bu,\bs})^{<\alpha>}\Big)\gamma^{\|\cdot\|}(d\bs)+i\,\scal{\bu,\bm_{{\|\cdot\|}}^0}\Bigg\}.
	\end{align*}
	With the change of variable $\bs' = T_{\|\cdot\|}^{-1}(\bs)=\bs/\|\bs\|_e$,
	\begin{align*}
	\varphi_{\bX}(\bu) & = \exp\Bigg\{-\int_{C_d^{\|\cdot\|}}\Big(|\scal{\bu,\dfrac{\bs}{\|\bs\|_e}}|^\alpha-ia(\scal{\bu,\dfrac{\bs}{\|\bs\|_e}})^{<\alpha>}\Big)\|\bs\|_e^{\alpha}\gamma^{\|\cdot\|}(d\bs)+i\,\scal{\bu,\bm_{{\|\cdot\|}}^0}\Bigg\}\\
	& = \exp\Bigg\{-\int_{T_{\|\cdot\|}^{-1}(C_d^{\|\cdot\|})}\Big(|\scal{\bu,\bs'}|^\alpha-ia(\scal{\bu,\bs'})^{<\alpha>}\Big)\left\|\dfrac{\bs'}{\|\bs'\|}\right\|_e^{\alpha}\gamma^{\|\cdot\|}\circ T_{\|\cdot\|}(d\bs')+i\,\scal{\bu,\bm_{{\|\cdot\|}}^0}\Bigg\}\\
	& = \exp\Bigg\{-\int_{S_d\setminus K^{\|\cdot\|}}\Big(|\scal{\bu,\bs}|^\alpha-ia(\scal{\bu,\bs})^{<\alpha>}\Big)\|\bs\|^{-\alpha}\gamma^{\|\cdot\|}\circ T_{\|\cdot\|}(d\bs)+i\,\scal{\bu,\bm_{{\|\cdot\|}}^0}\Bigg\}\\
	& = \exp\Bigg\{-\int_{S_d\setminus K^{\|\cdot\|}}\Big(|\scal{\bu,\bs}|^\alpha-ia(\scal{\bu,\bs})^{<\alpha>}\Big)\gamma(d\bs)+i\,\scal{\bu,\bm_{{\|\cdot\|}}^0}\Bigg\},
	\end{align*}
	where $\gamma(d\bs):=\|\bs\|^{-\alpha}\gamma^{\|\cdot\|}\circ T_{\|\cdot\|}(d\bs)$.
	Letting now $\overline{\gamma}(A):=\gamma(A\cap (S_d\setminus K^{\|\cdot\|}))$ for any Borel set $A$ of $S_d$, we have
	\begin{align*}
	\varphi_{\bX}(\bu) & = \exp\Bigg\{-\int_{S_d}\Big(|\scal{\bu,\bs}|^\alpha-ia(\scal{\bu,\bs})^{<\alpha>}\Big)\overline{\gamma}(d\bs)+i\,\scal{\bu,\bm_{{\|\cdot\|}}^0}\Bigg\}.
	\end{align*}
	By the unicity of the spectral representation of $\bX$ on $S_d$, we necessarily have $(\Gamma,\bmu^0)=(\overline{\gamma},\bm_{{\|\cdot\|}}^0)$. 
	Thus, $\overline{\gamma}$ and $\Gamma$ have to coincide, and in particular
	\begin{align*}
		\Gamma(K^{\|\cdot\|}) = \overline{\gamma}(K^{\|\cdot\|}) = \gamma(K^{\|\cdot\|}\cap (S_d\setminus K^{\|\cdot\|})) = \gamma(\emptyset)=0.
	\end{align*}
	Given that $\Gamma=\overline{\gamma}$ and $\Gamma(K^{\|\cdot\|}) = 0$, we can follow the initial steps of the proof to show that
	$\gamma^{\|\cdot\|}=\Gamma^{\|\cdot\|}$.\\
	
	\noindent Consider now the case where $\alpha=1$ and $\bX$ is not symmetric.
	Assume first that $\int_{S_d}\Big|\ln\|\bs\|\hspace{0.03cm}\Big|\Gamma(d\bs)<+\infty$, that is, $\Gamma(K^{\|\cdot\|})=0$ and $\int_{S_d\setminus K^{\|\cdot\|}}\Big|\ln\|\bs\|\hspace{0.03cm}\Big|\Gamma(d\bs)<+\infty$.
	With $a = 2/\pi$,
	\begin{align*}
		\varphi_{\bX}(\bu) & = \exp\Bigg\{-\int_{S_d}\Big(|\scal{\bu,\bs}|+ia\scal{\bu,\bs}\ln|\scal{\bu,\bs}|\Big)\Gamma(d\bs)+i\,\scal{\bu,\bmu^0}\Bigg\}\\
		& = \exp\Bigg\{-\int_{S_d\setminus K^{\|\cdot\|}}\bigg(|\scal{\bu,\dfrac{\bs}{\|\bs\|}}|+ia\scal{\bu,\dfrac{\bs}{\|\bs\|}}\ln|\scal{\bu,\dfrac{\bs}{\|\bs\|}}|\bigg)\|\bs\|\Gamma(d\bs)\\
		& \hspace{5cm} +i\,\scal{\bu,\bmu^0}-ia\int_{S_d\setminus K^{\|\cdot\|}}\scal{\bu,\bs}\ln\|\bs\|\Gamma(d\bs)\Bigg\}.
	\end{align*}
	We have $\int_{S_d\setminus K^{\|\cdot\|}}\scal{\bu,\bs}\ln\|\bs\|\Gamma(d\bs)=\sum_{i=1}^du_i\int_{S_d\setminus K^{\|\cdot\|}}s_i\ln\|\bs\|\Gamma(d\bs)=\scal{\bu,\tilde{\bmu}}$, and thus, 
	$$
	i\,\scal{\bu,\bmu^0}-ia\int_{S_d\setminus K^{\|\cdot\|}}\scal{\bu,\bs}\ln\|\bs\|\Gamma(d\bs) = i\scal{\bu,\bmu_{\|\cdot\|}^0}.
	$$
	The condition $ \int_{S_d\setminus K^{\|\cdot\|}}\Big|\ln\|\bs\|\Big|\Gamma(d\bs) < +\infty$, ensures that $|\bmu_{\|\cdot\|}^0|<+\infty$.
	Again with the change of variable $\bs' = T_{\|\cdot\|}(\bs)=\bs/\|\bs\|$, we get
	\begin{align*}
	\varphi_{\bX}(\bu) & = \exp\Bigg\{-\int_{T_{\|\cdot\|}(S_d\setminus K^{\|\cdot\|})}\bigg(|\scal{\bu,\bs'}|+ia\scal{\bu,\bs'}\ln|\scal{\bu,\bs'}|\bigg)\left\|\dfrac{\bs'}{\|\bs'\|_e}\right\|^{\alpha}\Gamma\circ T_{\|\cdot\|}^{-1}(d\bs') + i\scal{\bu,\bmu_{\|\cdot\|}^0}\Bigg\}\\
	& = \exp\Bigg\{-\int_{C_d^{\|\cdot\|}}\bigg(|\scal{\bu,\bs}|+ia\scal{\bu,\bs}\ln|\scal{\bu,\bs}|\bigg)\underbrace{\|\bs\|_e^{-\alpha}\Gamma\circ T_{\|\cdot\|}^{-1}(d\bs)}_{\Gamma^{\|\cdot\|}(d\bs)} + i\scal{\bu,\bmu_{\|\cdot\|}^0}\Bigg\}
	\end{align*}
	
    Reciprocally, assume there exists a measure $\gamma^{\|\cdot\|}$ on $C_d^{\|\cdot\|}$ satisfying \eqref{asu:integ_complique} and a non-random vector $\bm_{{\|\cdot\|}}^0\in\mathbb{R}^d$ such that 
	\begin{align*}
	\varphi_{\bX}(\bu) & = \exp\Bigg\{-\int_{C_d^{\|\cdot\|}}\Big(|\scal{\bu,\bs}|+ia\scal{\bu,\bs}\ln|\scal{\bu,\bs}|\Big)\gamma^{\|\cdot\|}(d\bs)+i\,\scal{\bu,\bm_{{\|\cdot\|}}^0}\Bigg\}.
	\end{align*}
	First, we can see that
	\begin{align*}
		\varphi_{\bX}(\bu) & = \exp\Bigg\{-\int_{C_d^{\|\cdot\|}}\bigg[\Big(|\scal{\bu,\dfrac{\bs}{\|\bs\|_e}}|+ia\scal{\bu,\dfrac{\bs}{\|\bs\|_e}}\ln|\scal{\bu,\dfrac{\bs}{\|\bs\|_e}}|\Big)\|\bs\|_e+ ia\scal{\bu,\bs}\ln\|\bs\|_e\bigg]\gamma^{\|\cdot\|}(d\bs)\\
		& \hspace{12cm}  +i\,\scal{\bu,\bm_{{\|\cdot\|}}^0}\Bigg\}.
	\end{align*}
	We will later show the following result:
	 \begin{lem}\label{le:a1_est_complique}
    Let $\gamma^{\|\cdot\|}$ a Borel measure on $C_d^{\|\cdot\|}$ satisfying \eqref{asu:integ_complique}. Then, 
    \begin{align}
     \int_{C_d^{\|\cdot\|}}\|\bs\|_e\Big|\ln\|\bs\|_e\hspace{0.03cm}\Big|\gamma^{\|\cdot\|}(d\bs)<+\infty.  \label{eq:complique2}
    \end{align}
    \end{lem}
	Assuming Lemma \ref{le:a1_est_complique} holds, then by the Cauchy-Schwarz inequality, we have $\int_{C_d^{\|\cdot\|}}|\scal{\bu,\bs}|\Big|\ln\|\bs\|_e\hspace{0.03cm}\Big|\gamma^{\|\cdot\|}(d\bs)<+\infty$, and thus
	\begin{align*}
		\varphi_{\bX}(\bu) & = \exp\Bigg\{-\int_{C_d^{\|\cdot\|}}\Big(|\scal{\bu,\dfrac{\bs}{\|\bs\|_e}}|+ia\scal{\bu,\dfrac{\bs}{\|\bs\|_e}}\ln|\scal{\bu,\dfrac{\bs}{\|\bs\|_e}}|\Big)\|\bs\|_e\gamma^{\|\cdot\|}(d\bs)\\
		& \hspace{5cm} +i\,\scal{\bu,\bm_{{\|\cdot\|}}^0}-ia\int_{C_d^{\|\cdot\|}}\scal{\bu,\bs}\ln\|\bs\|_e\gamma^{\|\cdot\|}(d\bs)\Bigg\},\\
		& = \exp\Bigg\{-\int_{S_d\setminus K^{\|\cdot\|}}\Big(|\scal{\bu,\bs'}|+ia\scal{\bu,\bs'}\ln|\scal{\bu,\bs'}|\Big)\gamma(d\bs')\\
		& \hspace{5cm}
		+i\,\scal{\bu,\bm_{{\|\cdot\|}}^0}-ia\int_{S_d\setminus K^{\|\cdot\|}}\scal{\bu,\bs'}\ln\|\bs'\|\gamma(d\bs')\Bigg\},
	\end{align*}
	where we used the change of variable $\bs' = T_{\|\cdot\|}^{-1}(\bs)=\bs/\|\bs\|_e$, and  $\gamma(d\bs):=\|\bs\|^{-1}\gamma^{\|\cdot\|}\circ T_{\|\cdot\|}(d\bs)$. 
	Letting then $\overline{\gamma}(A):=\gamma(A\cap (S_d\setminus K^{\|\cdot\|}))$ for any Borel set $A$ of $S_d$ and $\tilde{\bm} :=(\tilde{m}_i)$ with $\tilde{m}_i = \int_{S_d\setminus K^{\|\cdot\|}}s_i\ln\|\bs\| \overline{\gamma}(d\bs)$, $j=1,\ldots,d$, we get
	\begin{align*}
		\varphi_{\bX}(\bu) & = \exp\Bigg\{-\int_{S_d}\Big(|\scal{\bu,\bs}|+ia\scal{\bu,\bs}\ln|\scal{\bu,\bs}|\Big)\overline{\gamma}(d\bs)+i\,\scal{\bu,\bm_{{\|\cdot\|}}^0-a\tilde{\bm}}\Bigg\},
	\end{align*}
	and $\bX$ admits the pair $(\overline{\gamma},\bm_{{\|\cdot\|}}^0-a\tilde{\bm})$ for spectral representation on the Euclidean unit sphere.
	The unicity of the spectral representation of $\bX$ on $S_d$ implies that $(\Gamma,\bmu^0)=(\overline{\gamma},\bm_{{\|\cdot\|}}^0-a\tilde{\bm})$. Thus, $\overline{\gamma}$ and $\Gamma$ have to coincide, and in particular
	\begin{align*}
		\Gamma(K^{\|\cdot\|}) & = \overline{\gamma}(K^{\|\cdot\|}) = \gamma(K^{\|\cdot\|}\cap (S_d\setminus K^{\|\cdot\|})) = \gamma(\emptyset)=0,\\
	    \tilde{m}_i & = \int_{S_d\setminus K^{\|\cdot\|}}s_i\ln\|\bs\|\Gamma(d\bs), \hspace{0.3cm} i=1,\ldots,d.
	\end{align*}
    Last, as $\int_{C_d^{\|\cdot\|}}\|\bs\|_e\Big|\ln\|\bs\|_e\hspace{0.03cm}\Big|\gamma^{\|\cdot\|}(d\bs)<+\infty$ (Lemma \ref{le:a1_est_complique}) and $\Gamma(K^{\|\cdot\|})=0$, we have by a change of variable 
    \begin{align*}
    \int_{C_d^{\|\cdot\|}}\|\bs\|_e\Big|\ln\|\bs\|_e\hspace{0.03cm}\Big|\gamma^{\|\cdot\|}(d\bs) & = \int_{S_d\setminus K^{\|\cdot\|}}\Big|\ln\|\bs\|\hspace{0.03cm}\Big|\|\bs\|^{-1}\gamma^{\|\cdot\|}\circ T_{\|\cdot\|}(d\bs)\\
    & = \int_{S_d\setminus K^{\|\cdot\|}}\Big|\ln\|\bs\|\hspace{0.03cm}\Big|\gamma(d\bs)\\
    & =  \int_{S_d}\Big|\ln\|\bs\|\hspace{0.03cm}\Big|\Gamma(d\bs)\\
    & <+\infty,
    \end{align*}
    which concludes the proof of Proposition \ref{prop238}.\\
    
    \textbf{Proof of Lemma \ref{le:a1_est_complique}}\\  
    Notice that there exists a positive real number $b$ such that for all $\bs\in C_d^{\|\cdot\|}$, $\|\bs\|_e\ge b$ because $\|\bs\|=1$.
    Letting $M>0$, we have for all $\bu\in\mathbb{R}^d$
    \begin{align*}
    \int_{C_d^{\|\cdot\|}}\|\bs\|_e\Big|\ln\|\bs\|_e\hspace{0.03cm}\Big|\gamma^{\|\cdot\|}(d\bs) & = \int_{C_d^{\|\cdot\|}\cap  \{b\le\|\bs\|_e\le M\}} + \int_{C_d^{\|\cdot\|}\cap  \{\|\bs\|_e > M\}} := I_1 + I_2.
    \end{align*}
    We will show that both $I_1$ and $I_2$ are finite. Focus first on $I_2$.
    From \eqref{asu:integ_complique}, we know that for all $\bu\in\mathbb{R}^d$
    \begin{align}\label{dem:24}
    \int_{C_d^{\|\cdot\|}}|\scal{\bu,\bs}|\Big|\ln|\scal{\bu,\bs}|\hspace{0.03cm}\Big|\gamma^{\|\cdot\|}(d\bs) & = \int_{C_d^{\|\cdot\|}\cap  \{b\le\|\bs\|_e\le M\}} + \int_{C_d^{\|\cdot\|}\cap  \{\|\bs\|_e > M\}} <+\infty.
    \end{align}
    and thus, in particular
    \begin{align}
    & \int_{\{\bs'\in C_d^{\|\cdot\|}:\hspace{0.1cm} \|\bs'\|_e > M\}}|\scal{\bu,\bs}|\Big|\ln|\scal{\bu,\bs}|\hspace{0.03cm}\Big|\gamma^{\|\cdot\|}(d\bs)\nonumber\\ 
    & \hspace{4cm} = \int_{\{\bs'\in C_d^{\|\cdot\|}:\hspace{0.1cm} \|\bs'\|_e > M\}} |\scal{\bu,\bs}|\Big|\ln\|\bs\|_e +\ln|\scal{\bu,\dfrac{\bs}{\|\bs\|_e}}| \hspace{0.03cm}\Big|\gamma^{\|\cdot\|}(d\bs)<+\infty.\label{dem:prebeforepart}
    \end{align}
    By the triangular inequality, for all $\bu\in\mathbb{R}^d$,
    \begin{align}
    & \int_{\{\bs'\in C_d^{\|\cdot\|}:\hspace{0.1cm} \|\bs'\|_e > M\}} |\scal{\bu,\bs}|\Big|\ln\|\bs\|_e +\ln|\scal{\bu,\dfrac{\bs}{\|\bs\|_e}}| \hspace{0.03cm}\Big|\gamma^{\|\cdot\|}(d\bs) \nonumber\\
    & \hspace{4cm} = \int_{\{\bs'\in C_d^{\|\cdot\|}:\hspace{0.1cm} \|\bs'\|_e > M\}} |\scal{\bu,\bs}|\Big|\ln\|\bs\|_e\Big|\bigg|1 + \dfrac{\ln|\scal{\bu,\bs/\|\bs\|_e}|}{\ln\|\bs\|_e} \hspace{0.03cm}\bigg|\gamma^{\|\cdot\|}(d\bs)\nonumber\\
    & \hspace{4cm} \ge \int_{\{\bs'\in C_d^{\|\cdot\|}:\hspace{0.1cm} \|\bs'\|_e > M\}} |\scal{\bu,\bs}|\Big|\ln\|\bs\|_e\Big|\Bigg|1 - \bigg|\dfrac{\ln|\scal{\bu,\bs/\|\bs\|_e}|}{\ln\|\bs\|_e} \hspace{0.03cm}\bigg|\Bigg|\gamma^{\|\cdot\|}(d\bs)\label{dem:beforepart}
    \end{align}
    Let us now partition the space $\mathbb{R}^d$ into subsets $R_1,\ldots,R_d$ such that, for any $i=1,\ldots,d$ and any $\bs=(s_1,\ldots,s_d)\in R_i$, $\sup\limits_j |s_j| = |s_i|$.\footnote{Strictly speaking, $(R_1,\ldots,R_d)$ is not a partition of $\mathbb{R}^d$ as the $R_i$'s may intersect because of ties in the components of vectors.
    This will not affect the proof.}
    We have by \eqref{dem:prebeforepart}-\eqref{dem:beforepart} that for any $i=1,\ldots,d$, any $\bu\in\mathbb{R}^d$,
    \begin{align*}
    \int_{\{\bs'\in C_d^{\|\cdot\|}:\hspace{0.1cm} \|\bs'\|_e > M\}\cap R_i} |\scal{\bu,\bs}|\Big|\ln\|\bs\|_e\Big|\Bigg|1 - \bigg|\dfrac{\ln|\scal{\bu,\bs/\|\bs\|_e}|}{\ln\|\bs\|_e} \hspace{0.03cm}\bigg|\Bigg|\gamma^{\|\cdot\|}(d\bs) < +\infty.
    \end{align*}
    Denoting $(\be_1,\ldots,\be_d)$ the canonical orthonormal basis of $\mathbb{R}^d$, evaluate now the above at $\bu=\be_i$. We get that
    \begin{align}
    \int_{\{\bs'\in C_d^{\|\cdot\|}:\hspace{0.1cm} \|\bs'\|_e > M\}\cap R_i} |\scal{\be_i,\bs}|\Big|\ln\|\bs\|_e\Big|\Bigg|1 - \bigg|\dfrac{\ln|\scal{\be_i,\bs/\|\bs\|_e}|}{\ln\|\bs\|_e} \hspace{0.03cm}\bigg|\Bigg|\gamma^{\|\cdot\|}(d\bs) < +\infty.\label{dem:afterpart_beforeboundcos}
    \end{align}
    Let us show that $\bs \longmapsto \ln|\scal{\be_i,\bs/\|\bs\|_e}|$ is a bounded function for $\bs\in \{\bs'\in C_d^{\|\cdot\|}:\hspace{0.1cm} \|\bs'\|_e > M\}\cap R_i$.
    \textit{Ad absurdum}, if it is not bounded, then for any $A>0$, there exists $\bs\in \{\bs'\in C_d^{\|\cdot\|}:\hspace{0.1cm} \|\bs'\|_e > M\}\cap R_i$ such that
    $$
    \Big|\ln|\scal{\be_i,\bs/\|\bs\|_e}|\Big|>A.
    $$
    Taking the sequence $A_n = n$ for any $n\ge1$, we get that there exists a sequence $(\bs_n)$, $\bs_n\in \{\bs'\in C_d^{\|\cdot\|}:\hspace{0.1cm} \|\bs'\|_e > M\}\cap R_i$ such that
    \begin{align*}
    \Big|\ln|\scal{\be_i,\bs_n/\|\bs_n\|_e}|\Big|>n.
    \end{align*}
    Thus, for all $n\ge1$
    \begin{align*}
    0 \le |\scal{\be_i,\bs_n/\|\bs_n\|_e}| \le e^{-n}.
    \end{align*}
    and
    $$
    |\scal{\be_i,\bs_n/\|\bs_n\|_e}|  \underset{n\rightarrow+\infty}{\longrightarrow} 0.
    $$
    Consider now the decomposition of $\bs_n/\|\bs_n\|_e$ in the orthonormal basis $(\be_1,\ldots,\be_d)$,
    $$
    \bs_n/\|\bs_n\|_e = \sum_{j=1}^d \scal{\be_j,\bs_n/\|\bs_n\|_e}\be_j.
    $$
    As $\bs_n\in R_i$ for all $n\ge 1$, we also have that $\bs_n/\|\bs_n\|_e \in R_i$ for all $n\ge1$, and thus, for any $j=1,\ldots,d$
    \begin{align*}
    0 \le |\scal{\be_j,\bs_n/\|\bs_n\|_e}| \le |\scal{\be_i,\bs_n/\|\bs_n\|_e}| \underset{n\rightarrow+\infty}{\longrightarrow} 0.
    \end{align*}
    Hence, $\bs_n/\|\bs_n\|_e\underset{n\rightarrow+\infty}{\longrightarrow} 0$, which is impossible since $\Big\|\bs_n/\|\bs_n\|_e \Big\|_e = 1$ for all $n\ge1$.
    The function $\bs \longmapsto \ln|\scal{\be_i,\bs/\|\bs\|_e}|$ is thus bounded on $\{\bs\in C_d^{\|\cdot\|}:\hspace{0.1cm} \|\bs\|_e > M\}\cap R_i$, say $\Big|\ln|\scal{\be_i,\bs/\|\bs\|_e}|\Big| \le A$ for some $A>0$.
    Provided $M$ is taken large enough (e.g., $M>2A$), we will have in \eqref{dem:afterpart_beforeboundcos}
    $$
   \Bigg|1 - \bigg|\dfrac{\ln|\scal{\be_i,\bs/\|\bs\|_e}|}{\ln\|\bs\|_e} \hspace{0.03cm}\bigg|\Bigg| = 1 - \bigg|\dfrac{\ln|\scal{\be_i,\bs/\|\bs\|_e}|}{\ln\|\bs\|_e}\bigg| \ge 1 - \dfrac{A}{M} > 0,
    $$
    which thus yields for all $i=1,\ldots,d$
    \begin{align*}
    \int_{\{\bs'\in C_d^{\|\cdot\|}:\hspace{0.1cm} \|\bs'\|_e > M\}\cap R_i} |\scal{\be_i,\bs}|\Big|\ln\|\bs\|_e\Big|\gamma^{\|\cdot\|}(d\bs) < +\infty.
    \end{align*}
    As $|\scal{\be_i,\bs}| \ge \|\bs\|_e e^{-A}$, we further get that
    \begin{align*}
    \int_{\{\bs'\in C_d^{\|\cdot\|}:\hspace{0.1cm} \|\bs'\|_e > M\}\cap R_i} \|\bs\|_e\Big|\ln\|\bs\|_e\Big|\gamma^{\|\cdot\|}(d\bs) < +\infty,
    \end{align*}
    and because $\bigcup\limits_{i=1,\ldots,d} R_i = \mathbb{R}^d$, 
    \begin{align*}
    & I_2 = \int_{\{\bs'\in C_d^{\|\cdot\|}:\hspace{0.1cm} \|\bs'\|_e > M\}} \|\bs\|_e\Big|\ln\|\bs\|_e\Big|\gamma^{\|\cdot\|}(d\bs)\\
    & \hspace{3cm} \le \sum_{i=1}^d \int_{\{\bs'\in C_d^{\|\cdot\|}:\hspace{0.1cm} \|\bs'\|_e > M\}\cap R_i} \|\bs\|_e\Big|\ln\|\bs\|_e\Big|\gamma^{\|\cdot\|}(d\bs) < +\infty.
    \end{align*}

    \vspace{1cm}
    
    Let us now show that $I_1$ is finite. Assuming for a moment that 
    $$
    \gamma^{\|\cdot\|}\Big(\{\bs'\in C_d^{\|\cdot\|}:\hspace{0.1cm} b\le\|\bs'\|_e \le M\}\Big)<+\infty,
    $$
    we get
    \begin{align*}
    I_1 & = \int_{\{\bs'\in C_d^{\|\cdot\|}:\hspace{0.1cm} b\le\|\bs'\|_e  \le M\}}\|\bs\|_e\Big|\ln\|\bs\|_e\hspace{0.03cm}\Big|\gamma^{\|\cdot\|}(d\bs) \\
    & \hspace{3cm} 
    \le \Big({\max_{x\in[b,M]} x|\ln x|}\Big) \hspace{0.2cm} \gamma^{\|\cdot\|}\Big(\{\bs'\in C_d^{\|\cdot\|}:\hspace{0.1cm} b\le\|\bs'\|_e  \le M\}\Big),
    \end{align*}
    because $x\longmapsto x|\ln x|$ is a bounded function on $[b,M]$, and thus $I_1<+\infty$.
    We now show that $\gamma^{\|\cdot\|}$ is indeed finite on the set $\{\bs'\in C_d^{\|\cdot\|}:\hspace{0.1cm} b\le\|\bs'\|_e \le M\}$.\\
    Proceeding as in the case of $I_2$, it can be obtained that for $i=1,\ldots,d$, the function $\bs \longmapsto \ln|\scal{\be_i,\bs/\|\bs\|_e}|$ is bounded on the set $\{\bs'\in C_d^{\|\cdot\|}:\hspace{0.1cm} b\le\|\bs'\|_e \le M\}\cap R_i$.
    Say, again, that $\Big|\ln|\scal{\be_i,\bs/\|\bs\|_e}|\Big| \le A$ for some $A>0$.
    Then, $|\scal{\be_i,\bs}|\ge \|\bs\|_e e^{-A}$, and for any $\lambda>2b^{-1}e^{A}$, we have
    $$
    |\scal{\lambda\be_i,\bs}|\ge 2,
    $$
    for any $i=1,\ldots,d$, $\bs\in\{\bs'\in C_d^{\|\cdot\|}:\hspace{0.1cm} b\le\|\bs'\|_e \le M\}\cap R_i$.
    From \eqref{dem:24}, we have for any $\bu\in\mathbb{R}^d$
    $$
    \int_{\{\bs'\in C_d^{\|\cdot\|}:\hspace{0.1cm} b\le\|\bs'\|_e \le M\}} |\scal{\bu,\bs}|\Big|\ln|\scal{\bu,\bs}|\hspace{0.03cm}\Big|\gamma^{\|\cdot\|}(d\bs) <+\infty,
    $$
    and thus, for any $\bu\in\mathbb{R}^d$,
    $$
    \int_{\{\bs'\in C_d^{\|\cdot\|}:\hspace{0.1cm} b\le\|\bs'\|_e \le M\}\cap R_i} |\scal{\bu,\bs}|\Big|\ln|\scal{\bu,\bs}|\hspace{0.03cm}\Big|\gamma^{\|\cdot\|}(d\bs) <+\infty,
    $$
    for any $i=1,\ldots,d$.
    Evaluating the above in particular at $\bu=\lambda\be_i$, for any $\lambda>2b^{-1}e^{A}$, we get
    $$
    \int_{\{\bs'\in C_d^{\|\cdot\|}:\hspace{0.1cm} b\le\|\bs'\|_e \le M\}\cap R_i} |\scal{\lambda\be_i,\bs}|\Big|\ln|\scal{\lambda\be_i,\bs}|\hspace{0.03cm}\Big|\gamma^{\|\cdot\|}(d\bs) <+\infty.
    $$
    Noticing that $x\longmapsto x |\ln x|$ is increasing on $[1,+\infty)$ and that $|\scal{\lambda\be_i,\bs}|\ge 2$ for any $\bs$ in the domain of integration, we have $|\scal{\bu,\bs}|\Big|\ln|\scal{\bu,\bs}|\Big| \ge 2 \ln 2$, and
    \begin{align*}
    \int_{\{\bs'\in C_d^{\|\cdot\|}:\hspace{0.1cm} b\le\|\bs'\|_e \le M\}\cap R_i} \gamma^{\|\cdot\|}(d\bs) <+\infty,
    \end{align*}
    for any $i=1,\ldots,d$. Hence,
    $$
    \int_{\{\bs'\in C_d^{\|\cdot\|}:\hspace{0.1cm} b\le\|\bs'\|_e \le M\}} \gamma^{\|\cdot\|}(d\bs) \le \sum_{i=1}^d \int_{\{\bs'\in C_d^{\|\cdot\|}:\hspace{0.1cm} b\le\|\bs'\|_e \le M\}\cap R_i} \gamma^{\|\cdot\|}(d\bs)  <+\infty,
    $$
    and $\gamma^{\|\cdot\|}\Big(\{\bs'\in C_d^{\|\cdot\|}:\hspace{0.1cm} b\le\|\bs'\|_e \le M\}\Big)$ is finite.
    \cqfd
	
	\subsection{Proof of Proposition \ref{prop:cond_tail}}
	
	The proposition is an immediate consequence of Bayes formula and of the following result, which is an adaptation of Theorem 4.4.8 by Samorodnitsky and Taqqu (1994) \cite{st94} to seminorms.
	\begin{prop}\label{theo448}
		Let $\bX=(X_1,\ldots,X_d)$ be an $\alpha$-stable random vector and let $\|\cdot\|$ be a seminorm on $\mathbb{R}^d$ such that $\bX$ is representable on $C_d^{\|\cdot\|}$. Then, for every Borel set $A\subseteq C_d^{\|\cdot\|}$ with $\Gamma^{\|\cdot\|}(\partial A)=0$,
		\begin{align}
			\lim_{x\rightarrow+\infty} x^{\alpha} \mathbb{P}\Big(\|\bX\|>x,\dfrac{\bX}{\|\bX\|}\in A\Big)=C_\alpha \Gamma^{\|\cdot\|}(A),
		\end{align}
		with $C_\alpha = \dfrac{1-\alpha}{\Gamma(2-\alpha)\cos(\pi\alpha/2)}$ if $\alpha\ne1$, and $C_1 =2/\pi$.
	\end{prop}
	\noindent \textit{Proof.}\\
		We follow the proof of Theorem 4.4.8 by Samorodnitsky and Taqqu (1994) \cite{st94}. The main hurdle is to show that, with $\|\cdot\|$ a semi-norm, $K^{\|\cdot\|}=\{\bs\in S_d:\,\|\bs\|=0\}$, and $\Gamma^{\|\cdot\|}(K^{\|\cdot\|})=0$, we have the series representation of $\bX$, $(X_1,\ldots,X_d)\stackrel{d}{=}(Z_1,\ldots,Z_d)$ where 
	\begin{align}\label{eq:seriesrepres}
		Z_k = (C_\alpha\Gamma^{\|\cdot\|}(C_d^{\|\cdot\|}))^{1/\alpha}\sum_{i=1}^{\infty} [\Gamma_i^{-1/\alpha}S_i^{(k)}-b_{i,k}(\alpha)],\hspace{0.2cm}k=1,\ldots,d,
	\end{align}
	with $\bS_i = (S_i^{(1)},\ldots,S_i^{(d)})$, $i\ge1$, are i.i.d. $C_d^{\|\cdot\|}$-valued random vectors with common law $\Gamma^{\|\cdot\|}/\Gamma^{\|\cdot\|}(C_d^{\|\cdot\|})$ and the $b_{i,k(\alpha)}$'s are constants.\\
	By Proposition \ref{prop238}, we know that $\bX$ admits a characteristic function of the form \eqref{def:char_fun_vec}. This allows to restate the integral representation Theorem 3.5.6 in \cite{st94} on the semi-norm unit cylinder as follows: with the measurable space $(E, \mathcal{E})=(C_d^{\|\cdot\|}, \text{Borel}\,\sigma\!\text{-algebra on}\, C_d^{\|\cdot\|})$, let $M$ be an $\alpha$-stable random measure on $(E,\mathcal{E})$ with control measure $m=\Gamma^{\|\cdot\|}$, skewness intensity $\beta(\,\cdot\,)\equiv1$ (see Definition 3.3.1 in \cite{st94} for details). Letting also $f_j:C_d^{\|\cdot\|}\longrightarrow\mathbb{R}$ defined by $f_j\Big((s_1,\ldots,s_d)\Big)=s_j$, $j=1,\ldots,d$, then
	\begin{align*}
		\bX \stackrel{d}{=} \bigg(\int_{C_d^{\|\cdot\|}}f_1(\bs)M(d\bs),\ldots,\int_{C_d^{\|\cdot\|}}f_d(\bs)M(d\bs)\bigg) + \bmu^{\|\cdot\|}.
	\end{align*}
	This representation can be checked directly by comparing the characteristic functions of the left-hand and right-hand sides. We can now apply Theorem 3.10.1 in \cite{st94} to the above integral representation with $(E, \mathcal{E},m)$ the measure space as described before, and $\hat{m}=\Gamma^{\|\cdot\|}/\Gamma^{\|\cdot\|}(C_d^{\|\cdot\|})$. This establishes \eqref{eq:seriesrepres}. The rest of the proof is similar to that of Theorem 4.4.8 in \cite{st94}. We rely on the triangle inequality property of semi-norms and the fact that any norm  is finer than any semi-norm in finite dimension.\footnote{We say that a norm $N$ is finer than a semi-norm $N_s$ if there is a positive constant $C$ such that $N_s(x) \le C N(x)$ for any $x\in\mathbb{R}^d$.}
	\cqfd
	
	\subsection{Proof of Lemma \ref{le:ma_representable}}

	From Proposition \ref{prop238}, we know that a necessary condition for the representability of $\bXt$ on $C_{m+h+1}^{\|\cdot\|}$  is
	$\Gamma(K^{\|\cdot\|})=0$, where $K^{\|\cdot\|}=\{\bs\in S_{m+h+1}:\|\bs\|=0\}$. 
	This condition is also sufficient when either $\alpha\ne1$ or $\alpha=1$, $\beta=0$.
	Using the fact that $\Gamma$ only charges discrete atoms on $C_{m+h+1}^{\|\cdot\|}$,
	\begin{align*}
		\Gamma(K^{\|\cdot\|})=0
		& \iff \{\bs\in S_{m+h+1}:\Gamma(\{\bs\})>0\} \cap K^{\|\cdot\|} = \emptyset\\
		& \iff \forall \bs\in {S}_{m+h+1}, \quad \Big[\Gamma(\{\bs\})>0 \Longrightarrow \|\bs\| > 0\Big]\\
		& \iff \forall k \in \mathbb{Z}, \quad \Big[\|\bd_k\|_e>0 \Longrightarrow \|\bd_k\|>0\Big]\\
		& \iff \forall k \in \mathbb{Z}, \quad \Big[\|\bd_k\|=0  \Longrightarrow \|\bd_k\|_e=0 \Big]\\
		& \iff \forall k \in \mathbb{Z}, \quad \Big[\|\bd_k\|=0  \Longrightarrow \bd_k=0 \Big]\\
		& \iff \forall k \in \mathbb{Z}, \quad \Big[(d_{k+m},\ldots,d_{k})=\boldsymbol{0} \Longrightarrow (d_{k+m},\ldots,d_{k-h})=\boldsymbol{0} \Big],
	\end{align*}
	by \eqref{asu:semin}. Now assume that the following holds:
	\begin{align}\label{dem:rec}
		\forall k \in \mathbb{Z}, \quad \Big[(d_{k+m},\ldots,d_{k})=\boldsymbol{0} \Longrightarrow (d_{k+m},\ldots,d_{k-h})=\boldsymbol{0} \Big].
	\end{align}
	Then, if for some particular $k_0\in\mathbb{Z}$, we have
	\begin{align*}
		(d_{k_0+m},\ldots,d_{k_0}) = \boldsymbol{0}. 
	\end{align*}
	It implies that
	\begin{align*}
		(d_{k_0+m},\ldots,d_{k_0-h}) = \boldsymbol{0},
	\end{align*}
	and especially, as we assume $h\ge1$,
	\begin{align*}
		(d_{(k_0-1)+m},\ldots,d_{k_0-1}) = \boldsymbol{0}.
	\end{align*}
	Invoking \eqref{dem:rec}, we deduce by recurrence that for any $n\ge0$, 
	$$
	(d_{(k_0-n)+m},\ldots,d_{k_0-n}) = \boldsymbol{0}.
	$$
	Therefore, \eqref{dem:rec} implies
	\begin{align*}
		\forall k \in \mathbb{Z}, \quad \Big[(d_{k+m},\ldots,d_{k})=\boldsymbol{0} \Longrightarrow \forall \ell \le k-1, \quad d_\ell = 0 \Big]
	\end{align*}
	The reciprocal is clearly true.
	This establishes that \eqref{eq:ma_representable} is a necessary and sufficient condition for $\bX_t$ to be representable on $C_d^{\|\cdot\|}$ in the cases where either $\alpha\ne1$, or $\alpha=1$, $\beta=0$.\\
	
	In the case $\alpha=1$, $\beta\ne0$, Proposition \ref{prop238} states that the necessary and sufficient condition for representability reads $\int_{S_d}\Big|\ln\|\bs\|\hspace{0.03cm}\Big|\Gamma(d\bs) < +\infty$. That is
	\begin{align*}
    \Gamma(K^{\|\cdot\|})=0 \hspace{1cm} \text{and} \hspace{1cm} \int_{S_d\setminus K^{\|\cdot\|}}\Big|\ln\|\bs\|\hspace{0.03cm}\Big|\Gamma(d\bs) < +\infty.
	\end{align*}
	Substituting $\Gamma$ by its expression in \eqref{def:spectral}, the above condition holds if and only if \eqref{eq:ma_representable} is true and
	\begin{align*}
	\sigma\sum_{\vartheta\in S_1} \sum_{k\in\mathbb{Z}} w_\vartheta \|\bd_k\|_e \Bigg|\ln\bigg\|\dfrac{\vartheta\bd_k}{\|\bd_k\|_e}\bigg\|\hspace{0.05cm}\Bigg| < +\infty,
	\end{align*}
	the latter being equivalent to
	\begin{align*}
	\sum_{k\in\mathbb{Z}} \|\bd_k\|_e \Bigg|\ln\dfrac{\|\bd_k\|}{\hspace{0.2cm}\|\bd_k\|_e}\hspace{0.05cm}\Bigg| < +\infty.
	\end{align*}
	
	\subsection{Proof of Proposition \ref{prop:ma_pastrepres}}
	
	By Definition \ref{def:ma_representable}, $(X_t)$ is past-representable if and only if there exists $m\ge0$, $h\ge1$ such that the vector $(X_{t-m},\ldots,X_t,X_{t+1},\ldots,X_{t+h})$ is representable on $C_{m+h+1}^{\|\cdot\|}$.
	Consider first point $(\iota)(a)$, that is, the case $\alpha\ne1$, $(\alpha,\beta)=(1,0)$.
    By Lemma \ref{le:ma_representable}, 
	\begin{align*}
		(X_t) \hspace{0.2cm} \text{is past-representable}  & \hspace{0.1cm}\iff\hspace{0.1cm} \text{There exist} \hspace{0.1cm} m\ge0, \hspace{0.1cm} h\ge1, \hspace{0.1cm} \text{such that \eqref{eq:ma_representable} holds} \\
		& \hspace{0.1cm}\iff\hspace{0.1cm}  \exists m\ge0, \forall k \in\mathbb{Z}, \Big[d_{k+m}=\ldots=d_k=0  \hspace{0.2cm}\Longrightarrow \hspace{0.2cm} \forall  \ell \le k-1, \hspace{0.2cm} d_\ell=0\Big].
	\end{align*}
	Thus,
	\begin{align*}
		(X_t) \hspace{0.2cm} \text{not past-representable}  & \hspace{0.1cm}\iff\hspace{0.1cm} \forall m\ge0, \exists k \in\mathbb{Z}, d_{k+m}=\ldots=d_k=0  \hspace{0.2cm}\text{and} \hspace{0.2cm} \exists  \ell \le k-1, \hspace{0.2cm} d_\ell\ne0\\
		& \hspace{0.1cm}\iff\hspace{0.1cm} \forall m\ge0, \exists k \in\mathbb{Z}, d_{k+m}=\ldots=d_k=0  \hspace{0.2cm}\text{and} \hspace{0.2cm} d_{k-1}\ne0\\
		& \hspace{0.1cm}\iff\hspace{0.1cm} \forall m\ge1, \exists k \in\mathbb{Z}, d_{k+m}=\ldots=d_{k+1}=0  \hspace{0.2cm}\text{and} \hspace{0.2cm} d_{k}\ne0\\
		& \hspace{0.1cm}\iff\hspace{0.1cm} \sup\{m\ge1: \hspace{0.2cm} \exists k \in\mathbb{Z}, \hspace{0.2cm} d_{k+m}=\ldots=d_{k+1}=0 , \hspace{0.2cm} d_{k}\ne0\} = +\infty,
	\end{align*}
	hence \eqref{eq:ma_pastrepres}.\\
	Regarding the last statement of point $(\iota)(a)$, assume first that $m_0<+\infty$ and $m\ge m_0$. 
   Property \eqref{eq:ma_representable} necessarily holds with $m_0$.
    Indeed, if it did not, there would exist $k\in\mathbb{Z}$ such that
	$$
	d_{k+m_0}=\ldots=d_k= 0, \hspace{0.3cm} \text{and} \hspace{0.3cm} d_\ell \ne0, \hspace{0.3cm} \text{for some} \hspace{0.3cm} \ell \le k-1,
	$$
	and we would have found a sequence of consecutive zero values of length at least $m_0+1$ preceded by a non-zero value, contradicting the fact that
	$$
	m_0=\sup\{m\ge1:\hspace{0.1cm}\exists\hspace{0.1cm}k\in\mathbb{Z}, \hspace{0.1cm} d_{k+m}=\ldots=d_{k+1}= 0, \hspace{0.3cm} \text{and}  \hspace{0.3cm} d_k \ne0\}.
	$$
	As \eqref{eq:ma_representable} holds with $m_0$, it holds \textit{a fortiori} for any $m'\ge m_0$.
	Thus, $\bX_t=(X_{t-m},\ldots,X_t,X_{t+1},\ldots,X_{t+h})$ is representable for any $m'\ge m_0$, $h\ge1$ by Lemma \ref{le:ma_representable}, and $(X_t)$ is in particular $(m,h)$-past-representable.\\
	Reciprocally let $m\ge0$, $h\ge1$ and assume that $(X_t)$ is $(m,h)$-past-representable. 
	The process $(X_t)$ is thus in particular past-representable, which as we have shown previously, implies that $m_0 < +\infty$.
	\textit{Ad absurdum}, suppose now that $0\le m < m_0 < +\infty$.
	If $m_0=0$, there is nothing to do. 
	Otherwise if $m_0\ge1$, by definition, there exists a $k\in\mathbb{Z}$ such that
	\begin{equation}\label{eq:m0}
		d_{k+m_0}=\ldots=d_{k+1} = 0, \quad \text{and} \quad  d_k \ne 0.
	\end{equation}
	Because $(X_t)$ is $(m,h)$-past-representable, we have by Lemma \ref{le:ma_representable} that
	\eqref{eq:ma_representable} holds with $m$. As $m<m_0$ and  $d_{k+m_0}=\ldots=d_{k+1} = 0$, we thus have that  $d_\ell=0$ for all $\ell\le k+1$, and in particular $d_k=0$, hence the contradiction. We conclude that $m\ge m_0$.\\ 
	
	\noindent Consider now point $(\iota)(b)$, i.e., the case $\alpha=1$ and $\beta\ne0$.
	From Lemma \ref{le:ma_representable},
	\begin{align*}
		(X_t) \hspace{0.2cm} \text{is past-representable}  & \hspace{0.1cm}\iff\hspace{0.1cm} \text{There exist} \hspace{0.1cm} m\ge0, \hspace{0.1cm} h\ge1, \hspace{0.1cm} \text{such that \eqref{eq:ma_representable} and \eqref{eq:additionallourd} hold}
	\end{align*}
	From the previous proof, we moreover have that
	$$
	 \exists\hspace{0.1cm} m\ge0, \hspace{0.1cm} 
	 \text{such that \eqref{eq:ma_representable} holds} \hspace{0.1cm} \iff \hspace{0.1cm} m_0<+\infty \hspace{0.1cm} \iff 
	 \hspace{0.1cm} \left\{ \begin{array}{l}
	     m_0 < +\infty \\
	     \forall \hspace{0.1cm} m'\ge m_0, \hspace{0.1cm} 
	     \text{\eqref{eq:ma_representable} holds}\\
	     \forall \hspace{0.1cm} m'< m_0, \hspace{0.1cm} \text{\eqref{eq:ma_representable} does not hold}
		\end{array}\right.
	$$
	Hence 
	\begin{align*}
		& \exists \hspace{0.1cm} m\ge0, \hspace{0.1cm} h\ge1, \hspace{0.1cm} \text{such that \eqref{eq:ma_representable} and \eqref{eq:additionallourd} hold}\\
		& \hspace{4cm}\iff\hspace{0.1cm}
		\left\{ \begin{array}{l}
	     m_0 < +\infty \\
	     \forall \hspace{0.1cm} m'\ge m_0, \hspace{0.1cm} 
	     \text{\eqref{eq:ma_representable} holds}\\
	     \forall \hspace{0.1cm} m'< m_0, \hspace{0.1cm} \text{\eqref{eq:ma_representable} does not hold}\\
	     \exists \hspace{0.1cm} m\ge0, \hspace{0.1cm} h\ge1, \hspace{0.1cm} \text{such that \eqref{eq:ma_representable} and \eqref{eq:additionallourd} hold.}
		\end{array}\right. 
	\end{align*}
	The latter in particular implies $m_0<+\infty$ and the existence of $m\ge m_0$, $h\ge1$ such that \eqref{eq:additionallourd} holds. 
	Reciprocally, 
		\begin{align*}
		& \left\{ \begin{array}{l}
	     m_0 < +\infty \\
	     \exists \hspace{0.1cm} m\ge m_0, \hspace{0.1cm} h\ge1, \hspace{0.1cm} \text{such that \eqref{eq:additionallourd} holds}
		\end{array}\right. \\
		& \hspace{5cm} \Longrightarrow 		
		\left\{ \begin{array}{l}
	     m_0 < +\infty \\
	     \forall \hspace{0.1cm} m'\ge m_0, \hspace{0.1cm} 
	     \text{\eqref{eq:ma_representable} holds}\\
	     \exists \hspace{0.1cm} m\ge m_0, \hspace{0.1cm} h\ge1, \hspace{0.1cm} \text{such that \eqref{eq:additionallourd} holds,}
		\end{array}\right.
	\end{align*}
	which in particular implies that there exists $m\ge m_0$, $h\ge1$ such that both \eqref{eq:ma_representable} and \eqref{eq:additionallourd} hold.
	Hence the past-representability of $(X_t)$.\\
	
	\noindent In view of Definition \ref{def:ma_representable}, point $(\iota\iota)$ is a direct consequence of the second part of Proposition \ref{prop238}.
	
	\subsection{Proof of Corollary \ref{cor:nonanticiAR}}
	
	Letting $k_0$ be the greatest integer such that $d_{k_0}\ne0$ (such an index exists by \eqref{assump:coef}), then immediately, for any $m\ge1$, $d_{k_0+m}=\ldots=d_{k_0+1}=0$ and therefore $m_0=+\infty$.
	
	\subsection{Proof of Corollary \ref{cor:mar}}
	
	We first show that $\deg(\psi)\ge1$ if and only if $m_0<+\infty$.\\
	
	Clearly, if $\deg(\psi)=0$, then $X_t = \sum_{k=-\infty}^{k_0}d_k\varepsilon_{t+k}$ for some $k_0$ in $\mathbb{Z}$ and $m_0=+\infty$. \\
	
	\noindent Reciprocally, assume $\deg(\psi)=p\ge1$. 
	Let us first show that \eqref{eq:ma_pastrepres} holds.\\
	Denote $\psi(F)\phi(B)=\sum_{i=-q}^p \varphi_i F^i$ and $\Theta(F)H(B)=\sum_{k=-r}^s \theta_i F^i$, for any non-negative degrees $q=\deg(\phi)$, $r=\deg(H)$, $s=\deg(\Theta)$. 
	From the recursive equation satisfied by $(X_t)$, we have that
	\begin{align}
	&&  \sum_{i=-q}^p \varphi_i X_{t+i} & = \sum_{k=-r}^s \theta_k \varepsilon_{t+k}\nonumber\\
	 \iff && \sum_{i=-q}^p \varphi_i \sum_{k\in\mathbb{Z}}d_k \varepsilon_{t+k+i} & = \sum_{k=-r}^s \theta_k \varepsilon_{t+k}\nonumber\\
	\iff && \sum_{k\in\mathbb{Z}} \bigg(\sum_{i=-q}^p \varphi_i d_{k-i}\bigg) \varepsilon_{t+k} & = \sum_{k=-r}^s \theta_k \varepsilon_{t+k}.\label{eq:recsurdk}
	\end{align}
	Proceeding by identification using the uniqueness of representation of heavy-tailed moving averages (see \cite{gz15}), we get that for $|k|>\max(r,s)$,
	\begin{align}\label{eq:rec}
		\sum_{i=-q}^p \varphi_i d_{k-i}=0.
	\end{align}
	\textit{Ad absurdum}, if $(X_t)$ is not past-representable, then by Proposition \ref{prop:ma_pastrepres} 
	$$
	\sup \{m\ge1: \hspace{0.2cm} \exists k\in\mathbb{Z}, \hspace{0.2cm} d_{k+m}=\ldots=d_{k+1}=0, \hspace{0.2cm} d_{k}\ne0\}=+\infty.
	$$
	Thus, there exists a sequence $\{m_n: \hspace{0.2cm n\ge 0}\}$, $m_n\ge1$, $\lim_{n\rightarrow+\infty}=+\infty$, satisfying: for any $n\ge0$, there is an index $k\in\mathbb{Z}$ such that 
	$$
	d_{k-p}\ne0 \hspace{0.3cm} \text{and} \hspace{0.3cm} d_{k-p+1}=d_{k-p+2}=\ldots=d_{k+m_n}=0.
	$$
	We can therefore construct a sequence $(k_n)$ such that the above relation holds for all $n\ge0$.
	This sequence of integers in $\mathbb{Z}$ is either bounded or unbounded.
	We will show that both cases lead to a contradiction.\\
	
	\noindent \textbf{First case: $\boldsymbol{\sup\{|k_n|: n\ge0\}=+\infty}$}
	
	There are two subsequences such that $m_{g(n)}\longrightarrow+\infty$ and $|k_{g(n)}|\longrightarrow+\infty$. For some $n$ large enough such that \eqref{eq:rec} holds and $m_{g(n)}\ge p+q$, we have both
	$$
	\sum_{i=-q}^p \varphi_i d_{k_{g(n)}-i}=0.
	$$
	and 
	$$
	d_{k_{g(n)}-p}\ne0, \quad d_{k_{g(n)}-p+1}=\ldots=d_{k_{g(n)}+q}=0.
	$$
	Hence, 
	$$
	\varphi_pd_{k_{g(n)}-p}=0,
	$$
	which is impossible given that $d_{k_{g(n)}-p}\ne0$ and $\varphi_{p}\ne0$.
	Indeed, denoting $\psi(z)=1+\psi_1z+\ldots+\psi_pz^p$, $\psi_p\ne0$ because $\deg(\psi)=p$, it can be shown that $\varphi_{p}=\psi_p$.\\ 
	
	\noindent \textbf{Second case: $\boldsymbol{\sup\{|k_n|: n\ge0\}<+\infty}$}
	
	Given that $(k_n)$ is a bounded sequence, there exists by the Bolzano-Weierstrass theorem a convergent subsquence $(k_{g(n)})$. 
	As $(k_{g(n)})$ takes only discrete values, it necessarily holds that $(k_{g(n)})$ reaches its limit at a finite integer $n_0\ge1$, that is, for all $n\ge n_0$, $k_{g(n)}=\lim_{n\rightarrow +\infty}k_{g(n)}:=\Bar{k}\in\mathbb{Z}$.
	Thus, for all $n\ge n_0$
	\begin{equation*}
	d_{\Bar{k}}\ne0, \hspace{0.3cm} \text{and} \hspace{0.3cm} d_{\Bar{k}+m_{g(n)}}=0, 
	\end{equation*}
	and as $m_{g(n)}\rightarrow+\infty$, we deduce that
	\begin{equation*}
	d_{\Bar{k}}\ne0, \hspace{0.3cm} \text{and} \hspace{0.3cm} d_{\Bar{k}+\ell}=0, \hspace{0.4cm} \text{for all} \hspace{0.3cm} \ell\ge1.
	\end{equation*}
	The process $(X_t)$ hence admit a moving average representation of the form
	\begin{equation}\label{dem:notinf}
	X_t = \sum_{k=-\infty}^{\Bar{k}} d_k \varepsilon_{t+k}, \hspace{1cm} t\in\mathbb{Z}.
	\end{equation}
	However, we also have by partial fraction decomposition
	\begin{align*}
	X_t & = \dfrac{\Theta(F)H(B)}{\psi(F)\phi(B)} \varepsilon_t\\
	& = \Theta(F)H(B) \dfrac{B^p}{B^p\psi(F)\phi(B)} \varepsilon_t\\
	& = \Theta(F)H(B)B^p \Bigg[ \dfrac{b_1(B)}{B^p\psi(F)} + \dfrac{b_2(B)}{\phi(B)} \Bigg]\varepsilon_t\\
	& = \Theta(F)H(B) \Bigg[ \dfrac{b_1(B)}{\psi(F)} + \dfrac{B^pb_2(B)}{\phi(B)} \Bigg]\varepsilon_t,
	\end{align*}
	for some polynomials $b1$ and $b_2$ such that $0\le\deg(b_1)\le p-1$, $0\le\deg(b_2)\le q-1$ and $\phi(B)b_1(B)+B^pb_2(B)\psi(F)=1$.
	We can write in general
	\begin{align*}
   \dfrac{\Theta(F)H(B)b_1(B)}{\psi(F)} & = \sum_{k=-\ell_1}^{+\infty}c_k \varepsilon_{t+k},\\
	\dfrac{\Theta(F)H(B)B^pb_2(B)}{\phi(B)} & =  \sum_{k=-\infty}^{\ell_2}e_k \varepsilon_{t+k},
	\end{align*}
	for some sequences of coefficients $(c_k)$, $(e_k)$, and where $\ell_1$ is the degree of the largest order monomial in $B$ of $\Theta(F)H(B)b_1(B)$ (recall that $F=B^{-1}$) and $\ell_2$ is the degree of the largest monomial in $F$ of $B^p\Theta(F)H(B) b_2(B)$.
	By \eqref{dem:notinf}, we deduce by identification that there is some $\Bar{\ell}\in \mathbb{Z}$ such that $c_k=0$ for all $k\ge\Bar{\ell}+1$ and
	$$
	\dfrac{\Theta(F)H(B)b_1(B)}{\psi(F)} = \sum_{k=-\ell_1}^{\Bar{\ell}}c_k F^k.
	$$
	Necessarily, $\Bar{\ell}\ge\ell_1$, otherwise $\Theta(F)H(B)b_1(B)\psi^{-1}(F)=0$ which is impossible as all the polynomials involved have non-negative degrees.
	Thus, we deduce that there exist two polynomials $P$ and $Q$ of non-negative degrees such that
	$$
	\dfrac{\Theta(z^{-1})H(z)b_1(z)}{\psi(z^{-1})} = \sum_{k=-\ell_1}^{\Bar{\ell}}c_kz^k := P(z^{-1})+Q(z), \hspace{1cm} z\in\mathbb{C}.,
	$$
	which yields
	\begin{align}\label{dem:pol_eq}
    \Theta(z^{-1})H(z)b_1(z) = \psi(z^{-1})(P(z^{-1})+Q(z)),\hspace{1cm} z\in\mathbb{C}.
	\end{align}
	As $\deg(\psi)=p$ and $\psi(z)=0$ if and only if $|z|>1$, we know that there are $p$ complex numbers $z_1,\ldots,z_p$ such that $0<|z_i|<1$ and $\psi(z_i^{-1})=0$ for $i=1,\ldots,p$.
	Evaluating \eqref{dem:pol_eq} at the $z_i$'s, we get that 
	\begin{align*}
	\Theta(z_i^{-1})b_1(z_i) = 0, \hspace{1cm} \text{for} \hspace{0.5cm} i=1,\ldots,p,
	\end{align*}
	because $H$ has no roots inside the unit circle and $P$ and $Q$ are of finite degrees.
	From the fact that $\deg(b_1)\le p-1$, we also know that for some $z_{i_0}$, $b(z_{i_0})\ne0$ which finally yields 
	$$
	\Theta(z_{i_0}^{-1}) = 0.
	$$
	We therefore obtain that $\psi$ and $\Theta$ have a common root, which is ruled out by assumption, hence the contradiction. The sequence $(k_n)$ can thus be neither bounded nor unbounded, which is absurd.
	We conclude that 
	$$
	m_0 = \sup \{m\ge1: \hspace{0.2cm} \exists k\in\mathbb{Z}, \hspace{0.2cm} d_{k+m}=\ldots=d_{k+1}=0, \hspace{0.2cm} d_{k}\ne0\}<+\infty.
	$$
	Hence the equivalence between $(\iota\iota)$ and $(\iota\iota\iota)$.\\
	
	\noindent Let us now show that whenever $m_0<+\infty$, then \eqref{eq:additionallourd} holds for any $m\ge m_0$.\\

	As $m_0<+\infty$, we have that for any $m\ge m_0$ and $h\ge1$, $\|\bd_k\|>0$ as soon as $\bd_k\ne\boldsymbol{0}$, for all $k\in\mathbb{Z}$ (recall $\bd_k=(d_{k+m},\ldots,d_k,d_{k+1},\ldots,d_{k-h})$).
	For ARMA processes, the non-zero coefficients $d_k$ of the moving average necessarily decay geometrically (times a monomial) as $k\rightarrow\pm\infty$. 
	To fix ideas, say $d_k \underset{k\rightarrow\pm\infty}{\sim} a k^b \lambda^k$, for constants $a\ne0$, $b$ a non-negative integer, and $0<|\lambda|<1$, which may change according to whether $k\rightarrow+\infty$ or $k\rightarrow-\infty$ (if $\deg(\phi)=0$, then $d_{-k}=0$ for $k\ge0$ large enough, however, since we assume $\deg(\psi)\ge1$, it always holds that $|d_k| \underset{k\rightarrow+\infty}{\sim} a k^b \lambda^k$, for the non-zero terms $d_k$).
    Hence, 
    \begin{align*}
    \bd_k \underset{k\rightarrow\pm\infty}{\sim} a k^b \lambda^k \bd_\ast,
    \end{align*}
    for some constant vector $\bd_\ast$ such that $\|\bd_\ast\|>0$ (which may change according to whether $k\rightarrow+\infty$ or $k\rightarrow-\infty$).
    We then have that
    $$
    \dfrac{\|\bd_k\|}{\|\bd_k\|_e} \underset{k\rightarrow\pm\infty}{\longrightarrow} \dfrac{\|\bd_\ast\|}{\|\bd_\ast\|_e} > 0,
    $$
    and
    \begin{align*}
    \|\bd_k\|_e \bigg|\ln\Big(\|\bd_k\|/\|\bd_k\|_e\Big)\bigg| \underset{k\rightarrow\pm\infty}{\sim} \text{const} \hspace{0.1cm} k^b \lambda^k.
    \end{align*}
	Therefore, for any $m\ge m_0$, $h\ge1$,
	$$
	\sum_{k\in\mathbb{Z}}\|\bd_k\|_e \bigg|\ln\Big(\|\bd_k\|/\|\bd_k\|_e\Big)\bigg| < +\infty
	$$

	The equivalence between $(\iota)$ and $(\iota\iota\iota)$ is now clear: on the one hand, if $m_0<+\infty$, then \eqref{eq:additionallourd} holds for all $m\ge m_0$, $h\ge1$, which yields the $(m,h)$-past-representability of $(X_{t-m},\ldots,X_t,X_{t+1},\ldots,X_{t+h})$ for any $m\ge m_0$, $h\ge1$, by Lemma \ref{le:ma_representable}.
	In particular, $(X_t)$ is past-representable.
	On the other hand, assuming $(X_t)$  is past-representable, then necessarily $m_0<+\infty$. \\ 
	
	Regarding the last statement, it follows from the above proof that the condition $m_0<+\infty$ and $m\ge m_0$ is sufficient for $(m,h)$-past-representability.
	It is also necessary, as \eqref{eq:ma_representable} never holds with $m<m_0$ (\textit{a fortiori}, with $m<m_0=+\infty$), concluding the proof.
	
	\subsection{Proof of Lemma \ref{le:spec_repres_agg}}
	
	Denote $\bX_{j,t}=(X_{j,t-m},\ldots,X_{j,t},X_{j,t+1},\ldots,X_{j,t+h})$ the paths of the moving averages $(X_{j,t})$, for $j=1,\ldots,J$.
	The $\boldsymbol{X}_{j,t}$'s are independent $\alpha$-stable random vectors with spectral representations $(\Gamma_j,\bmu^0_j)$ of the form \eqref{def:spectral}.
	We consider only the more delicate case $\alpha=1$ and $\beta_j\in[-1,1]$ for $j=1,\ldots,J$.
	Because of the independence between $\boldsymbol{X}_{1,t},\ldots,\boldsymbol{X}_{J,t}$, we have with $a=2/\pi$
	\begin{align*}
		& \mathbb{E}\Big[e^{i\scal{\bu,\bXt}}\Big] = \mathbb{E}\Big[e^{i\scal{\bu,\sum_{j=1}^J \pi_j\boldsymbol{X}_{j,t}}}\Big] = \prod_{j=1}^J \mathbb{E}\Big[e^{i\scal{\pi_j\bu,\boldsymbol{X}_{j,t}}}\Big] \\
		& \hspace{0.5cm} = \prod_{j=1}^J \exp\Bigg\{-\int_{S_{m+h+1}}\Big(|\scal{\pi_j\bu,\bs}|+ia\scal{\pi_j\bu,\bs}\ln|\scal{\pi_j\bu,\bs}|\Big)\Gamma_j(d\bs)+i\scal{ \pi_j\bu,\boldsymbol{\mu_j}^0}\Bigg\}\\
		& \hspace{0.5cm} =  \exp\Bigg\{-\int_{S_{m+h+1}}\Big(|\scal{\bu,\bs}|+ia\scal{\bu,\bs}\ln|\scal{\bu,\bs}|\Big)\sum_{j=1}^J \pi_j\Gamma_j(d\bs)\\
		& \hspace{5cm}+i\sum_{j=1}^J\bigg(\scal{ \bu,\pi_j\boldsymbol{\mu_j}^0}-a\pi_j\ln\pi_j\int_{S_{m+h+1}}\scal{\bu,\bs}\Gamma_j(d\bs)\bigg)\Bigg\}.
	\end{align*}
	Focusing on the shift vector, we have
	\begin{align*}
	\sum_{j=1}^J\bigg(\scal{ \bu,\pi_j\boldsymbol{\mu_j}^0}-a\pi_j\ln\pi_j\int_{S_{m+h+1}}\scal{\bu,\bs}\Gamma_j(d\bs)\bigg) = \scal{\bu,\sum_{j=1}^J\pi_j(\bmu_j^0 - a \ln\pi_j\tilde{\bmu}_{j})},
	\end{align*}
    with $\tilde{\bmu}_{j}=(\tilde{\mu}_{j,\ell})$ and $\tilde{\mu}_{j,\ell}=\int_{S_{m+h+1}}s_\ell\Gamma_j(d\bs)$, $\ell=-m,\ldots,0,1,\ldots,h$.
    Using the form of $\Gamma_j$ in \eqref{def:spectral}, i.e., $\Gamma_j=\sum_{\vartheta\in S_1}\sum_{k\in\mathbb{Z}}w_{j,\vartheta}\|\bd_{j,k}\|_e \delta_{\left\{\frac{\vartheta \bd_{j,k}}{|\bd_{j,k}\|_e}\right\}}$, we get
    \begin{align*}
    \tilde{\mu}_{j,\ell} & = \int_{S_{m+h+1}}s_\ell\Gamma_j(d\bs)  = \sum_{\vartheta\in S_1}\sum_{k\in\mathbb{Z}}w_{j,\vartheta}\|\bd_{j,k}\|_e \dfrac{\vartheta d_{j,k+\ell}}{\|\bd_{j,k}\|_e}
     = \beta_j\sum_{k\in\mathbb{Z}}d_{j,k+\ell}, \hspace{0.3cm} \ell=-m,\ldots,h.
    \end{align*}
    Hence, $\tilde{\bmu}_{j}=\beta_j\sum_{k\in\mathbb{Z}}\bd_{j,k}$, and using the form of $\bmu_j^0$ as given in \eqref{def:spectral},
    \begin{align*}
    \sum_{j=1}^J\pi_j(\bmu_j^0 - a \ln\pi_j\tilde{\bmu}_{j}) & =  \sum_{j=1}^J\pi_j \bigg(\beta_j\sum_{k\in\mathbb{Z}}\bd_{j,k}\ln\|\bd_{j,k}\|_e -a \ln\pi_j \beta_j\sum_{k\in\mathbb{Z}}\bd_{j,k}  \bigg)\\
    & = -a\sum_{j=1}^J\sum_{k\in\mathbb{Z}}\pi_j\beta_j\bd_{j,k}\ln\|\pi_j\bd_{j,k}\|_e\\
    & := \bmu^0.
    \end{align*}
    Therefore,
    \begin{align*}
    \mathbb{E}\Big[e^{i\scal{\bu,\bXt}}\Big] = \exp\Bigg\{-\int_{S_{m+h+1}}\Big(|\scal{\bu,\bs}|+ia\scal{\bu,\bs}\ln|\scal{\bu,\bs}|\Big)\sum_{j=1}^J \pi_j\Gamma_j(d\bs) + i\scal{\bu,\bmu^0}\Bigg\},
    \end{align*}
	and the random vector $\bX_t$ is $1$-stable with spectral measure 
	\begin{align*}
   \sum_{j=1}^J \pi_j\Gamma_j = \sum_{j=1}^J \sum_{\vartheta\in S_1}\sum_{k\in\mathbb{Z}} w_{j,\vartheta}\pi_j^\alpha \|\bd_{j,k}\|_e^\alpha \delta_{\left\{\dfrac{\vartheta\bd_{j,k}}{\|\bd_{j,k}\|_e}\right\}},
	\end{align*}
	by \eqref{def:spectral}, and shift vector as announced in the lemma.
	
	\subsection{Proof of Lemma \ref{le:spec_repres_agg_vec}}
	
	With the usual notations, let the $\bX_{j,t}$'s be the paths of the moving averages $(X_{j,t})$'s and let $\Gamma_j$, $j=1,\ldots,J$, their spectral measures on the Euclidean unit sphere.
	Let $\Gamma$ the spectral measure of $\bX_t$.
	By Lemma \eqref{le:spec_repres_agg}, $\Gamma=\sum_{j=1}^J\pi_j^\alpha \Gamma_j$.
	Thus, by Proposition \ref{prop238}, in the cases where either $\alpha\ne1$ or $\bX_t$ is symmetric, the vector $\bX_t$ is representable on $C_{m+h+1}^{\|\cdot\|}$ if and only if 
    \begin{align*}
    \Gamma(K^{\|\cdot\|})=0 & \hspace{0.3cm} \iff \hspace{0.3cm} \sum_{j=1}^J\pi_j^\alpha\Gamma_j(K^{\|\cdot\|})=0\\
    & \hspace{0.3cm} \iff \hspace{0.3cm} \Gamma_j(K^{\|\cdot\|})=0, \hspace{0.2cm}\forall \hspace{0.1cm} j =1,\ldots,J.
    \end{align*}
	Given that the $\Gamma_j$'s are the spectral measures of paths of non-aggregated moving averages, it has been shown in the proof of Lemma \ref{le:ma_representable} $\Gamma_j(K^{\|\cdot\|})$ if and only if \eqref{eq:ma_representable} holds for $m$ and the sequence $(d_{j,k})_k$.
	The conclusion in that case follows. 
	The case $\alpha\ne1$ and $\bX_t$ asymmetric is similar.
	
	\subsection{Proof of Proposition \ref{prop:aggma}}
	
	If $\alpha\ne1$, we have by Lemma \ref{le:ma_representable} and the proof of Proposition \ref{prop:ma_pastrepres},
	\begin{align*}
	(X_t) \hspace{0.1cm} \text{past-representable}	& \hspace{0.1cm} \iff \hspace{0.1cm} \exists \hspace{0.1cm} m\ge0,\hspace{0.1cm}\text{\eqref{eq:ma_representable} holds with}\hspace{0.1cm}m\hspace{0.1cm}\text{for all sequences}\hspace{0.1cm}(d_{j,k})_k\\
	& \hspace{0.1cm} \iff \hspace{0.1cm} \forall \hspace{0.1cm} j=1,\ldots,J,\hspace{0.1cm}m_{0,j}<+\infty\\
	& \hspace{0.1cm} \iff \hspace{0.1cm} \forall \hspace{0.1cm} j=1,\ldots,J,\hspace{0.1cm}(X_{j,t})\hspace{0.1cm}\text{past-representable}.
	\end{align*}
	For a given series $(d_{j,k})_k$, \eqref{eq:ma_representable} holds with $m\ge m_{0,j}$ and does not hold with $m<m_{0,j}$.
	Regarding the last statement, we know that for $(X_t)$ $(m,h)$-past-representable, \eqref{eq:ma_representable} holds with the same $m$ for all the sequences $(d_{j,k})_k$, $j=1,\ldots,J$. 
	This holds if $m\ge \max\limits_{j}m_{0,j}$ and cannot hold for if $m<\max\limits_{j}m_{0,j}$.\\
	
	In the case where $\alpha=1$, again by Lemma \ref{le:ma_representable} and denoting generically by $\bX_{t}$ a vector $(X_{t-m},\ldots,X_t,X_{t+1},\ldots,X_{t+h})$ of size $m+h+1$,
	\begin{align*}
		& (X_t) \hspace{0.1cm} \text{past-representable}\\
		& \hspace{0.1cm} \iff \hspace{0.1cm} \exists \hspace{0.1cm} m\ge0,h\ge1,\hspace{0.1cm}\left\{\begin{array}{l}
		     \bX_{t}\hspace{0.1cm}S1S\hspace{0.1cm} \text{and}\hspace{0.1cm} \text{\eqref{eq:ma_representable} holds with}\hspace{0.1cm}m\hspace{0.1cm}\text{for all sequences}\hspace{0.1cm}(d_{j,k})_k \\
		     \text{or}\\
		     \bX_{t}\hspace{0.1cm}\text{asymmetric}\hspace{0.1cm} \text{and}\hspace{0.1cm} \text{\eqref{eq:ma_representable}-\eqref{eq:additionallourd} hold with}\hspace{0.1cm}m,h\hspace{0.1cm}\text{for all sequences}\hspace{0.1cm}(d_{j,k})_k
		\end{array}\right.\\
		& \hspace{0.1cm} \iff \hspace{0.1cm}  \forall\hspace{0.1cm}j=1,\ldots,J,\hspace{0.1cm}m_{0,j}<+\infty, \hspace{0.1cm} \text{and}\hspace{0.1cm}\exists\hspace{0.1cm} m\ge0,h\ge1,\hspace{0.1cm}\left\{\begin{array}{l}
		     \bX_{t}\hspace{0.1cm}S1S\\
		     \text{or}\\
		     \bX_{t}\hspace{0.1cm}\text{asymmetric}\hspace{0.1cm} \text{and}\hspace{0.1cm} \text{\eqref{eq:additionallourd} hold}\\
		     \hspace{0cm}\text{with}\hspace{0.1cm}m,h\hspace{0.1cm}\text{for all sequences}\hspace{0.1cm}(d_{j,k})_k
		\end{array}\right.
	\end{align*}
	We conclude again by noting that the necessary condition \eqref{eq:ma_representable} holds for $m\ge \max\limits_{j}m_{0,j}$ and is violated for $m< \max\limits_{j}m_{0,j}$.
	
	\subsection{Proof of Corollary \ref{cor:agg_arma}}
	
	The equivalence between $(\iota\iota)$ and $(\iota\iota\iota)$ follows from Corollary \eqref{cor:mar}. 
	From the proof of Corollary \eqref{cor:mar}, we also know that, for any $j$, if $m_{0,j}<+\infty$, then \eqref{eq:additionallourd} holds for the sequence $(d_{j,k})_k$ for any $m\ge m_{0,j}$.
	Hence,
	\begin{align*}
	\sup_j m_{0,j} <+\infty & \hspace{0.1cm} \Longrightarrow \hspace{0.1cm} \text{\eqref{eq:additionallourd} holds for any sequence} \hspace{0.1cm} (d_{j,k})_k \hspace{0.1cm} \text{for any} \hspace{0.1cm} m\ge m_{0,j}\\
	& \hspace{0.1cm} \Longrightarrow \hspace{0.1cm} \text{\eqref{eq:additionallourd} holds for any sequence} \hspace{0.1cm} (d_{j,k})_k \hspace{0.1cm} \text{for any} \hspace{0.1cm} m\ge\max_j m_{0,j}m_{0,j}
	\end{align*}
	Thus, $(\iota\iota\iota)$ implies $(\iota)$. The reciprocal is clear.\\
	
	Regarding the last statement, notice that $(X_t)$ if $(m,h)$-past-representable for some $m<\max\limits_j m_{0,j}$, there would then exists some $j$ such that $m<m_{0,j}$. 
	Hence, \eqref{eq:ma_representable} does not holds with $m$ for some particular sequence $(d_{j,k})_k$, which is impossible by Lemma \ref{le:spec_repres_agg_vec}.
	
	\subsection{Proof of Proposition \ref{prop:relativprob}}
	
	By Proposition \ref{prop:cond_tail}
	\begin{align*}
		\mathbb{P}_x^{\|\cdot\|}(\bXt,A|B) \underset{x\rightarrow+\infty}{\longrightarrow} \dfrac{\Gamma^{\|\cdot\|}(A\cap B(V))}{\Gamma^{\|\cdot\|}(B(V))}.
	\end{align*}
	The conclusion follows by considering the points of $B(V)$ and $A\cap B(V)$ that are charged by the spectral measure $\Gamma^{\|\cdot\|}$ in \eqref{eq:specmeasu}.
	
	\subsection{Proof of Lemma \ref{le:djk_aggar1}}
	
	By Proposition \ref{prop:aggma}, we have
	\begin{align*}
		\Gamma^{\|\cdot\|} & =  \sum_{j=1}^J \sum_{\vartheta\in S_1} \sum_{k\in\mathbb{Z}} w_{j,\vartheta}\pi_j^\alpha\|\bd_{j,k}\|^\alpha \delta_{\left\{\frac{\vartheta \bd_{j,k}}{\|\bd_{j,k}\|}\right\}},
	\end{align*}
	with $\bd_{j,k}=(\rho_j^{k+m}\mathds{1}_{\{k+m\ge0\}},\ldots,\rho_j^{k-h}\mathds{1}_{\{k-h\ge 0\}})$ for any $j=1,\ldots,J$ and $k\in\mathbb{Z}$. 
	Thus, for any $j\in\{1,\ldots,J\}$
	\begin{align*}
		\bd_{j,k} & =  \left\{
		\begin{array}{ll}
			\boldsymbol{0}, & \text{if} \quad k \le -m -1,\\ (\rho_j^{k+m},\ldots,\rho_{j},1,0,\ldots,0), & \text{if} \quad -m\le k \le h,\\
			\rho_j^{k-h} \bd_{j,h}, & \text{if} \quad k\ge h.\\
		\end{array}
		\right.
	\end{align*}
	Therefore,
	\begin{align*}
		\Gamma^{\|\cdot\|} & =  \sum_{j=1}^J\sum_{\vartheta\in S_1} w_{j,\vartheta} \pi_j^\alpha  \Bigg[\sum_{k=-m}^{h-1}\|\bd_{j,k}\|^\alpha \delta_{\left\{\frac{\vartheta \bd_{j,k}}{\|\bd_{j,k}\|}\right\}} + \sum_{k=h}^{+\infty}|\rho_j|^{\alpha(k-h)} \|\bd_{j,h}\|^\alpha \delta_{\left\{\frac{\vartheta \rho_j^{k-h} \bd_{j,h}}{|\rho_j|^{k-h}\|\bd_{j,h}\|}\right\}} \Bigg].
	\end{align*}
	Moreover,
	\begin{align*}
		&  \sum_{j=1}^J \sum_{\vartheta\in S_1} w_{j,\vartheta}\pi_j^\alpha  \sum_{k=h}^{+\infty}|\rho_j|^{\alpha (k-h)}\|\bd_{j,h}\|^\alpha \delta_{\left\{\sign(\rho_j)^{k-h}\frac{\vartheta\bd_{j,h}}{\|\bd_{j,h}\|}\right\}} \\
		& \hspace{3cm} =  \sum_{j=1}^J\sum_{\vartheta\in S_1} \pi_j^\alpha  \|\bd_{j,h}\|^\alpha \dfrac{1}{2}\Bigg[\sum_{k=h}^{+\infty}|\rho_j|^{\alpha (k-h)} + \vartheta \beta_j \sum_{k=h}^{+\infty}(\rho_j^{<\alpha>})^{k-h}\Bigg] \delta_{\left\{\frac{\vartheta \bd_{j,h}}{\|\bd_{j,h}\|}\right\}}\\
		& \hspace{3cm} =  \sum_{j=1}^J\sum_{\vartheta\in S_1} \pi_j^\alpha \dfrac{1}{1-|\rho_j|^\alpha} \|\bd_{j,h}\|^\alpha \Bar{w}_{j,\vartheta} \delta_{\left\{\frac{\vartheta \, \bd_{j,h}}{\|\bd_{j,h}\|}\right\}}.
	\end{align*}
	Finally, noticing that for $k=-m$ and any $j\in\{1,\ldots,J\}$, $\bd_{j,k}=(1,0,\ldots,0)$,
	\begin{align*}
		\Gamma^{\|\cdot\|} & =  \sum_{j=1}^J\sum_{\vartheta\in S_1}  \pi_j^\alpha  \Bigg[ w_{j,\vartheta}\sum_{k=-m}^{h-1}\|\bd_{j,k}\|^\alpha \delta_{\left\{\dfrac{\vartheta \bd_{j,k}}{\|\bd_{j,k}\|}\right\}} + \dfrac{\Bar{w}_{j,\vartheta}}{1-|\rho_j|^\alpha} \|\bd_{j,h}\|^\alpha  \delta_{\left\{\dfrac{\vartheta \, \bd_{j,h}}{\|\bd_{j,h}\|}\right\}}\Bigg]\\
		 & =  \sum_{j=1}^J\sum_{\vartheta\in S_1}  \pi_j^\alpha  \Bigg[ w_{j,\vartheta}\bigg(\delta_{\{(\vartheta,0,\ldots,0)\}}+\sum_{k=-m+1}^{h-1}\|\bd_{j,k}\|^\alpha \delta_{\left\{\dfrac{\vartheta \bd_{j,k}}{\|\bd_{j,k}\|}\right\}}\bigg)  + \dfrac{\Bar{w}_{j,\vartheta}}{1-|\rho_j|^\alpha} \|\bd_{j,h}\|^\alpha  \delta_{\left\{\dfrac{\vartheta \, \bd_{j,h}}{\|\bd_{j,h}\|}\right\}}\Bigg]\\
		 & = \sum_{\vartheta\in S_1} \Bigg[w_\vartheta\delta_{\{(\vartheta,0,\ldots,0)\}} + \sum_{j=1}^J  \pi_j^\alpha  \bigg( w_{j,\vartheta}\sum_{k=-m+1}^{h-1}\|\bd_{j,k}\|^\alpha \delta_{\left\{\dfrac{\vartheta \bd_{j,k}}{\|\bd_{j,k}\|}\right\}}  + \dfrac{\Bar{w}_{j,\vartheta}}{1-|\rho_j|^\alpha} \|\bd_{j,h}\|^\alpha  \delta_{\left\{\dfrac{\vartheta \, \bd_{j,h}}{\|\bd_{j,h}\|}\right\}}\bigg)\Bigg].
	\end{align*}
	
	
	\subsection{Proof of Proposition \ref{prop:aggar1_pred}}
	
	\begin{lem}
		\label{le:prepropaggar1}
		Let $\Gamma^{\|\cdot\|}$ be the spectral measure given in Lemma \ref{le:djk_aggar1} and assume that the $\rho_j$'s are all positive. \\
		Letting $(\vartheta_0,j_0,k_0)\in\mathcal{I}$, consider
		\begin{align*}
			I_0 := \left\{\dfrac{\vartheta' \bd_{j',k'}}{\|\bd_{j',k'}\|}:\hspace{0.3cm} \dfrac{\vartheta' f(\bd_{j',k'})}{\|\bd_{j',k'}\|}=\dfrac{\vartheta_0 f(\bd_{j_0,k_0})}{\|\bd_{j_0,k_0}\|} \hspace{0.2cm}\text{for} \hspace{0.2cm}(\vartheta',j',k')\in\mathcal{I}\right\}.
		\end{align*}
		For $m\ge1$, and $0 \le k_0 \le h$, then
		\begin{align*}
			I_0 = \left\{ \dfrac{\vartheta_0 \bd_{j_0,k'}}{\|\bd_{j_0,k'}\|}:\hspace{0.3cm} 0\le k'\le h\right\}. 
		\end{align*}	
		For $m\ge1$, and $-m \le k_0 \le -1$, then
	    \begin{align*}
			I_0 = \left\{ 
            \begin{array}{ll}
                \left\{ \dfrac{\vartheta_0 \bd_{j_0,k_0}}{\|\bd_{j_0,k_0}\|} \right\},  & \text{if} \hspace{0.3cm} -m+1\le k_0 \le -1 \\
                & \\
                \left\{ \dfrac{\vartheta_0 \bd_{0,k_0}}{\|\bd_{0,k_0}\|} \right\} = \left\{ (\vartheta_0,0,\ldots,0)\right\}, & \text{if} \hspace{0.3cm} k_0=-m.
            \end{array}			
			\right.
		\end{align*}
		For $m=0$, then
		\begin{align*}
			I_0 = \left\{\dfrac{\vartheta_0 \bd_{j',k'}}{\|\bd_{j',k'}\|}:\hspace{0.3cm} (j',k')\in\{1,\ldots,J\}\times\{1,\ldots,h\}\cup\{(0,0)\}\right\}. 
		\end{align*}
	\end{lem}
	\textit{Proof.}\\
	\textbf{Case $\boldsymbol{m\ge1}$ and $\boldsymbol{k_0\in\{0,\ldots,h\}}$}\\
	If $k' \in\{-m,\ldots,-1\}$, 
	the $(m+1)$\textsuperscript{th} component of $f(\bd_{j',k'})$ is zero, whereas the $(m+1)$\textsuperscript{th} component of $f(\bd_{j_0,k_0})$ is $\rho_{j_0}^{k_0}\ne0$.
	Necessarily, $\vartheta' f(\bd_{j',k'})/\|\bd_{j',k'}\|\ne\vartheta_0 f(\bd_{j_0,k_0})/\|\bd_{j_0,k_0}\|$ and
	$$
	I_0 = \left\{\dfrac{\vartheta' \bd_{j',k'}}{\|\bd_{j',k'}\|}:\hspace{0.3cm} \dfrac{\vartheta' f(\bd_{j',k'})}{\|\bd_{j',k'}\|}=\dfrac{\vartheta_0 f(\bd_{j_0,k_0})}{\|\bd_{j_0,k_0}\|} \hspace{0.2cm}\text{for} \hspace{0.2cm}(\vartheta',j',k')\in\{-1,+1\}\times\{1,\ldots,J\}\times\{0,\ldots,h\}\right\}.
	$$
	Now, with $k'\in\{0,\ldots,h\}$, we have that
	\begin{align*}
	f(\bd_{j',k'})& =(\rho_{j'}^{k'+m},\ldots,\rho_{j'}^{k'+1},\rho_{j'}^{k'}),\\
	f(\bd_{j_0,k_0})&=(\rho_{j_0}^{k_0+m},\ldots,\rho_{j_0}^{k_0+1},\rho_{j_0}^{k_0}),
	\end{align*}
	and by \eqref{asu:semin} we also have that
	\begin{align*}
	\|\bd_{j',k'}\|&=\|(\rho_{j'}^{k'+m},\ldots,\rho_{j'}^{k'+1},\rho_{j'}^{k'},\overbrace{0,\ldots,0}^h)\|,\\
	\|\bd_{j_0,k_0}\|&=\|(\rho_{j_0}^{k_0+m},\ldots,\rho_{j_0}^{k_0+1},\rho_{j_0}^{k_0},\underbrace{0,\ldots,0}_h)\|.
	\end{align*}
	Thus,
	\begin{align*}
		& \dfrac{\vartheta' f( \bd_{j',k'})}{\|\bd_{j',k'}\|}=\dfrac{\vartheta_0 f(\bd_{j_0,k_0})}{\|\bd_{j_0,k_0}\|}\\
		& \hspace{3cm} \iff  \dfrac{\vartheta'\rho_{j'}^{k'}f\left( \bd_{j',0}\right)}{|\rho_{j'}|^{k'}\|\bd_{j',0}\|}=\dfrac{\vartheta_0\rho_{j_0}^{k_0}f\left(\bd_{j_0,0}\right)}{|\rho_{j_0}|^{k_0}\|\bd_{j_0,0}\|}\\
		& \hspace{3cm} \iff  \dfrac{\vartheta'\rho_{j'}^{\ell}}{\|\bd_{j',0}\|} =\dfrac{\vartheta_0\rho_{j_0}^{\ell}}{\|\bd_{j_0,0}\|}, \hspace{0.3cm} \hspace{0.3cm} \ell=0,\ldots,m\\
		& \hspace{3cm} \iff  \vartheta'\vartheta_0\dfrac{\|\bd_{j_0,0}\|}{\|\bd_{j',0}\|} =\left(\dfrac{\rho_{j_0}}{\rho_{j'}}\right)^{\ell},  \hspace{0.3cm} \ell=0,\ldots,m\\
		& \hspace{3cm} \iff \rho_{j'}=\rho_{j_0}  \hspace{0.3cm} \text{and}  \hspace{0.3cm} \vartheta'\vartheta_0 = 1\\
		& \hspace{3cm} \iff j'=j_0 \hspace{0.3cm} \text{and}  \hspace{0.3cm} \vartheta' = \vartheta_0,
	\end{align*}
	because the $\rho_j$'s are assumed to be non-zero and distinct.\\ 
	\textbf{Case $\boldsymbol{m\ge1}$ and $\boldsymbol{k_0\in\{-m,\ldots,-1\}}$}\\
	By comparing the place of the first zero component, it is easy to see that
	$$
	\dfrac{\vartheta' f(\bd_{j',k'})}{\|\bd_{j',k'}\|} = \dfrac{\vartheta_0 f(\bd_{j_0,k_0})}{\|\bd_{j_0,k_0}\|} \Longrightarrow k'=k_0.
	$$
	\begin{align*}
	f(\bd_{j',k'})& =(\overbrace{\rho_{j'}^{k'+m},\ldots,\rho_{j'},1,0,\ldots,0}^{m+1},\overbrace{0,\ldots,0}^h),\\
	f(\bd_{j_0,k_0})&=(\underbrace{\rho_{j_0}^{k_0+m},\ldots,\rho_{j_0},1,0,\ldots,0}_{m+1},\underbrace{0,\ldots,0}_h),
	\end{align*}
	and we also have that
	\begin{align*}
	\|\bd_{j',k'}\|&=\|(\overbrace{\rho_{j'}^{k'+m},\ldots,\rho_{j'},1,0,\ldots,0}^{m+1},\overbrace{0,\ldots,0}^h)\|,\\
	\|\bd_{j_0,k_0}\|&=\|(\underbrace{\rho_{j_0}^{k_0+m},\ldots,\rho_{j_0},1,0,\ldots,0}_{m+1},\underbrace{0,\ldots,0}_h)\|.
	\end{align*}
	As $ k'=k_0\le-1$,
	\begin{align*}
    	& \dfrac{\vartheta' f(\bd_{j',k'})}{\|\bd_{j',k'}\|}=\dfrac{\vartheta_0 f(\bd_{j_0,k_0})}{\|\bd_{j_0,k_0}\|}\\
		& \hspace{2cm} \iff \dfrac{\vartheta'\rho_{j'}^{\ell}}{\|\bd_{j',k_0}\|} =\dfrac{\vartheta_0\rho_{j_0}^{\ell}}{\|\bd_{j_0,k_0}\|}, \hspace{0.3cm} \ell=0,\ldots,m+k_0, \hspace{0.2cm}\text{and}\hspace{0.2cm}k'=k_0\\
		& \hspace{2cm} \iff \vartheta'\vartheta_0 \dfrac{\|\bd_{j_0,k_0}\|}{\|\bd_{j',k_0}\|}=\left(\dfrac{\rho_{j_0}}{\rho_{j'}}\right)^{\ell},  \hspace{0.3cm} \ell=0,\ldots,m+k_0, \hspace{0.2cm}\text{and}\hspace{0.2cm}k'=k_0. 
	\end{align*}
	Now if $-m+1\le k_0 \le -1$, 
	\begin{align*}
	& \vartheta'\vartheta_0 \dfrac{\|\bd_{j_0,k_0}\|}{\|\bd_{j',k_0}\|}=\left(\dfrac{\rho_{j_0}}{\rho_{j'}}\right)^{\ell},  \hspace{0.3cm} \ell=0,1,\ldots,m+k_0, \hspace{0.2cm}\text{and}\hspace{0.2cm}k'=k_0\\
	& \hspace{5cm} \iff \vartheta'=\vartheta_0 \hspace{0.3cm} \text{and}  \hspace{0.3cm} j'=j_0 \hspace{0.3cm} \text{and}  \hspace{0.3cm}k'=k_0.
	\end{align*}
	If $k_0=-m$, given that $(\vartheta_0,j_0,k_0)\in\mathcal{I}=S_1\times\Big(\{1,\ldots,J\}\times\{-m,\ldots,-1,0,1,\ldots,h\}\cup \{(0,-m)\}\Big)$, then necessarily $j_0=0$. 
	Furthermore, as $k'=k_0=-m$, we similarly have that $j'=j_0=0$ and thus $\bd_{j',k_0}=\bd_{j_0,k_0}=\bd_{0,-m}=(1,0,\ldots,0)$.
	Hence
	\begin{align*}
	& \vartheta'\vartheta_0 \dfrac{\|\bd_{j_0,k_0}\|}{\|\bd_{j',k_0}\|}=\left(\dfrac{\rho_{j_0}}{\rho_{j'}}\right)^{\ell},  \hspace{0.3cm} \ell=0, \hspace{0.2cm}\text{and}\hspace{0.2cm}k'=k_0=-m \hspace{0.2cm}\text{and}\hspace{0.2cm}j'=j_0=0,\\
	& \hspace{2cm} \iff \vartheta' =\vartheta_0 \hspace{0.2cm}\text{and}\hspace{0.2cm}k'=k_0=-m \hspace{0.2cm}\text{and}\hspace{0.2cm} j'=j_0=0
	\end{align*}
	\textbf{Case $\boldsymbol{m=0}$}\\
	If $k_0\in\{1,\ldots,h\}$ then $f(\bd_{j_0,k_0})=\rho_{j_0}^{k_0}$ and by \eqref{asu:semin}, $\|\bd_{j_0,k_0}\|=|\rho_{j_0}|^{k_0}$.
	Thus, $\vartheta_0 f(\bd_{j_0,k_0})/\|\bd_{j_0,k_0}\|=\vartheta_0$.
	If $k_0=-m=0$, then $j_0=0$ and $f(\bd_{j_0,k_0})=1$ and $\vartheta_0 f(\bd_{j_0,k_0})/\|\bd_{j_0,k_0}\|=\vartheta_0$.
	The same holds for $(\vartheta',j',k')\in \mathcal{I}$ and we obtain that 
	\begin{align*}
		& \dfrac{\vartheta' f(\bd_{j',k'})}{\|\bd_{j',k'}\|}=\dfrac{\vartheta_0 f(\bd_{j_0,k_0})}{\|\bd_{j_0,k_0}\|} \iff \vartheta' = \vartheta_0.
	\end{align*}
	\cqfd
	\noindent Let us now prove Proposition \ref{prop:aggar1_pred}. 
	By Proposition \ref{prop:relativprob},
	\begin{align}\label{dem:ar1agg_ratioprob}
		\mathbb{P}_x^{\|\cdot\|}\Big(\bXt,A_{\vartheta,j,k}\Big|B(V_0)\Big) & \underset{x\rightarrow\infty}{\longrightarrow}
		\dfrac{\Gamma^{\|\cdot\|}\Bigg(\bigg\{\dfrac{\vartheta'\bd_{j',k'}}{\|\bd_{j',k'}\|}\in A_{\vartheta,j,k}: \hspace{0.2cm} \dfrac{\vartheta' f(\bd_{j',k'})}{\|\bd_{j',k'}\|}\in V_0\bigg\}\Bigg)}{\Gamma^{\|\cdot\|}\Bigg(\bigg\{\dfrac{\vartheta'\bd_{j',k'}}{\|\bd_{j',k'}\|}\in C_{m+h+1}^{\|\cdot\|}: \hspace{0.2cm} \dfrac{\vartheta' f(\bd_{j',k'})}{\|\bd_{j',k'}\|}\in V_0\bigg\}\Bigg)}.
	\end{align}
	Focusing on the denominator, we have by \eqref{eq:v0}
	\begin{align*}
    \Gamma^{\|\cdot\|}\Bigg(\bigg\{\dfrac{\vartheta'\bd_{j',k'}}{\|\bd_{j',k'}\|}\in C_{m+h+1}^{\|\cdot\|}: \hspace{0.2cm} \dfrac{\vartheta' f(\bd_{j',k'})}{\|\bd_{j',k'}\|}\in V_0\bigg\}\Bigg) & = \Gamma^{\|\cdot\|}\Bigg(\bigg\{\dfrac{\vartheta'\bd_{j',k'}}{\|\bd_{j',k'}\|}\in C_{m+h+1}^{\|\cdot\|}: \hspace{0.2cm} \dfrac{\vartheta' f(\bd_{j',k'})}{\|\bd_{j',k'}\|} = \dfrac{\vartheta_0 f(\bd_{j_0,k_0})}{\|\bd_{j_0,k_0}\|}\bigg\}\Bigg)
	\end{align*}
	We will now distinguish the cases arising from the application of Lemma \ref{le:prepropaggar1}.
	Recall that we assume for this proposition that the $\rho_j$'s are positive.
	Thus, $\sign(\rho_j)=1$ and $\Bar{\beta}_j = \beta_j \dfrac{1-|\rho_j|^\alpha}{1-\rho_j^{<\alpha>}}=\beta_j$ and $\Bar{w}_{j,\vartheta}=w_{j,\vartheta}$ in \eqref{eq:specaggar1} for all $j$'s and $\vartheta\in\{-1,+1\}$.\\
	\noindent \textbf{Case} $\boldsymbol{m\ge1}$ \textbf{and} $\boldsymbol{0\le k_0 \le h}$\\
	By Lemma \ref{le:prepropaggar1},
	\begin{align*}
    & \Gamma^{\|\cdot\|}\Bigg(\bigg\{\dfrac{\vartheta'\bd_{j',k'}}{\|\bd_{j',k'}\|}\in C_{m+h+1}^{\|\cdot\|}: \hspace{0.2cm} \dfrac{\vartheta' f(\bd_{j',k'})}{\|\bd_{j',k'}\|} = \dfrac{\vartheta_0 f(\bd_{j_0,k_0})}{\|\bd_{j_0,k_0}\|}\bigg\}\Bigg)\\
    & \hspace{3cm} = \Gamma^{\|\cdot\|}\Bigg(\bigg\{ \dfrac{\vartheta_0 \bd_{j_0,k'}}{\|\bd_{j_0,k'}\|}:\hspace{0.3cm} 0\le k'\le h \bigg\}\Bigg)\\
    & \hspace{3cm} = \pi_{j_0}^\alpha  \Bigg[ w_{{j_0},\vartheta_0}\sum_{k'=0}^{h-1}\|\bd_{{j_0},k'}\|^\alpha  + \dfrac{\Bar{w}_{{j_0},\vartheta_0}}{1-|\rho_j|^\alpha} \|\bd_{{j_0},h}\|^\alpha  \Bigg]
	\end{align*}
	By \eqref{asu:semin}, for $k'\in\{0,1,\ldots,h\}$
	\begin{align*}
    \|\bd_{j_0,k'}\| & = \|(\rho_{j_0}^{k'+m},\ldots,\rho_{j_0}^{k'+1},\rho_{j_0}^{k'},\underbrace{0,\ldots,0}_h)\|\\
    & = |\rho_{j_0}|^{k'-h} \|(\rho_{j_0}^{m+h},\ldots,\rho_{j_0}^{h+1},\rho_{j_0}^{h},\underbrace{0,\ldots,0}_h)\|\\
    & = |\rho_{j_0}|^{k'-h} \|\bd_{j_0,h}\|.
	\end{align*}
	Thus,
	\begin{align*}
	\Gamma^{\|\cdot\|}\Bigg(\bigg\{\dfrac{\vartheta'\bd_{j',k'}}{\|\bd_{j',k'}\|}\in C_{m+h+1}^{\|\cdot\|}: \hspace{0.2cm} \dfrac{\vartheta' f(\bd_{j',k'})}{\|\bd_{j',k'}\|} = \dfrac{\vartheta_0 f(\bd_{j_0,k_0})}{\|\bd_{j_0,k_0}\|}\bigg\}\Bigg) & =  \pi_{j_0}^\alpha w_{{j_0},\vartheta_0} \|\bd_{{j_0},h}\|^\alpha \Bigg[ \sum_{k'=0}^{h-1} \rho_{j_0}^{\alpha(k'-h)} + \dfrac{1}{1-|\rho_j|^\alpha} \Bigg]\\
	&  = \pi_{j_0}^\alpha w_{{j_0},\vartheta_0} \|\bd_{{j_0},h}\|^\alpha \dfrac{|\rho_j|^{-\alpha h}}{1-|\rho_j|^\alpha}.
	\end{align*}
	Similarly for the numerator in \eqref{dem:ar1agg_ratioprob}, by \eqref{eq:avjk},
	\begin{align*}
		& \Gamma^{\|\cdot\|}\Bigg(\bigg\{\dfrac{\vartheta'\bd_{j',k'}}{\|\bd_{j',k'}\|}\in A_{\vartheta,j,k}: \hspace{0.2cm} \dfrac{\vartheta' f(\bd_{j',k'})}{\|\bd_{j',k'}\|}\in V_0\bigg\}\Bigg)\\
		& \hspace{2cm} = \Gamma^{\|\cdot\|}\Bigg(\bigg\{ \dfrac{\vartheta_0 \bd_{j_0,k'}}{\|\bd_{j_0,k'}\|} \in A_{\vartheta,j,k}:\hspace{0.3cm} 0\le k'\le h \bigg\}\Bigg)\\
		& \\
    	& \hspace{2cm} = \left\{\begin{array}{ll}
			\Gamma^{\|\cdot\|}\Bigg(\bigg\{ \dfrac{\vartheta_0 \bd_{j_0,k}}{\|\bd_{j_0,k}\|} \bigg\}\Bigg), &\hspace{0.3cm} \text{if} \hspace{0.3cm} j= j_0 \hspace{0.3cm} \text{and} \hspace{0.3cm} \vartheta=\vartheta_0,\\
			\Gamma^{\|\cdot\|}(\emptyset), &\hspace{0.3cm} \text{if} \hspace{0.3cm} j\ne j_0 \hspace{0.3cm} \text{or} \hspace{0.3cm} \vartheta\ne\vartheta_0,
		\end{array}\right.\\
		& \\
		& \hspace{2cm} = \left\{
		\begin{array}{ll}
			\pi_{j_0}^\alpha w_{j_0,\vartheta_0} \|\bd_{j_0,h}\|^\alpha |\rho_{j_0}|^{\alpha(k-h)}\delta_{\{\vartheta_0\}}(\vartheta)\delta_{\{j_0\}}(j), &\hspace{0.3cm} \text{if} \hspace{0.3cm} 0\le k \le h-1,\\
			\pi_{j_0}^\alpha w_{j_0,\vartheta_0} \|\bd_{j_0,h}\|^\alpha \dfrac{1}{1-|\rho_{j_0}|^\alpha}\delta_{\{\vartheta_0\}}(\vartheta)\delta_{\{j_0\}}(j), &\hspace{0.3cm} \text{if} \hspace{0.3cm} k=h.
		\end{array}
		\right.
	\end{align*}
	The conclusion follows.\\
	\noindent \textbf{Case} $\boldsymbol{m\ge1}$ \textbf{and} $\boldsymbol{-m\le k_0 \le -1}$  \hfill
	
	We have by Lemma \ref{le:prepropaggar1}
	
	\begin{align*}
    & \Gamma^{\|\cdot\|}\Bigg(\bigg\{\dfrac{\vartheta'\bd_{j',k'}}{\|\bd_{j',k'}\|}\in C_{m+h+1}^{\|\cdot\|}: \hspace{0.2cm} \dfrac{\vartheta' f(\bd_{j',k'})}{\|\bd_{j',k'}\|} = \dfrac{\vartheta_0 f(\bd_{j_0,k_0})}{\|\bd_{j_0,k_0}\|}\bigg\}\Bigg) = \Gamma^{\|\cdot\|}\Bigg(\bigg\{ \dfrac{\vartheta_0 \bd_{j_0,k_0}}{\|\bd_{j_0,k_0}\|}\bigg\}\Bigg).
    & \hspace{3cm} = 
	\end{align*}
	If $-m+1\le k_0 \le -1$,
	\begin{align*}
	\Gamma^{\|\cdot\|}\Bigg(\bigg\{ \dfrac{\vartheta_0 \bd_{j_0,k_0}}{\|\bd_{j_0,k_0}\|}\bigg\}\Bigg) = \pi_{j_0}^\alpha w_{j_0,\vartheta_0} \|\bd_{j_0,k_0}\|^\alpha,
	\end{align*}
	and
	\begin{align*}
		& \Gamma^{\|\cdot\|}\Bigg(\bigg\{\dfrac{\vartheta'\bd_{j',k'}}{\|\bd_{j',k'}\|}\in A_{\vartheta,j,k}: \hspace{0.2cm} \dfrac{\vartheta' f(\bd_{j',k'})}{\|\bd_{j',k'}\|}\in V_0\bigg\}\Bigg)\\
		& \hspace{5cm} = \Gamma^{\|\cdot\|}\Bigg(A_{\vartheta,j,k}\cap\bigg\{ \dfrac{\vartheta_0 \bd_{j_0,k_0}}{\|\bd_{j_0,k_0}\|}\bigg\}\Bigg)\\
		& \\
    	& \hspace{5cm} = \left\{\begin{array}{ll}
			\Gamma^{\|\cdot\|}\Bigg( \bigg\{ \dfrac{\vartheta_0 \bd_{j_0,k_0}}{\|\bd_{j_0,k_0}\|}\bigg\}\Bigg), &\hspace{0.2cm} \text{if} \hspace{0.2cm} j= j_0 \hspace{0.2cm} \text{and} \hspace{0.2cm} \vartheta=\vartheta_0, \hspace{0.2cm} \text{and} \hspace{0.2cm} k=k_0,\\
			\Gamma^{\|\cdot\|}(\emptyset), &\hspace{0.2cm} \text{if} \hspace{0.2cm} j\ne j_0 \hspace{0.2cm} \text{or} \hspace{0.2cm} \vartheta\ne\vartheta_0 \hspace{0.2cm} \text{or} \hspace{0.2cm} k\ne k_0,
		\end{array}\right.\\
		& \\
		& \hspace{5cm} = \pi_{j_0}^\alpha w_{j_0,\vartheta_0} \|\bd_{j_0,k_0}\|^\alpha \delta_{\{\vartheta_0\}}(\vartheta)\delta_{\{j_0\}}(j)\delta_{\{k_0\}}(k).
	\end{align*}
	If $k_0=-m$, then $\bd_{j_0,k_0}=\bd_{0,-m}=(1,0,\ldots,0)$, and
	\begin{align*}
	\Gamma^{\|\cdot\|}\Bigg(\bigg\{ \dfrac{\vartheta_0 \bd_{j_0,k_0}}{\|\bd_{j_0,k_0}\|}\bigg\}\Bigg) = \Gamma^{\|\cdot\|}\Big( \{\vartheta_0 (1,0,\ldots,0)\}\Big) = w_{\vartheta_0}, 
	\end{align*}
	and
	\begin{align*}
		& \Gamma^{\|\cdot\|}\Bigg(\bigg\{\dfrac{\vartheta'\bd_{j',k'}}{\|\bd_{j',k'}\|}\in A_{\vartheta,j,k}: \hspace{0.2cm} \dfrac{\vartheta' f(\bd_{j',k'})}{\|\bd_{j',k'}\|}\in V_0\bigg\}\Bigg)\\
		& \hspace{2cm} = \Gamma^{\|\cdot\|}\Bigg(A_{\vartheta,j,k}\cap\bigg\{ \dfrac{\vartheta_0 \bd_{j_0,k_0}}{\|\bd_{j_0,k_0}\|}\bigg\}\Bigg)\\
		& \\
    	& \hspace{2cm} = \left\{\begin{array}{ll}
			\Gamma^{\|\cdot\|}\Big( A_{\vartheta,j,k}\cap\{\vartheta_0 (1,0,\ldots,0)\}\Big), &\hspace{0.2cm} \text{if} \hspace{0.2cm} \vartheta=\vartheta_0, \hspace{0.2cm} \text{and} \hspace{0.2cm} k=k_0=-m, \hspace{0.2cm} \text{and} \hspace{0.2cm} j=j_0=0\\
			\Gamma^{\|\cdot\|}(\emptyset), &\hspace{0.2cm} \text{if} \hspace{0.2cm} \vartheta\ne\vartheta_0 \hspace{0.2cm} \text{or} \hspace{0.2cm} k\ne k_0, \hspace{0.2cm} \text{or} \hspace{0.2cm} j\ne j_0
		\end{array}\right.\\
		& \\
		& \hspace{2cm} =
		w_{\vartheta_0}\delta_{\{\vartheta_0\}}(\vartheta)\delta_{\{j_0\}}(j)\delta_{\{k_0\}}(k).
	\end{align*}
		
	The conclusion follows as previously.\\
	\noindent \textbf{Case} $\boldsymbol{m=0}$ 
	
	By Lemma \ref{le:prepropaggar1}, as the $\rho_j$'s are positive 
	\begin{align*}
	& \Gamma^{\|\cdot\|}\Bigg(\bigg\{\dfrac{\vartheta'\bd_{j',k'}}{\|\bd_{j',k'}\|}\in C_{m+h+1}^{\|\cdot\|}: \hspace{0.2cm} \dfrac{\vartheta' f(\bd_{j',k'})}{\|\bd_{j',k'}\|} = \dfrac{\vartheta_0 f(\bd_{j_0,k_0})}{\|\bd_{j_0,k_0}\|}\bigg\}\Bigg) \\
	& \hspace{2cm} = \Gamma^{\|\cdot\|}\Bigg(\bigg\{\dfrac{\vartheta_0\bd_{j',k'}}{\|\bd_{j',k'}\|}\in C_{m+h+1}^{\|\cdot\|}: \hspace{0.2cm} (j',k')\in\{1,\ldots,J\}\times\{0,\ldots,h\}\cup\{(0,0)\}\bigg\}\Bigg)
	\end{align*}
	Given that $w_{\vartheta_0}=\sum_{j'=1}^J\pi_{j}^\alpha w_{j',\vartheta_0}$ and $\|\bd_{j',k'}\|=|\rho_{j'}|^{k'}$, for any $1 \le j' \le J$, $1\le k' \le h$,
	\begin{align*}
	& \Gamma^{\|\cdot\|}\Bigg(\bigg\{\dfrac{\vartheta'\bd_{j',k'}}{\|\bd_{j',k'}\|}\in C_{m+h+1}^{\|\cdot\|}: \hspace{0.2cm} \dfrac{\vartheta' f(\bd_{j',k'})}{\|\bd_{j',k'}\|} = \dfrac{\vartheta_0 f(\bd_{j_0,k_0})}{\|\bd_{j_0,k_0}\|}\bigg\}\Bigg) \\
	& \hspace{6cm} = w_{\vartheta_0} + \sum_{j'=1}^J\pi_{j}^\alpha w_{j',\vartheta_0} \Bigg[\sum_{k'=1}^{h-1}\|\bd_{j',k'}\|^\alpha + \dfrac{\|\bd_{j',h}\|^\alpha}{1-|\rho_{j'}|^\alpha}\Bigg]\\
	& \hspace{6cm} =  \sum_{j'=1}^J\pi_{j}^\alpha w_{j',\vartheta_0} \Bigg[1+\sum_{k'=1}^{h-1}|\rho_{j'}|^{\alpha k'} + \dfrac{|\rho_{j'}|^{\alpha h}}{1-|\rho_{j'}|^\alpha}\Bigg]\\
	& \hspace{6cm} =  \sum_{j'=1}^J\pi_{j'}^\alpha w_{j',\vartheta_0} \Bigg[\dfrac{1-|\rho_{j'}|^{\alpha h}}{1-|\rho_{j'}|^{\alpha}} + \dfrac{|\rho_{j'}|^{\alpha h}}{1-|\rho_{j'}|^\alpha}\Bigg]\\
	& \hspace{6cm} =  \sum_{j'=1}^J\pi_{j'}^\alpha w_{j',\vartheta_0} \dfrac{1}{1-|\rho_{j'}|^\alpha}.
	\end{align*}
	Similarly, by \eqref{eq:avjk},
	\begin{align*}
	& \Gamma^{\|\cdot\|}\Bigg(\bigg\{\dfrac{\vartheta'\bd_{j',k'}}{\|\bd_{j',k'}\|}\in A_{\vartheta,j,k}: \hspace{0.2cm} \dfrac{\vartheta' f(\bd_{j',k'})}{\|\bd_{j',k'}\|}\in V_0\bigg\}\Bigg)\\
	& \hspace{2cm} = \Gamma^{\|\cdot\|}\Bigg(A_{\vartheta,j,k}\cap\bigg\{\dfrac{\vartheta_0\bd_{j',k'}}{\|\bd_{j',k'}\|}\in C_{m+h+1}^{\|\cdot\|}: \hspace{0.2cm} (j',k')\in\{1,\ldots,J\}\times\{0,\ldots,h\}\cup\{(0,0)\}\bigg\}\Bigg)\\
	& \\
	& \hspace{2cm} = \left\{\begin{array}{ll}
			\Gamma^{\|\cdot\|}\Bigg( \bigg\{ \dfrac{\vartheta_0 \bd_{j,k}}{\|\bd_{j,k}\|}\bigg\}\Bigg), &\hspace{0.2cm} \text{if} \hspace{0.2cm} \vartheta=\vartheta_0,\\
			\Gamma^{\|\cdot\|}(\emptyset), &\hspace{0.2cm} \text{if} \hspace{0.2cm} \vartheta\ne\vartheta_0,
		\end{array}\right.\\
    & \\
	 & \hspace{2cm} = \left\{
		\begin{array}{ll}
		     \sum_{j'=1}^J \pi_{j'}^\alpha w_{j',\vartheta_0} \delta_{\{\vartheta_0\}}(\vartheta), & \quad \text{if} \quad k=0,\\
			\pi_{j}^\alpha w_{j,\vartheta_0}  |\rho_{j}|^{\alpha k}\delta_{\{\vartheta_0\}}(\vartheta), &\quad \text{if} \quad 1\le k \le h-1,\\
			\pi_{j}^\alpha w_{j,\vartheta_0} \dfrac{|\rho_{j}|^{\alpha h}}{1-|\rho_{j}|^{\alpha}}\delta_{\{\vartheta_0\}}(\vartheta), &\quad \text{if} \quad k=h.
		\end{array}
		\right.
	\end{align*}
	The conclusion follows.
	
	\subsection{Proof of Proposition \ref{prop:ar2_pred}}
	
	\begin{lem}\label{lem:k=lar2}
		Let $X_t$ be the $\alpha$-stable anticipative AR(2) (resp. fractionally integrated AR) as in \eqref{def:ar2} (resp. \eqref{def:frac}). 
		With $f$ as in \eqref{def:f}, and for any $m\ge1$, $h\ge0$,
		\begin{align*}
			\forall k,\ell\ge -m, \hspace{0.2cm} \forall \vartheta_1,\vartheta_2\in S_1, \quad \bigg[\dfrac{f(\vartheta_1\bd_k)}{\|\bd_k\|} = \dfrac{f(\vartheta_2 \bd_\ell)}{\|\bd_\ell\|} \hspace{0.2cm} \Longrightarrow \hspace{0.2cm} k=\ell \hspace{0.2cm} \text{and} \hspace{0.2cm} \vartheta_1=\vartheta_2\bigg].
		\end{align*}
	\end{lem}
	\noindent \textit{Proof.}
	
	\noindent The result is clear for both processes for $-m\le k,\ell \le -1$. For  $k,\ell\ge0$,
	\begin{align}
		\dfrac{f(\vartheta_1\bd_k)}{\|\bd_k\|} = \dfrac{f(\vartheta_2\bd_\ell)}{\|\bd_\ell\|}  & \iff  \bigg[\forall\,i=0,\ldots,m, \hspace{0.3cm} \dfrac{\vartheta_1d_{k+i}}{\|\bd_k\|} = \dfrac{\vartheta_2d_{\ell+i}}{\|\bd_\ell\|}\bigg]\nonumber\\
		& \iff  \dfrac{d_{k}}{d_{\ell}} = \dfrac{d_{k+1}}{d_{\ell+1}} = \ldots = \vartheta_1\vartheta_2\dfrac{\|\bd_k\|}{\|\bd_\ell\|}.\label{eq:ratio_eq}
	\end{align}
	The last statement in particular implies that $\dfrac{d_{k}}{d_{\ell}} = \dfrac{d_{k+1}}{d_{\ell+1}}$.\\
	For the anticipative AR(2), if $\lambda_1\ne\lambda_2$, we then have
	\begin{align*}
		\dfrac{d_{k}}{d_{\ell}} = \dfrac{d_{k+1}}{d_{\ell+1}} & \iff  \dfrac{\lambda_1^{k+1}-\lambda_2^{k+1}}{\lambda_1^{\ell+1}-\lambda_2^{\ell+1}} = \dfrac{\lambda_1^{k+2}-\lambda_2^{k+2}}{\lambda_1^{\ell+2}-\lambda_2^{\ell+2}}\\
		& \iff  \lambda_1^{k-\ell} = \lambda_2^{k-\ell}\\
		& \iff  k =\ell.
	\end{align*}
	This case $\lambda_1=\lambda_2=\lambda$ is similar. For the anticipative fractionally integrated AR, given that $\Gamma(z+1)=z\Gamma(z)$ for any $z\in\mathbb{C}$, we have
	\begin{align*}
		\dfrac{d_{k}}{d_{\ell}} = \dfrac{d_{k+1}}{d_{\ell+1}} & \iff \dfrac{\Gamma(k+d)\Gamma(\ell+1)}{\Gamma(\ell+d)\Gamma(k+1)}  = \dfrac{\Gamma(k+d+1)\Gamma(\ell+2)}{\Gamma(\ell+d+1)\Gamma(k+2)}\\
		& \iff \dfrac{\Gamma(\ell+d+1)\Gamma(k+2)}{\Gamma(\ell+d)\Gamma(k+1)}  = \dfrac{\Gamma(k+d+1)\Gamma(\ell+2)}{\Gamma(k+d)\Gamma(\ell+1)}\\
		& \iff (k-\ell)(d-1) = 0\\
		& \iff k=\ell.
	\end{align*}
	Therefore, in all cases,
	\begin{align*}
		\dfrac{d_{k}}{d_{\ell}} = \dfrac{d_{k+1}}{d_{\ell+1}} = \ldots = \vartheta_1\vartheta_2\dfrac{\|\bd_k\|}{\|\bd_\ell\|} & \hspace{0.3cm}  \Longrightarrow \hspace{0.3cm}  k=\ell \hspace{0.2cm} \text{and} \hspace{0.2cm} \vartheta_1\vartheta_2=1.
	\end{align*}
	
	\cqfd
	
	\noindent Let us now prove Proposition \ref{prop:ar2_pred}.
	The spectral measure of $\bX_t$ writes
	\begin{align*}
	\Gamma^{\|\cdot\|} & = \sigma^\alpha\sum_{\vartheta\in S_1}\sum_{k\in\mathbb{Z}} w_\vartheta \|\bd_{k}\|^\alpha \delta_{\left\{\frac{\vartheta\bd_{k}}{\|\bd_{k}\|}\right\}},
	\end{align*}
	where the sequences $(d_k)$ are given respectively by \eqref{eq:ar2dk} and \eqref{eq:fracdk} for the anticipative AR(2) and fractionally integrated processes.
	By Proposition \ref{prop:cond_tail},
	\begin{align*}
		\mathbb{P}_x^{\|\cdot\|}\Big(\bXt,A\Big|B(V_0)\Big) & \underset{x\rightarrow\infty}{\longrightarrow} \dfrac{\Gamma^{\|\cdot\|}(A\cap B(V_0))}{\Gamma^{\|\cdot\|}(B(V_0))}.
	\end{align*}
	On the one hand, we have by definition of $B(V_0)$, $V_0$ and Lemma \ref{lem:k=lar2},
	\begin{align*}
		\Gamma^{\|\cdot\|}(B(V_0)) & = \Gamma^{\|\cdot\|}\Bigg(\bigg\{\dfrac{\vartheta\bd_{k}}{\|\bd_{k}\|}\in B(V_0): \hspace{0.1cm}(\vartheta,k)\in\{-1,+1\}\times\mathbb{Z}\bigg\}\Bigg)\\
		& = \Gamma^{\|\cdot\|}\Bigg(\bigg\{\dfrac{\vartheta\bd_{k}}{\|\bd_{k}\|}\in C_{m+h+1}^{\|\cdot\|}:\hspace{0.1cm}\dfrac{\vartheta f(\bd_{k})}{\|\bd_{k}\|}\in V_0, \hspace{0.1cm}(\vartheta,k)\in\{-1,+1\}\times\mathbb{Z}\bigg\}\Bigg)\\
		& = \Gamma^{\|\cdot\|}\Bigg(\bigg\{\dfrac{\vartheta\bd_{k}}{\|\bd_{k}\|}\in C_{m+h+1}^{\|\cdot\|}:\hspace{0.1cm}\dfrac{\vartheta f(\bd_{k})}{\|\bd_{k}\|} = \dfrac{\vartheta_0f(\bd_{k_0})}{\|\bd_{k_0}\|}, \hspace{0.1cm}(\vartheta,k)\in\{-1,+1\}\times\mathbb{Z}\bigg\}\Bigg)\\
		& = \Gamma^{\|\cdot\|}\Bigg(\bigg\{ \dfrac{\vartheta_0\bd_{k_0}}{\|\bd_{k_0}\|}\bigg\}\Bigg).
	\end{align*}
	Similarly, it is easily shown that
	\begin{align*}
    \Gamma^{\|\cdot\|}(A\cap B(V_0)) & = \Gamma^{\|\cdot\|}\Bigg(A\cap\bigg\{ \dfrac{\vartheta_0\bd_{k_0}}{\|\bd_{k_0}\|}\bigg\}\Bigg).
	\end{align*}
	The conclusion follows.

	\subsection{Proof of Proposition \ref{prop:reprbility_multivar}}
	
	We start with a lemma showing that $\boldsymbol{\underline{X}}_t$ is indeed $S \alpha S$ and providing the form of its spectral measure on the Euclidean unit sphere. 
	The representability condition on the unit cylinder $C_4^{\|\cdot\|}$ will follow. 
	\begin{lem}\label{lem:xunderbar_spec} Let $(\bX_t)$ as in \eqref{def:bivarAR1}.
		Then, the vector $\boldsymbol{\underline{X}}_t=(X_{1,t},X_{2,t},X_{1,t+1},X_{2,t+1})'$ is $S\alpha S$ with zero shift vector and spectral measure given by
		$$
		\Gamma_4 = \Delta + \Gamma_{4,1} + \Gamma_{4,2}.
		$$
		Here,
		\begin{align*}
			\Delta & = \sum_{i=1,2}\dfrac{\sigma_i^\alpha}{2}(1+\rho_i^2)^{\alpha/2} \dfrac{|\rho_i|^\alpha}{1-|\rho_i|^\alpha}(\delta_{\left\{\bx_i/\|\bx_i\|_e\right\}}+\delta_{\left\{-\bx_i/\|\bx_i\|_e\right\}}),
		\end{align*}
		with $\sigma_i^\alpha:=\int_{S_2}|s_i|^\alpha\Gamma_2(d\bs)$, points $\bx_1 = (1,0,\rho_1^{-1},0)$, $\bx_2 = (0,1,0,\rho_2)$,
		\begin{align*}
			\Gamma_{4,1} (d\bs)& = \|B\bs\|_e^{-\alpha}\tilde{\Gamma}_{4,1}\circ T_{B}(d\bs),\\
			\Gamma_{4,2}(d\bs) & = \|C\bs\|_e^{-\alpha}\tilde{\Gamma}_{4,2}\circ T_{C}(d\bs),
		\end{align*}
		with $T_{A}:S_4\rightarrow S_4$ is defined by $T_{A}(\boldsymbol{\underline{s}}) = A\boldsymbol{\underline{s}}/\|A\boldsymbol{\underline{s}}\|_e$, for any invertible matrix $A$ of dimension 4,
		$$
		B = \begin{pmatrix}
		1 & 0 & 0 & 0\\
		0 & 1 & 0 & 0\\
		0 & 0 & 1 & 0\\
		0 & -\rho_2 & 0 & 1
		\end{pmatrix}, \quad \text{and} \quad C = \begin{pmatrix}
		1 & 0 & -\rho_1 & 0\\
		0 & 1 & 0 & 0\\
		0 & 0 & 1 & 0\\
		0 & 0 & 0 & 1
		\end{pmatrix},
		$$
		and $\tilde{\Gamma}_{4,i}(\,\cdot\,)=\Gamma_2 \circ h_i (S_{4,i}\cap \,\cdot\,)$, $i=1,2$, where $S_{4,1}=\{(s_1,s_2,0,0)\in S_4:\hspace{0.2cm} (s_1,s_2)\in S_2\}$, $S_{4,2}=\{(0,0,s_3,s_4)\in S_4:\hspace{0.2cm} (s_3,s_4)\in S_2\}$, and $h_1,h_2: S_4\rightarrow S_2$ are the functions defined by $h_1((s_1,s_2,s_3,s_4))=(s_1,s_2)$ and $h_2((s_1,s_2,s_3,s_4))=(s_3,s_4)$.
	\end{lem}
	\noindent \textit{Proof.}\\
	Let $\boldsymbol{\underline{u}}=(\bu_{0}',\bu_{1}')'\in\mathbb{R}^4$ with $\bu_i=(u_{1,i},u_{2,i})'$, $i=0,1$. 
	The characteristic function of $\boldsymbol{\underline{X}}_t$ reads
	\begin{align*}
		\varphi(\boldsymbol{\underline{u}}) & := \mathbb{E}[\exp\{i\scal{\boldsymbol{\underline{u}},\boldsymbol{\underline{X}}_t}\}]=\mathbb{E}[\exp\{i\sum_{j=0}^1\scal{\bu_j,\bX_{t+j}}\}] = \mathbb{E}[\exp\{i\sum_{k\in\mathbb{Z}}\sum_{j=0}^1\scal{\bu_j,\bA_k\bveps_{t+k+j}}\}]\\
		& = \prod_{k\in\mathbb{Z}}\mathbb{E}[\exp\{i\scal{\sum_{j=0}^1\bA_{k-j}'\bu_j,\bveps_{t+k}}\}],
	\end{align*}
	where for all $k\in\mathbb{Z}$
	$$
	A_k = 
	\begin{pmatrix}
	\rho_1^k \mathds{1}_{\{k\ge0\}} & 0\\
	0 & \rho_2^{-k} \mathds{1}_{\{k\le0\}}
	\end{pmatrix}.
	$$
	Thus,
	\begin{align}
		-\ln\,\varphi (\boldsymbol{\underline{u}}) & = \sum_{k\in\mathbb{Z}}\int_{S_2} |\scal{A_k\bu_0+A_{k-1}\bu_1,\bs}|^\alpha\Gamma_2(d\bs)\nonumber\\
		& = \sum_{k\le-1} \int_{S_2}|\rho_2^{-k}(u_{2,0}+\rho_2u_{2,1})s_2|^\alpha\Gamma_2(d\bs) + \sum_{k\ge2} \int_{S_2}|\rho_1^{k-1}(\rho_1u_{1,0}+u_{1,1})s_1|^\alpha\Gamma_2(d\bs)\nonumber\\
		& \hspace*{0.5cm } + \int_{S_2}|u_{1,0}s_1+(u_{2,0}+\rho_2u_{2,1})s_2|^\alpha\Gamma_2(d\bs) + \int_{S_2}|(\rho_1u_{1,0}+u_{1,1})s_1+u_{2,1}s_2|^\alpha\Gamma_2(d\bs)\nonumber\\
		& = \sigma_2^\alpha \dfrac{|\rho_2|^\alpha}{1-|\rho_2|^\alpha} |u_{2,0}+\rho_2u_{2,1}|^\alpha + \sigma_1^\alpha \dfrac{|\rho_1|^\alpha}{1-|\rho_1|^\alpha}|\rho_1u_{1,0}+u_{1,1}|^\alpha\nonumber\\
		& \hspace*{0.5cm } + \int_{S_2}|u_{1,0}s_1+(u_{2,0}+\rho_2u_{2,1})s_2|^\alpha\Gamma_2(d\bs) + \int_{S_2}|(\rho_1u_{1,0}+u_{1,1})s_1+u_{2,1}s_2|^\alpha\Gamma_2(d\bs),\label{eq:pre_charac_fun}
	\end{align}
	where $\sigma_i^\alpha:=\int_{S_2}|s_i|^\alpha\Gamma_2(d\bs)$, $i=1,2$.
	We notice that the characteristic function of $\boldsymbol{\underline{X}}_t$ is real. 
	Hence, $\boldsymbol{\underline{X}}_t$ being $\alpha$-stable is equivalent to $\boldsymbol{\underline{X}}_t$ being symmetric $\alpha$-stable, and therefore, by Theorem 2.4.3 in \cite{st94}, $\boldsymbol{\underline{X}}_t$ will be $\alpha$-stable if and only if there exists a unique symmetric finite measure $\Gamma_4$ on the Euclidean unit sphere such that
	\begin{equation}\label{dem:charactoshow}
		-\ln\,\varphi (\boldsymbol{\underline{u}}) = \int_{S_4} |\scal{\boldsymbol{\underline{u}},\boldsymbol{\underline{s}}}|\Gamma_4(d\bs).    
	\end{equation}
	We will thus rewrite \eqref{eq:pre_charac_fun} to exhibit such a symmetric measure.
	The two first terms are easily rewritten with charged atoms on $S_4$: for all $\boldsymbol{\underline{u}}\in\mathbb{R}^4$,
	\begin{align}
		\label{eq:mu_meas}
		\sigma_2^\alpha \dfrac{|\rho_2|^\alpha}{1-|\rho_2|^\alpha} |u_{2,0}+\rho_2u_{2,1}|^\alpha + \sigma_1^\alpha \dfrac{|\rho_1|^\alpha}{1-|\rho_1|^\alpha}|\rho_1u_{1,0}+u_{1,1}|^\alpha = \int_{S_4}|\scal{\boldsymbol{\underline{u}},\boldsymbol{\underline{s}}}|^\alpha \Delta(d\boldsymbol{\underline{s}}),
	\end{align}
	where $\Delta = \sum_{i=1,2}\dfrac{\sigma_i^\alpha}{2}(1+\rho_i^2)^{\alpha/2} \dfrac{|\rho_i|^\alpha}{1-|\rho_i|^\alpha}(\delta_{\left\{\bx_i/\|\bx_i\|_e\right\}}+\delta_{\left\{-\bx_i/\|\bx_i\|_e\right\}})$.
	The third and fourth terms in \eqref{eq:pre_charac_fun} can also be rewritten as integral over $S_4$. Starting with the third term, notice that the integral over $S_2$ can be seen as an integral over $S_4$ with a spectral measure $\tilde{\Gamma}_{4,1}$ coinciding with $\Gamma_2$ on $S_{4,1}=\{(s_1,s_2,0,0)\in S_4:\hspace{0.2cm} (s_1,s_2)\in S_2\}$ and having zero mass outside:
	\begin{align*}
		\int_{S_2}|u_{1,0}s_1+(u_{2,0}+\rho_2u_{2,1})s_2|^\alpha\Gamma_2(d\bs) & = \int_{S_4}|u_{1,0}s_1+(u_{2,0}+\rho_2u_{2,1})s_2 + u_{1,1}s_{3}+u_{2,1}s_4|^\alpha\tilde{\Gamma}_{4,1}(d\boldsymbol{\underline{s}}),
	\end{align*}
	with $\tilde{\Gamma}_{4,1}(\,\cdot\,)=\Gamma_2 \circ h_1 (S_{4,1}\cap \,\cdot\,)$, where $h_1:S_4\rightarrow S_2$ is the function defined by $h_1((s_1,s_2,s_3,s_4))=(s_1,s_2)$.
	Thus,
	\begin{align}
		\int_{S_2}|u_{1,0}s_1+(u_{2,0}+\rho_2u_{2,1})s_2|^\alpha\Gamma_2(d\bs) & = \int_{S_4}|\scal{\bb\boldsymbol{\underline{u}},\boldsymbol{\underline{s}}}|^\alpha\tilde{\Gamma}_{4,1} (d\boldsymbol{\underline{s}})\nonumber\\ & = \int_{S_4}|\scal{\boldsymbol{\underline{u}},\bb'\boldsymbol{\underline{s}}}|^\alpha\tilde{\Gamma}_{4,1} (d\boldsymbol{\underline{s}}) , \label{eq:third_term}
	\end{align}
	with 
	$$
	\bb = 
	\begin{pmatrix}
	1 & 0 & 0 & 0\\
	0 & 1 & 0 & \rho_2\\
	0 & 0 & 1 & 0\\
	0 & 0 & 0 & 1
	\end{pmatrix}.
	$$
	As the matrix $\bb'$ is invertible and $B=\bb'^{-1}$, where $B$ is as stated in the lemma, we notice that $T_{\bb'}^{-1} = T_{\bb'^{-1}}=T_B$, where $T_{\bb'}:S_4\rightarrow S_4$ is the transformation such that $T_{\bb'}(\boldsymbol{\underline{s}}) = \bb'\boldsymbol{\underline{s}}/\|\bb'\boldsymbol{\underline{s}}\|_e$. Performing a change of variable in \eqref{eq:third_term} using $T_{\bb'}$, we get
	\begin{align}
		\int_{S_2}|u_{1,0}s_1+(u_{2,0}+\rho_2u_{2,1})s_2|^\alpha\Gamma_2(d\bs) & = \int_{S_4}|\scal{\boldsymbol{\underline{u}},\boldsymbol{\underline{s}}}|^\alpha\|\bb'^{-1}\boldsymbol{\underline{s}}\|^{-\alpha}\tilde{\Gamma}_{4,1}\circ T_{\bb'^{-1}} (d\boldsymbol{\underline{s}})\nonumber\\
		& = \int_{S_4}|\scal{\boldsymbol{\underline{u}},\boldsymbol{\underline{s}}}|^\alpha\Gamma_{4,1}(d\boldsymbol{\underline{s}}).\label{eq:gamma41}
	\end{align}
	Similarly for the fourth term in \eqref{eq:pre_charac_fun}, we have
	\begin{align*}
		\int_{S_2}|(\rho_1u_{1,0}+u_{1,1})s_1+u_{2,1}s_2|^\alpha\Gamma_2(d\bs) & = \int_{S_4}|u_{1,0}s_1+u_{2,0}s_2+(\rho_1u_{1,0}+u_{1,1})s_3+u_{2,1}s_4|^\alpha\tilde{\Gamma}_{4,2}(d\bs)\\
		& = \int_{S_4}|\scal{\bc\boldsymbol{\underline{u}},\boldsymbol{\underline{s}}}|^\alpha\tilde{\Gamma}_{4,2}(d\boldsymbol{\underline{s}})\\
		& = \int_{S_4}|\scal{\boldsymbol{\underline{u}},\bc'\boldsymbol{\underline{s}}}|^\alpha\tilde{\Gamma}_{4,2}(d\boldsymbol{\underline{s}}),
	\end{align*}
	with $\tilde{\Gamma}_{4,2}(\,\cdot\,)=\Gamma_2 \circ h_2 (S_{4,2}\cap \,\cdot\,)$, where $h_2:S_4\rightarrow S_2$ is the function defined by $h_2((s_1,s_2,s_3,s_4))=(s_3,s_4)$, $S_{4,2}=\{(0,0,s_3,s_4)\in S_4:\hspace{0.2cm} (s_3,s_4)\in S_2\}$, and
	$$
	\bc = 
	\begin{pmatrix}
	1 & 0 & 0 & 0\\
	0 & 1 & 0 & 0\\
	\rho_1 & 0 & 1 & 0\\
	0 & 0 & 0 & 1
	\end{pmatrix}.
	$$
	With a change of variable using $T_{\bc'}$, and since $\bc'^{-1}=C$, 
	\begin{align}
		\int_{S_2}|(\rho_1u_{1,0}+u_{1,1})s_1+u_{2,1}s_2|^\alpha\Gamma_2(d\bs) & = \int_{S_4}|\scal{\boldsymbol{\underline{u}},\boldsymbol{\underline{s}}}|^\alpha\|\bc'^{-1}\boldsymbol{\underline{s}}\|^{-\alpha}\tilde{\Gamma}_{4,2}\circ T_{\bc'^{-1}} (d\boldsymbol{\underline{s}})\nonumber\\
		& = \int_{S_4}|\scal{\boldsymbol{\underline{u}},\boldsymbol{\underline{s}}}|^\alpha {\Gamma}_{4,2} (d\boldsymbol{\underline{s}}).\label{eq:gamma42}
	\end{align}
	Finally, by \eqref{eq:pre_charac_fun}, \eqref{eq:mu_meas}, \eqref{eq:gamma41} and \eqref{eq:gamma42}, we have that \eqref{dem:charactoshow} holds with $\Gamma_4=\Delta + {\Gamma}_{4,1} + {\Gamma}_{4,2}$.
	One can check that $\Gamma_4$ is indeed symmetric:  for any transformation $g$ among $\{T_B,T_C,h_1,h_2\}$, we have $g(-\bs)=-g(\bs)$ for any $\bs\in S_4$, and as $\Gamma_2$ is symmetric by assumption, it is easy to check that the measures ${\Gamma}_{4,1}$ and ${\Gamma}_{4,2}$ are also symmetric. The case of	$\Delta$ is obvious.
	\cqfd
	
	\paragraph{Return to the proof of Proposition \ref{prop:reprbility_multivar}}\hfill
	
	\noindent By Lemma \ref{lem:xunderbar_spec} and Proposition \ref{prop238}, we know that $\boldsymbol{\underline{X}}_t$ will be representable on $C_4^{\|\cdot\|}$ if and only if 
	$$
	\Gamma_4(K^{\|\cdot\|}) = 0,
	$$
	where $K^{\|\cdot\|} = \{\boldsymbol{\underline{s}}\in {S}^4: \hspace*{0.1cm} \|\boldsymbol{\underline{s}}\|=0\}=\{\boldsymbol{\underline{s}}\in {S}^4: \hspace*{0.1cm} \underline{s}_1=\underline{s}_2=0 \}=S_{4,2}$. We have
	\begin{align*}
		\Gamma_4(K^{\|\cdot\|}) = \Delta(K^{\|\cdot\|}) + \Gamma_{4,1}(K^{\|\cdot\|}) + \Gamma_{4,2}(K^{\|\cdot\|}),
	\end{align*}
	with $\Delta$, $\Gamma_{4,1}$ and $\Gamma_{4,2}$ are as in Lemma \ref{lem:xunderbar_spec}. Given the points charged by $\Delta$, it is easily seen that $\Delta(K^{\|\cdot\|})=0$. Turning to $\Gamma_{4,1}(K^{\|\cdot\|})$, we have
	\begin{align*}
		\Gamma_{4,1}(K^{\|\cdot\|}) = \int_{K^{\|\cdot\|}}\|B\bs\|_e^{-\alpha}\tilde{\Gamma}_{4,1}\circ T_B(d\bs) = \tilde{\Gamma}_{4,1}\circ T_B(K^{\|\cdot\|}).
	\end{align*}
	Given that $T_B(K^{\|\cdot\|})=K^{\|\cdot\|}$ and $S_{4,1}\cap K^{\|\cdot\|} =\emptyset $, we have $\tilde{\Gamma}_{4,1}=\Gamma_2\circ h_1(S_{4,1}\cap K^{\|\cdot\|})=0$. Last, we have
	\begin{align*}
		\Gamma_{4,2}(K^{\|\cdot\|}) = \int_{K^{\|\cdot\|}}\|C\bs\|_e^{-\alpha}\tilde{\Gamma}_{4,2}\circ T_C(d\bs) = \tilde{\Gamma}_{4,2}\circ T_C(K^{\|\cdot\|}),
	\end{align*}
	with
	\begin{align*}
		T_C(K^{\|\cdot\|}) & = \left\{\dfrac{(-\rho_1s_3,0,s_3,s_4)}{(1+\rho_1^2)s_3^2+s_4^2}: \hspace*{0.1cm} \forall (s_3,s_4)\in\mathbb{R}^2, \hspace*{0.1cm} s_3^2+s_4^2=1\right\}\\
		& = \left\{\dfrac{(-\rho_1s_3,0,s_3,s_4)}{(1+\rho_1^2)s_3^2+s_4^2}: \hspace*{0.1cm} \forall (s_3,s_4)\in\mathbb{R}^2, \hspace*{0.1cm} s_3^2+s_4^2=1 \hspace*{0.2cm} \text{and} \hspace*{0.2cm} s_3\ne0\right\} \cup \{(0,0,0,\pm1)\}\\
		& := K' \cup \{(0,0,0,-1),(0,0,0,+1)\}
	\end{align*}
	Given that $K'\cap S_{4,2}=\emptyset$, we have by the $\sigma$-additivity of $\tilde{\Gamma}_{4,2}$,
	\begin{align*}
		\tilde{\Gamma}_{4,2}\Big(T_C(K^{\|\cdot\|})\Big) & = \tilde{\Gamma}_{4,2}(K') + \tilde{\Gamma}_{4,2}\Big(\{(0,0,0,-1),(0,0,0,+1)\}\Big)\\
		& = \Gamma_2\circ h_2(S_{4,2}\cap K') + \Gamma_2 \circ h_2(\{(0,0,0,-1),(0,0,0,+1)\})\\
		& = \Gamma_2 (\{(0,-1),(0,+1)\})
	\end{align*}
	Hence, $\Gamma_4(K^{\|\cdot\|}) = \Gamma_2\Big(\{(0,-1),(0,+1)\}\Big)$ and the representability condition for $\boldsymbol{\underline{X}}_t$ follows. 
	
	\subsection{Proof of Proposition \ref{prop:biar1_asymp}}
	
	To prove Proposition \ref{prop:biar1_asymp}, will make use of Proposition \ref{prop:cond_tail}. We first provide the form of the spectral representation of $\boldsymbol{\underline{X}}_t$ on $C_4^{\|\cdot\|}$ in the next lemma.
	\begin{lem}
		\label{lem:support_specmeas}
		Under the assumptions of Proposition \eqref{prop:reprbility_multivar} and assuming in addition that \eqref{asu:zeropoles} holds, then the characteristic function of the random vector $\boldsymbol{\underline{X}}_t$ can be written as
		$$
		\mathbb{E}[e^{i\scal{\boldsymbol{\underline{u}},\boldsymbol{\underline{X}}_t}}] = \exp\left\{-\int_{C_4^{\|\cdot\|}}|\scal{\boldsymbol{\underline{u}},\boldsymbol{\underline{s}}}|^\alpha\Gamma_4^{\|\cdot\|}(d\boldsymbol{\underline{s}})\right\}, \hspace{1cm} \text{for all} \hspace{0.3cm} \boldsymbol{\underline{u}}\in\mathbb{R}^4,
		$$
		where
		\begin{equation}\label{lem:spec_biv_oncylinder}
			\Gamma_4^{\|\cdot\|} = \Delta^{\|\cdot\|} + \Gamma_{4,1}^{\|\cdot\|} + \Gamma_{4,2}^{\|\cdot\|}. 
		\end{equation}
		Here,
		\begin{align}\label{lem:spec_biv_oncylinder_delta}
			\Delta^{\|\cdot\|} & = \dfrac{\sigma_1^\alpha}{2}\dfrac{ |\rho_1|^{2\alpha}}{1-|\rho_1|^\alpha}(\delta_{\{\bx_1\}}+\delta_{\{-\bx_1\}}) + \dfrac{\sigma_2^\alpha}{2}\dfrac{ |\rho_2|^{\alpha}}{1-|\rho_2|^\alpha}(\delta_{\{\bx_2\}}+\delta_{\{-\bx_2\}}),
		\end{align}
		with points $\bx_1,\bx_2$ as in Lemma \ref{lem:xunderbar_spec}, and
		\begin{align}
			\Gamma_{4,1}^{\|\cdot\|}(d\bs) & = \|B\bs\|_e^{-\alpha}\tilde{\Gamma}_{4,1}\circ T_B \circ T_{\|\cdot\|}^{-1}(d\bs),\label{lem:g41_interm}\\
			\Gamma_{4,2}^{\|\cdot\|}(d\bs) & = \|C\bs\|_e^{-\alpha}\tilde{\Gamma}_{4,2}\circ T_C \circ T_{\|\cdot\|}^{-1}(d\bs).\label{lem:g42_interm}
		\end{align}
		Moreover, for any Borel set $A\subset C_4^{\|\cdot\|}$,
		\begin{align}
			\Gamma_{4,1}^{\|\cdot\|}(A) & = \Gamma_2 \circ h_1 \Big(T_B\circ T_{\|\cdot\|}^{-1}(A\cap C_{4,1}^{\|\cdot\|})\Big), \label{lem:form_g1}\\
			\Gamma_{4,2}^{\|\cdot\|}(A) & = |\rho_1|^\alpha\int_{   T_C\circ T_{\|\cdot\|}^{-1}(A\cap C_{4,2}^{\|\cdot\|})}|s_3|^\alpha \Gamma_2 \circ h_2(d\bs),\label{lem:form_g2}
		\end{align}
		where
		\begin{align}
			C_{4,1}^{\|\cdot\|} & = \{(s_1,s_2,0,\rho_2s_2)\in C_4^{\|\cdot\|}: \hspace*{0.2cm} (s_1,s_2)\in S_2\},\label{lem:c41}\\
			C_{4,2}^{\|\cdot\|} & = \{\vartheta'(1,0,\rho_1^{-1},s_4)\in C_4^{\|\cdot\|}: \hspace*{0.2cm} s_4\in \mathbb{R}, \hspace*{0.1cm} \vartheta'\in\{-1,+1\}\}\label{lem:c42}
		\end{align}
	\end{lem}

	\noindent \textit{Proof.}\\
	Starting from $\Gamma_4$ as given in Lemma \ref{lem:xunderbar_spec} and applying a change of variable using $T_{\|\cdot\|}$ yields \eqref{lem:spec_biv_oncylinder}-\eqref{lem:g42_interm}.
	To show \eqref{lem:form_g1}, consider
	\begin{align*}
		\Gamma_{4,1}^{\|\cdot\|} (A) & = \int_{A} \|B \bs\|_e^{-\alpha} \tilde{\Gamma}_{4,1} \circ T_B \circ T_{\|\cdot\|}^{-1} (d\bs),
	\end{align*}
	and perform the first change of variable $\bs' = T_{\|\cdot\|}^{-1}(\bs) = \bs/\|\bs\|_e$ and get
	\begin{align*}
		\Gamma_{4,1}^{\|\cdot\|} (A) & = \int_{T_{\|\cdot\|}^{-1}(A)} \|B \bs'\|_e^{-\alpha} \|\bs'\|^{\alpha}\tilde{\Gamma}_{4,1} \circ T_B (d\bs').
	\end{align*}
	With the second change of variable $\bs = T_B(\bs') = B\bs'/\|B\bs'\|_e$, we obtain 
	\begin{align*}
		\Gamma_{4,1}^{\|\cdot\|} (A) & = \int_{T_B\circ T_{\|\cdot\|}^{-1}(A)} \|B^{-1} \bs\|^{-\alpha}\tilde{\Gamma}_{4,1} (d\bs).
	\end{align*}
	Given that $\tilde{\Gamma}_{4,1} (\,\cdot\,) = \Gamma_2 \circ h_1 (\,\cdot\,\cap S_{4,1})$, 
	\begin{align*}
		\Gamma_{4,1}^{\|\cdot\|} (A) & = \int_{T_B\circ T_{\|\cdot\|}^{-1}(A)\cap S_{4,1}} \|B^{-1} \bs\|^{-\alpha}\tilde{\Gamma}_{4,1} (d\bs),
	\end{align*}
	and noticing that for any $\bs\in S_{4,1}$, $B^{-1}\bs=\bs$ and $\|\bs\|=1$, we get
	\begin{align*}
		\Gamma_{4,1}^{\|\cdot\|} (A) & = \Gamma_2 \circ h_1 \Big(T_B\circ T_{\|\cdot\|}^{-1}(A)\cap S_{4,1}\Big),
	\end{align*}
	and using the fact that $T_B\circ T_{\|\cdot\|}^{-1}$ is bijective, we have 
	\begin{align*}
		T_B\circ T_{\|\cdot\|}^{-1}(A)\cap S_{4,1} & = T_B\circ T_{\|\cdot\|}^{-1}\Big(A\cap \hspace{0.2cm}T_{\|\cdot\|}\circ T_B^{-1} (S_{4,1}) \Big) = T_B\circ T_{\|\cdot\|}^{-1}(A\cap C_{4,1}^{\|\cdot\|} ),
	\end{align*}
	and \eqref{lem:form_g1} follows. We proceed similarly for \eqref{lem:form_g2} using in addition the fact that $\Gamma_2\Big(\{(0,-1),(0,+1)\}\Big)=0$.
	\cqfd
	
	\paragraph{Return to the proof of Proposition \ref{prop:biar1_asymp}}\hfill
	
	\noindent $(\iota)$ As \eqref{asu:zeropoles} holds, we know by Proposition \ref{prop:reprbility_multivar} that $\boldsymbol{\underline{X}}_t$ is representable on $C_4^{\|\cdot\|}$. By Lemma \ref{lem:support_specmeas}, we further know that its spectral measure $\Gamma_4^{\|\cdot\|}$ satisfies \eqref{lem:spec_biv_oncylinder}-\eqref{lem:spec_biv_oncylinder_delta} and \eqref{lem:form_g1}-\eqref{lem:form_g2}.
	Thus, by Proposition \ref{prop:cond_tail},
	\begin{align*}
		\mathbb{P}_x^{\|\cdot\|}\Big(\bXt,A_{\theta,\eta,P}\Big|B(V_0)\Big)\underset{x\rightarrow+\infty}{\longrightarrow} \dfrac{\Gamma_4^{\|\cdot\|}\Big(A_{\theta,\eta,P}\cap B(V_0) \Big)}{\Gamma_4^{\|\cdot\|}\Big(B(V_0)\Big)},
	\end{align*}
	and
	\begin{align}\label{dem:numoverdenom}
		\dfrac{\Gamma_4^{\|\cdot\|}\Big(A_{\theta,\eta,P}\cap B(V_0) \Big)}{\Gamma_4^{\|\cdot\|}\Big(B(V_0)\Big)} & = \dfrac{\Big[\Delta^{\|\cdot\|} + \Gamma_{4,1}^{\|\cdot\|} + \Gamma_{4,2}^{\|\cdot\|}\Big]\Big(A_{\theta,\eta,P}\cap B(V_0) \Big)}{\Big[\Delta^{\|\cdot\|} + \Gamma_{4,1}^{\|\cdot\|} + \Gamma_{4,2}^{\|\cdot\|}\Big]\Big(B(V_0)\Big)}.
	\end{align}
	From \eqref{lem:spec_biv_oncylinder_delta}, \eqref{lem:form_g1}-\eqref{lem:form_g2}, we can see that for any Borel set $A\subset C_4^{\|\cdot\|}$,
	\begin{align*}
		\Delta^{\|\cdot\|}(A) & = \Delta^{\|\cdot\|}(A\cap\{\pm \bx_1,\pm \bx_2\}),\\
		\Gamma_{4,1}^{\|\cdot\|} (A) & = \Gamma_{4,1}^{\|\cdot\|} (A\cap C_{4,1}^{\|\cdot\|})\\
		\Gamma_{4,2}^{\|\cdot\|} (A) & = \Gamma_{4,2}^{\|\cdot\|} (A\cap C_{4,2}^{\|\cdot\|}),
	\end{align*}
	where $C_{4,1}^{\|\cdot\|}$, $C_{4,2}^{\|\cdot\|}$ are given in \eqref{lem:c41} and \eqref{lem:c42}.
	We thus proceed in three steps: $(1)$ we derive the form of sets in the right-hand side of the above equations in the case $A=B(V_0)$, $(2)$ we then consider $A=A_{\theta,\eta,P}\cap B(V_0)$, $(3)$ we finally evaluate the mass over the obtained sets to derive the numerator and denominator in \eqref{dem:numoverdenom}.\\
	Let us consider the denominator.
	Because we assume $V_0\cap \{(\pm1,0),(0,\pm1)\}=\emptyset$, it is easy to see that
	\begin{align*}
		B(V_0)\cap \{\pm \bx_1,\pm \bx_2\} & = \emptyset,\\
		B(V_0)\cap C_{4,1}^{\|\cdot\|} & = \{(\cos u,\sin u, 0, \rho_2 \sin u): \hspace{0.2cm} u \in [\theta_0-\eta_0,\theta_0+\eta_0]\},\\
		B(V_0)\cap C_{4,2}^{\|\cdot\|} & = \emptyset.
	\end{align*}
	Thus, by \eqref{lem:form_g1},
	\begin{align*}
		\Gamma_4^{\|\cdot\|}\Big(B(V_0)\Big) & = \Gamma_{4,1}^{\|\cdot\|}\Big(B(V_0)\cap C_{4,1}^{\|\cdot\|}\Big)\\
		& = \Gamma_2 \circ h_1 \Big(T_B \circ T_{\|\cdot\|}^{-1}((B(V_0)\cap C_{4,1}^{\|\cdot\|})\Big)\\
		& = \Gamma_2 \circ h_1 \Big(\{(\cos u,\sin u, 0, 0): \hspace{0.2cm} u \in [\theta_0-\eta_0,\theta_0+\eta_0]\}\Big)\\
		& = \Gamma_2 (V_0).
	\end{align*}
	For the numerator, we have
	\begin{align*}
		A_{\theta,\eta,P}\cap B(V_0)\cap C_{4,1}^{\|\cdot\|} = 
		\left\{
		\begin{array}{cc}
			V_{\theta,\eta}\cap V_0, & \text{if} \hspace{0.3cm } (0,0) \in P,\\
			\emptyset, & \text{if} \hspace{0.3cm } (0,0) \notin P.
		\end{array}\right.
	\end{align*}
	The conclusion follows.
	
	\noindent $(\iota\iota)$ We proceed as in point $(\iota)$. Given the assumptions on $V_0$, we have
	\begin{align*}
		B(V_0)\cap \{\pm \bx_1,\pm \bx_2\} & = \{\vartheta \bx_2\},\\
		B(V_0)\cap C_{4,1}^{\|\cdot\|} & = \{(\cos u,\sin u, 0, \rho_2 \sin u): \hspace{0.2cm} u \in [\theta_0-\eta_0,\theta_0+\eta_0]\},\\
		B(V_0)\cap C_{4,2}^{\|\cdot\|} & = \emptyset.
	\end{align*}
	Thus, 
	\begin{align*}
		\Gamma_4^{\|\cdot\|}\Big(B(V_0)\Big) & = \Delta^{\|\cdot\|}(\{\vartheta \bx_2\}) + \Gamma_{4,1}^{\|\cdot\|}\Big(B(V_0)\cap C_{4,1}^{\|\cdot\|}\Big)\\
		& = \dfrac{\sigma_2^\alpha}{2}\dfrac{|\rho_2|^{\alpha}}{1-|\rho_2|^{\alpha}} + \Gamma_2 (V_0).
	\end{align*}
	Turning to $\Gamma_4^{\|\cdot\|}\Big(A_{\theta,\eta,P}\cap B(V_0)\Big)$, consider
	\begin{align*}
		& A_{\theta,\eta,P}\cap B(V_0)\cap \{\pm x_1,\pm x_2\}  = A_{\theta,\eta,P}\cap\{\vartheta \bx_2\}\\
		& \hspace{0.5cm} = \{\vartheta (0,1,0,\rho_2)\}\cap\left\{ (\cos u, \sin u, 0, \rho_2 \sin u) + (0,0,x,y): \hspace{0.2cm} u \in [\theta-\eta,\theta+\eta] \hspace{0.1cm} \text{and} \hspace{0.1cm} (x,y)\in P\right\}.
	\end{align*}
	Noticing that for any $u$ such that $(\cos u,\sin u)\ne(0,\vartheta)$, necessarily 
	$$
	\vartheta (0,1,0,\rho_2)\ne(\cos u, \sin u, 0, \rho_2 \sin u) + (0,0,x,y),\hspace{0.3cm} \text{for all} \hspace{0.3cm} (x,y)\in P,
	$$
	we have
	\begin{align*}
		\{\vartheta \bx_2\} \cap A_{\theta,\eta,P}  & =
		\left\{
		\begin{array}{cc}
			\{\vartheta \bx_2\} \cap \left\{ \vartheta\bx_2+(0,0,x,y): \hspace{0.1cm} (x,y)\in P\right\}, & \begin{array}{c}
				\text{if} \hspace{0.3cm } (\cos u, \sin u) = (0,\vartheta),  \\
				\text{for some} \hspace{0.1cm} u \in [\theta-\eta,\theta+\eta],
			\end{array}   \\
			& \\
			\emptyset, & \text{otherwise}. 
		\end{array}\right. 
	\end{align*}
	Hence
	\begin{align*}
		A_{\theta,\eta,P}\cap B(V_0)\cap \{\pm x_1,\pm x_2\} & =
		\left\{
		\begin{array}{cc}
			\{\vartheta x_2\}, & \text{if} \hspace{0.3cm } (0,\vartheta)\in V_{\theta,\eta} \hspace{0.2cm } \text{and} \hspace{0.2cm } (0,0)\in P \\
			\emptyset, & \text{otherwise}.
		\end{array}\right.
	\end{align*}
	Similarly, we have
	\begin{align*}
		A_{\theta,\eta,P}\cap B(V_0)\cap C_{4,1}^{\|\cdot\|} & = \left\{
		\begin{array}{cc}
			V_{\theta,\eta}\cap V_0, & \text{if} \hspace{0.3cm } (0,0) \in P,\\
			\emptyset, & \text{if} \hspace{0.3cm } (0,0) \notin P,
		\end{array}\right.\\
		A_{\theta,\eta}\cap B(V_0)\cap C_{4,2}^{\|\cdot\|} & = \emptyset.
	\end{align*}
	The result follows by evaluating $\Gamma_4^{\|\cdot\|}$ on the above sets.
	
	\noindent $(\iota\iota\iota)$ Proceeding as above, we first have that
	\begin{align*}
		B(V_0)\cap \{\pm x_1,\pm x_2\} & = \{\vartheta \bx_1\},\\
		B(V_0)\cap C_{4,1}^{\|\cdot\|} & = \{(\cos u,\sin u, 0, \rho_2 \sin u): \hspace{0.2cm} u \in [\theta_0-\eta_0,\theta_0+\eta_0]\},\\
		B(V_0)\cap C_{4,2}^{\|\cdot\|} & = \{\vartheta(1,0,\rho_1^{-1},s_4): \hspace{0.2cm} s_4\in\mathbb{R} \}.
	\end{align*}
	Hence,
	\begin{align*}
		\Gamma_4^{\|\cdot\|}(B(V_0)) & = \dfrac{\sigma_1^\alpha}{2}\dfrac{|\rho_1|^{2\alpha}}{1-|\rho_1|^\alpha} + \Gamma_{2}(V_0) + \Gamma_{4,2}^{\|\cdot\|}(\{\vartheta(1,0,\rho_1^{-1},s_4): \hspace{0.2cm} s_4\in\mathbb{R} \}).
	\end{align*}
	Given that
	$$
	T_C\circ T_{\|\cdot\|}^{-1}\Big(B(V_0)\cap C_{4,2}^{\|\cdot\|}\Big) = \left\{\vartheta\dfrac{(0,0,\rho_1^{-1},s_4)}{\sqrt{\rho_1^{-2}+s_4^2}}:\hspace{0.2cm} s_4\in\mathbb{R}\right\} = \{(0,0,s_1,s_2):\hspace{0.2cm} (s_1,s_2)\in S_2, \hspace{0.2cm} s_1\vartheta \rho_1>0\},
	$$
	the third term in $\Gamma_{4}^{\|\cdot\|}$ can be rewritten using \eqref{lem:form_g2} as
	\begin{align*}
		\Gamma_{4,2}^{\|\cdot\|}(\{\vartheta(1,0,\rho_1^{-1},s_4): \hspace{0.2cm} s_4\in\mathbb{R} \}) = |\rho_1|^{\alpha}\int_{\{(s_1,s_2)\in S_2:\hspace{0.2cm}s_1\vartheta \rho_1>0\}} |s_1|^\alpha \Gamma_2(d\bs).
	\end{align*}
	Since, by assumption, $\Gamma_2$ is symmetric and does not charge masses at $(0,\pm1)$, 
	$$
	\int_{\{(s_1,s_2)\in S_2:\hspace{0.2cm}s_1\vartheta \rho_1>0\}} |s_1|^\alpha \Gamma_2(d\bs)=\frac{1}{2}\int_{S_2} |s_1|^\alpha \Gamma_2(d\bs)=\sigma_1^\alpha/2,
	$$
	and thus
	\begin{align*}
		\Gamma_4^{\|\cdot\|}(B(V_0)) & = \dfrac{\sigma_1^\alpha}{2}\dfrac{|\rho_1|^{2\alpha}}{1-|\rho_1|^\alpha} + \Gamma_2(V_0) + |\rho_1|^\alpha\dfrac{\sigma_1^\alpha}{2}\\
		& = \Gamma_2(V_0) + \dfrac{\sigma_1^\alpha}{2}\dfrac{|\rho_1|^\alpha}{1-|\rho_1|^\alpha}.
	\end{align*}
	We last turn to $\Gamma_4\Big(A_{\theta,\eta,P}\cap B(V_0)\Big)$. We have
	\begin{align*}
		& A_{\theta,\eta,P}\cap B(V_0)\cap \{\pm x_1,\pm x_2\}  = A_{\theta,\eta,P}\cap\{\vartheta \bx_1\}\\
		& \hspace{0.5cm} = \{\vartheta (1,0,\rho_1^{-1},0)\}\cap\left\{ (\cos u, \sin u, 0, \rho_2 \sin u) + (0,0,x,y): \hspace{0.2cm} u \in [\theta-\eta,\theta+\eta] \hspace{0.1cm} \text{and} \hspace{0.1cm} (x,y)\in P\right\}\\
		& \hspace{0.5cm} = \left\{
		\begin{array}{cc}
			\{\vartheta \bx_1\} \cap \left\{ \vartheta(1,0,0,0)+(0,0,x,y): \hspace{0.1cm} (x,y)\in P\right\}, & \text{if} \hspace{0.3cm } (\vartheta,0)\in V_{\theta,\eta}  \\
			\emptyset, & \text{otherwise}. 
		\end{array}\right. \\
		& \hspace{0.5cm} = \left\{
		\begin{array}{cc}
			\{\vartheta x_1\}, & \text{if} \hspace{0.3cm } (\vartheta,0)\in V_{\theta,\eta} \hspace{0.2cm } \text{and} \hspace{0.2cm } (\vartheta \rho_1^{-1},0)\in P, \\
			\emptyset, & \text{otherwise}.
		\end{array}\right.
	\end{align*}
	As previously, we also have that
	$$
	A_{\theta,\eta,P}\cap B(V_0)\cap C_{4,1}^{\|\cdot\|}  = \left\{
	\begin{array}{cc}
	V_{\theta,\eta}\cap V_0, & \text{if} \hspace{0.3cm } (0,0) \in P,\\
	\emptyset, & \text{if} \hspace{0.3cm } (0,0) \notin P,
	\end{array}\right.
	$$
	Finally, 
	\begin{align*}
		& A_{\theta,\eta,P}\cap B(V_0)\cap C_{4,2}^{\|\cdot\|}  = A_{\theta,\eta,P}\cap\{\vartheta(1,0,\rho_1^{-1},s_4): \hspace{0.2cm} s_4\in\mathbb{R} \}\\
		& \hspace{0.5cm} = \left\{
		\begin{array}{cc}
			\{\vartheta(1,0,\rho_1^{-1},s_4): \hspace{0.2cm} s_4\in\mathbb{R} \} \cap \left\{ \vartheta(1,0,0,0)+(0,0,x,y): \hspace{0.1cm} (x,y)\in P\right\}, & \text{if} \hspace{0.3cm } (\vartheta,0)\in V_{\theta,\eta}  \\
			\emptyset, & \text{otherwise}. 
		\end{array}\right. \\
		& \hspace{0.5cm} = \left\{
		\begin{array}{cc}
			\{\vartheta(1,0,\rho_1^{-1},y): \hspace{0.2cm} y\in P_2 \}, & \text{if} \hspace{0.3cm } (\vartheta,0)\in V_{\theta,\eta} \hspace{0.2cm } \text{and} \hspace{0.2cm } \vartheta \rho_1^{-1} \in P_1, \\
			\emptyset, & \text{otherwise},
		\end{array}\right. 
	\end{align*}
	and 
	$$
	T_C\circ T_{\|\cdot\|}^{-1}\Big(\{\vartheta(1,0,\rho_1^{-1},y): \hspace{0.2cm} y\in P_2 \}\Big) = \left\{\vartheta\dfrac{(0,0,\rho_1^{-1},y)}{\sqrt{\rho_1^{-2}+y^2}}:\hspace{0.2cm} y\in P_2\right\}.
	$$
	\cqfd
	
	\clearpage

	\begin{center}
		{\sc References}
	\end{center}
	
	

\end{document}